\documentclass[12pt,oneside,reqno]{article}
\usepackage{amsmath}
\usepackage{paralist}
\usepackage{graphics}
\usepackage{epsfig}
\usepackage{float}
\usepackage{graphicx}
\usepackage{lineno}
\usepackage{epstopdf}
\usepackage{amssymb}
\usepackage{cite}
\usepackage{mathrsfs}
\usepackage{amsthm}
\usepackage{xcolor}
\usepackage{caption}
\usepackage{subfigure}
\usepackage{txfonts}
\usepackage[colorlinks=true]{hyperref}
\hypersetup{urlcolor=blue, citecolor=blue}
\usepackage[utf8]{inputenc}
\usepackage[shortlabels]{enumitem}
\usepackage{cancel} 
\usepackage[margin=1.0in]{geometry}
\usepackage{overpic}
\usepackage[normalem]{ulem}

\numberwithin{equation}{section}

\newtheorem{theorem}{Theorem}[section]
\newtheorem{corollary}{Corollary}
\newtheorem{lemma}[theorem]{Lemma}
\newtheorem{proposition}{Proposition}

\newtheorem{physical conclusion}{Physical Conclusion}

\newtheorem{definition}[theorem]{Definition}
\newtheorem{remark}{Remark}

\definecolor{darkgreen}{rgb}{0,0.35,0}

\title{Global existence  and Rayleigh-Taylor instability
for the semi-dissipative Boussinesq system with Naiver boundary conditions}

\author{
Huafei Di$^{a}$
\thanks{dihuafei@gzhu.edu.cn}
\quad\quad
Liang Li$^{a}$
\thanks{Corresponding author:llbohou@gzhu.edu.cn }
\quad\quad
Xiaoming Peng$^{b}$
\thanks{pengxm@gdufe.edu.cn}
\quad\quad
Quan Wang$^{c}$
\thanks{xihujunzi@scu.edu.cn}
\\ \footnotesize $^a$ School of Mathematics and Information Science,
\footnotesize Guangzhou University,
\\\footnotesize Guangzhou,
 Guangdong, 510000, China
 \\
 \footnotesize $^{b}$ School of Statistics and Mathematics,
\footnotesize Guangdong University of Finance and Economics,
\\
  \footnotesize Guangzhou, Guangdong, 510320, China
  \\ \footnotesize $^{c}$ College of Mathematics, Sichuan University,
  \footnotesize
 Chengdu, Sichuan, 610065,  China
}
\medskip

\begin{document}
\maketitle
\begin{abstract}
    Considered herein is the global existence of weak, strong solutions and Rayleigh-Taylor (RT) instability for 2D semi-dissipative
    Boussinesq equations in an infinite strip domain $\Omega_{\infty}$ subject to Navier boundary conditions with non-positive slip coefficients.   We first prove the global existence of weak and strong solutions on bounded domain $\Omega_{R}$ via the Galerkin method,  characteristic analyzing technique and Stokes estimates etc. Based on above results, we further derive the uniform estimates, independent of the length of horizontal direction of $\Omega_{R}$, ensuring the global existence of weak and strong solutions in  unbounded case $\Omega_{\infty}$ by utilizing the domain expansion method. Moreover,
     when the steady temperature is higher with decreasing
height (i.e., RT steady-state) on certain region,
we demonstrate that the steady-state is linear unstable through the construction of energy functional and the settlement of a family of modified variational problems. Furthermore, with the help of unstable solutions constructed in linear instability and global existence theorems, we confirm the instability of nonlinear problem in a Lipschitz structural sense. Finally, we give a series of rigorous verification (see Appendix) including the spectra of Stokes equations with Navier boundary conditions, Sobolev embedding inequalities, trace inequalities, and Stokes estimates under Navier boundary conditions etc, used in the proof of main conclusions.
    \end{abstract}
    \newpage
    \tableofcontents
    \newpage
\section{Introduction}
This paper is interested in the following semi-dissipative Boussinesq system
\begin{align}\label{Model}
        \begin{aligned}
            \begin{cases}
         \frac{\partial \mathbf{u}}{\partial t}+(\mathbf{u}\cdot \nabla)\mathbf{u}+\nabla p
         =\mu\Delta\mathbf{u}+\theta
         \mathbf{g},
         \\
         \frac{\partial\theta}{\partial t}+(\mathbf{u}\cdot\nabla)\theta=0,
         \\
         \nabla\cdot\mathbf{u}=0,
            \end{cases}
        \end{aligned}
    \end{align}
    on $\Omega\times[0,T]$, where $\Omega$ is a domain of $R^{N}(N\geq 2)$ and $T$ is the maximum existence time, the unknown functions
    $\mathbf{u},\theta$ and $p$ are the velocity field, temperature field and pressure field, respectively,
    $\mu>0$ is the fluid viscosity. The symbol $\mathbf{g}$ represents the constant and upward-pointing gravity vector $(0,g)$, where $g$ is the (scalar) acceleration from gravity. The Boussinesq equations have been applied in the examination of
    buoyancy-driven flows, such as atmospheric
    fronts and ocean circulation\cite{Majda2003,Pedlosky2003} and then attracted wide attention
    from some researchers, see \cite{Adhikari2016,Biswas2017,He2022} and references therein.
    The fully dissipative Boussinesq system, characterized by positive constants for all viscous and
    diffusivity coefficients, is the global well-posedness\cite{Foias1987,Temam2009}.
    In contrast, the inviscid scenario (all viscosities and diffusivity are zero) poses an
    open question\cite{Majda2002}
    concerning the global well-posedness. As the intermediate cases between the
    fully dissipative and inviscid scenarios, the Boussinesq equations with partial dissipation, that is,
    with partial viscosities or partial diffusivity has obvious value in modeling the
    dynamics of geophysical flows in which the dissipation in
    specific directions dominates\cite{Majda1997,Chemin2000,Chemin2008}.
    The Boussinesq equations \eqref{Model} with viscosity and without heat diffusion,
    exhibiting dissipative characteristics in velocity but not in
    temperature, are categorized as semi-dissipative.

    As far as we know, the Boussinesq equation is mostly studied with periodic boundary condition\cite{Biswas2017,He2022}
     or Dirichlet
    condition\cite{Li2022} which means the fluid does not slip along the boundary. Here, the semi-dissipative Boussinesq system \eqref{Model} is equipped
    with the Naiver boundary conditions
\begin{align}\label{bianjie1230}
        \begin{aligned}
        &~~~~~~~~~~~~~~~~~~~~~~~\mathbf{u}\cdot\mathbf{n}=0,~\text{on~}\partial\Omega,
        \\
       &\left[\left(-p \mathbf{I}+\mu\left(\nabla\mathbf{u}+(\nabla\mathbf{u})^{\text{Tr}}\right)\right)\cdot\mathbf{n}\right]\cdot\tau=
        \alpha \mathbf{u}\cdot\tau,~\text{on~}\partial\Omega,
        \end{aligned}
\end{align}
proposed by the famous
    mathematician and physicist C. Naiver\cite{Naiver1827} in 1827 under taking the slip of the fluid on the boundary, such
    as hurricanes and tornadoes.
    Here $\mathbf{n}$ is the outward normal vector field to $\partial \Omega$, $\tau$ is the corresponding
tangential vector, Tr means matrix transposition, $\mathbf{I}$ is the $N\times N$ identity
matrix, $\alpha$ stands for physical meaning parameter which is either a constant or a function, even a
smooth matrix. Therefore,
$\mathbf{u}\cdot\mathbf{n}|_{\partial\Omega}$ means that the fluid will not penetrate the boundary, and the second condition signifies that the fluid slips along the boundary.
For more information about the work involving Naiver boundary conditions,
one can refer to \cite{Ding2017,Ding2020,Hu2018,Wang2021}.

In this paper, we consider the initial-boundary value problem for semi-dissipative Boussinesq system in a strip domain $\Omega_{\infty}=\mathbf{R}\times (0,1)$,
\begin{align}\label{model}
    \begin{aligned}
        \begin{cases}
     \frac{\partial \mathbf{u}}{\partial t}+(\mathbf{u}\cdot \nabla)\mathbf{u}+\nabla p
     =\mu\Delta\mathbf{u}+\theta
     \mathbf{g},
     \\
     \frac{\partial\theta}{\partial t}+(\mathbf{u}\cdot\nabla)\theta=0,
     \\
     \nabla\cdot\mathbf{u}=0,
     \\
     \mathbf{u}(x,0)=\mathbf{u}_{0},~\theta(x,0)=\theta_{0},
        \end{cases}
    \end{aligned}
\end{align}
with the Naiver boundary conditions
\begin{align}\label{bianjietiaojian1230}
    \begin{aligned}
       & ~~~~~~~~~~~~~~~~~~~~~~\mathbf{u}\cdot\mathbf{n}=0,~~\text{on}~\{x_{2}=0,1\},
        \\
        &\left[\left(-p \mathbf{I}+\mu\left(\nabla\mathbf{u}+(\nabla\mathbf{u})^{\text{Tr}}\right)\right)\cdot\mathbf{n}\right]\cdot\tau
        =k_{1}\mathbf{u}\cdot \tau,~\text{on}~\{x_{2}=1\},
        \\
        &\left[\left(-p \mathbf{I}+\mu\left(\nabla\mathbf{u}+(\nabla\mathbf{u})^{\text{Tr}}\right)\right)\cdot\mathbf{n}\right]\cdot\tau
        =k_{0}\mathbf{u}\cdot \tau,~\text{on}~\{x_{2}=0\},
    \end{aligned}
\end{align}
where $\mathbf{u}(x,t)=(u_{1}(x,t),u_{2}(x,t))$, $p=p(x,t)$ and $\theta=\theta(x,t)$,
$x=(x_{1},x_{2})$ is the spatial variable;
$\mathbf{I}$ is $2\times 2$ identity matrix, $\mathbf{n}=(0,1)$ on $\{x_{2}=1\}$ and
$\mathbf{n}=(0,-1)$ on $\{x_{2}=0\}$ are the outward unit normal vectors, $\tau=(1,0)$ on both $\{x_{2}=0,1\}$ is the tangential vector. Moreover, throughout this
paper, we assume that the slip coefficients $k_{0}$ and $k_{1}$ are always non-positive constants. In this case, the conditions \eqref{bianjietiaojian1230} can
be simplified  to forms
\begin{align}\label{naiverboundarycondition}
    \begin{aligned}
    &u_{2}(x_{1},0)=u_{2}(x_{1},1)=0,~\partial_{2} u_{1}(x_{1},1)=\frac{k_{1}}{\mu}u_{1}(x_{1},1),
    ~\partial_{2} u_{1}(x_{1},0)=-\frac{k_{0}}{\mu}u_{1}(x_{1},0),~x_{1}\in\mathbf{R},
    \end{aligned}
\end{align}
where $\partial_{i}u_{1}:=\frac{\partial u_{1}}{\partial x_{i}}$ for $i=1,2$. Furthermore, we can assume, without loss of generality, that
\begin{align}\label{boundarycondition}
\lim\limits_{|x_{1}|\rightarrow+\infty}\mathbf{u}(x_{1},x_{2})=\mathbf{0},
\end{align}
due to the Galilean invariance of fluid mechanics\cite{jun_choe_strong_2003}.

One of the purposes in the present paper is to investigate the global well-posedness of system \eqref{model} with
conditions \eqref{bianjietiaojian1230} and \eqref{boundarycondition}. It is worth mentioning that the unbounded nature of $\Omega_{\infty}$ causes some challenges in the study of its qualitative properties. For example, the
eigenvectors of the Stokes operator cannot constitute a basis. Additionally, even in bounded domain $\Omega_{R}$, it is difficult to perform a spectrum analysis of linear differential operator to \eqref{model} because of the presence of $\eqref{model}_{2}$.
 These difficulties hinder the direct application
of the Galerkin approximation method to problem \eqref{model}. Drawing inspiration from \cite{jun_choe_strong_2003},
we explore the case on a bounded domain $\Omega_{R}$, which approximates $\Omega_{\infty}$
as $R\rightarrow+\infty$ (see \autoref{quyutu1}). If the uniform estimates of the solution independent of $R$
are achieved, the global well-posedness on $\Omega_{\infty}$ can be effectively addressed. This approach is called the  domain expansion technique. Furthermore, motivated by \cite{kim_weak_1987}, we employ the semi-Galerkin
method to tackle problem \eqref{model} on $\Omega_{R}$. Specifically, we adopt the  Galerkin approximation method to $\eqref{model}_{1}$ and the characteristic analyzing technique to $\eqref{model}_{2}$.

It is known that the spectrum must be investigated when employing Galerkin approximation method. Therefore,
this paper establishes the spectrum analysis of the Stokes operator subjet to Navier boundary conditions.
These results  align with the conclusions under Dirichlet boundary conditions, more details can be seen
Lemma \ref{tezhenzhidewenti}. The spectrum of the Stokes operator with Navier boundary
conditions has been studied in previous works such as \cite{Clopeau1998,li_global_2021}.
However, we obtain a set of smooth basis vectors compared with \cite{Clopeau1998} (basis vectors just belong to $H^{3}$), and directly employ a classical variational method and bootstrap method for better understanding in contrast to \cite{li_global_2021} (an auxiliary eigenvalue problem was introduced).
The Stokes estimate plays a key role in establishing the existence of strong solutions. In \cite{Temam1977,Galdi2011},
the constant appearing in Stokes type inequality  depends on the bounded domain, which leads to a failure in direct application of Stokes
estimate for our case that $\Omega_{\infty}$ is unbounded in $x_{1}$-direction and bounded in $x_{2}-$direction. Hence, we establish the necessary Stokes estimate in
Lemma \ref{tishengzhengzexing} and Corollary \ref{tuilun1216}, which are independent of horizontal length. Particularly,
in the proof of Stokes estimate, we propose a
new perspective, that is, applying the difference method with the dual theory of uniform convex Banach space to achieve our goal, rather than utilizing harmonic analysis approach \cite{Temam1977,Galdi2011}. Note that our conclusion involving Stokes estimate is valid for any $q>1$,
not just for $q=2$ in \cite{Ding2017}. In addition, some classical Sobolev embedding inequalities and
trace theorem,
have been extensively
covered in references such as \cite{adams_sobolev_1975,noauthor_partial_nodate,Taylor2011}, where the constants in above inequalities are depending on domain. Therefore, to get the main results of this paper,
we verify the Sobolev embedding inequalities (see Lemmas \ref{sobolveqianrudingli12241}-\ref{jiushizheyang12242})
and a trace Theorem (see Lemma \ref{jidingli}), which are all independent of $x_1$-direction length $R$. All above prerequisites ensure the applicability of domain expansion technique.

In doing so, via the combination of semi-discrete Galerkin method, characteristic analyzing technique and
stokes estimates etc, we prove the global existence of weak and strong solutions on bounded domain $\Omega_{R}$,
cf. Theorems \ref{cunzairuojie}-\ref{kanqilaimaodun}.
Furthermore, after deriving the uniform estimates of the solution independent of the horizontal length of $\Omega_{R}$,
we obtain the global existence of weak and strong solutions in case $\Omega_{\infty}$ by the domain expansion technique, cf. Theorems \ref{cunzairuojie1}-\ref{qiangjie1}.

With the global existence of  strong solution in hand, the second goal of this paper is to study the RT instability of a smooth steady profile solution
$(\overline{\mathbf{u}},\overline{\theta},\overline{p})$ to
problem \eqref{model}-\eqref{boundarycondition}. We know that the RT instability driven by gravity occurs in a heavy fluid on top of a light one. Here,
$\overline{\theta}$ satisfies the conditions:
\begin{align}\label{cond1}
    \overline{\theta}=\overline{\theta}(x_{2})\in C^{\infty}([0,1])~
    \text{and~}D\overline{\theta}(x_{2}^0)<0,~~\text{for~some}~x_{2}^{0} \in (0,1),
\end{align}
where $D:=\frac{d}{dx_{2}}$. Clearly, this $\overline{\theta}$ with $\mathbf{\overline{u}}\equiv (0,0)$ forms
a steady state to \eqref{model}, provided
\begin{align}\label{cond20105}
    \overline{p}=\overline{p}(x_{2})~\text{and}~D \overline{p}=g\overline{\theta}.
\end{align}
\begin{remark}Note that the condition \eqref{cond1} shows that there is at least a part of region where the
    steady temperature solution has higher temperature with lower $x_{2}$ (height). It is well known
    that the higher temperature is, the smaller density is, which means
    in this region, the larger density fluid is on top of the
    smaller density fluid. As a result, this will lead to the classical
    RT instability. Our research will offer a new perspective on this instability phenomenon.
\end{remark}

Hence, we take the perturbation to be
\begin{align}
    \mathbf{u}=\mathbf{v}+\overline{\mathbf{u}},~~\theta=\Theta +\overline{\theta},~~
    p=\pi+\overline{p}.
\end{align}
Then, a direct calculation gives that $(\mathbf{v},\Theta,\pi)$ satisfies the perturbed equations
\begin{align}\label{henyoupinwei}
    \begin{cases}
        \partial_{t}\mathbf{v}+(\mathbf{v}\cdot\nabla)\mathbf{v}+\nabla \pi=\mu\Delta\mathbf{v}+
        \Theta \mathbf{g},
        \\
        \partial_{t}\Theta+(\mathbf{v}\cdot\nabla)\Theta+v_{2}D\overline{\theta}=0,
        \\
        \nabla\cdot\mathbf{v}=0.
    \end{cases}
\end{align}
To complete the statement of above perturbed problem, from the boundary condition \eqref{bianjietiaojian1230},
we can specify the boundary conditions and initial values:
\begin{align}\label{naiverboundarycondition123}
    \begin{aligned}
    &v_{2}(x_{1},0)=v_{2}(x_{1},1)=0,
    ~\partial_{2}v_{1}(x_{1},1)=\frac{k_{1}}{\mu}v_{1}(x_{1},1),
    \\
    &\partial_{{2}}v_{1}(x_{1},0)=-\frac{k_{0}}{\mu}v_{1}(x_{1},0),~x_{1}\in\mathbf{R},
    ~\lim\limits_{|x_{1}|\rightarrow\infty}\mathbf{v}(t,x)=\mathbf{0},
    \end{aligned}
\end{align}
and
\begin{align}\label{123456}
    (\mathbf{v},\Theta)|_{t=0}=(\mathbf{v}_{0},\Theta_{0}).
\end{align}
Moreover, the initial data should satisfy $\nabla\cdot\mathbf{v}_{0}=0$.

To analyze the linear instability of problem \eqref{henyoupinwei}-\eqref{123456}, we need to consider the linearized equations of \eqref{henyoupinwei}, which read as
\begin{align}\label{xianxingpart}
    \begin{aligned}
        \begin{cases}
        \partial_{t}\mathbf{v}+\nabla \pi=\mu\Delta\mathbf{v}+\Theta\mathbf{g},
        \\
        \partial_{t}\Theta+v_{2}D\overline{\theta}=0,
        \\
        \nabla\cdot\mathbf{v}=0.
        \end{cases}
    \end{aligned}
\end{align}

 For the RT instability,  Rayleigh\cite{Rayleigh1883} first introduced the instability of the linearized
problem (i.e. linear instability) for an incompressible fluid in 1883. The RT instability
in the Hadamard sense of 2D nonhomogeneous incompressible
inviscid fluid in strip domain with zero normal velocity on its boundary was explored by Hwang and
Guo in 2003\cite{Hwang2003}. However, for viscous fluid, note that it is difficult to construct exponentially growing solutions
for the linearized equation lacking a direct variational structure. To address this challenge, Guo and
Tice\cite{Guo2010} developed a general method in 2011, employing a modified variational problem and
fixed-point theory. In 2014, Jiang et al.\cite{Jiang2014} investigated 3D gravity-driven viscous flows
in a bounded domain, considering nonlinear RT instability in the Lipschitz structure sense.
More information on RT instability can be referred to \cite{Jiang2015,Jiang2022} and references therein.

In views of the works mentioned above, we remark that the existence of unstable solutions to the linearized and nonlinear problems
\eqref{henyoupinwei}-\eqref{123456} are all interesting and challenging, and it remains an open question up to now. Regarding the linear instability, the absence of a variational structure for $\lambda$ in \eqref{equ1} is
caused by the viscosity, and then this prompts us to draw inspiration from  \cite{Guo2010,Jiang2013}. Hence, we overcome this challenge by considering a family of modified problems with variational structures.
Afterwards, the growing solution to original problem is obtained by using the intermediate value theorem of continuous functions. Thus,  the linear instability of \eqref{henyoupinwei}-\eqref{123456} is completed, cf. Theorem \ref{xianxingbuwending}.
On the other hand, concerning nonlinear RT instability, we focus on initial values $(\mathbf{u}_{0}, \Theta_{0})$
that ensure the solvability of the linearized problem \eqref{naiverboundarycondition123}-\eqref{xianxingpart}.
We demonstrate that the solution to nonlinear problem \eqref{henyoupinwei}-\eqref{123456}, corresponding
to the same initial values $(\mathbf{u}_{0}, \Theta_{0})$, exhibits instability in the Lipschitz structure.
To be specific, this proof can be completed by employing the method of contradiction. Here, we construct a family of strong solutions to the
nonlinear problem \eqref{henyoupinwei}-\eqref{123456}. Leveraging the nonlinear estimate of strong
solutions (see Proposition \ref{feichangzhongyaode}), we establish that the limit of this family of
strong solutions satisfies the linearized equation \eqref{naiverboundarycondition123}-\eqref{xianxingpart}.
By exploiting the exponential growth rate of the solution to linearized problem, a contradiction is
obtained (see Lemma \ref{maodunchengli0105}). Thus, we conclude the result on nonlinear RT instability, cf. Theorem \ref{feixianxingbuwending1217}.

The plan of the rest of this paper is as follows: Section \ref{chubu0104} introduces some essential notations and lemmas including the function spaces and spectrum of the Stokes operator with Naiver boundary conditions.
In Section \ref{jiedecunzaixing0104}, we establish the global well-posedness to problem
\eqref{model}-\eqref{boundarycondition}. Sections \ref{RTbuwendingxing0104} and \ref{feixian0331} address the linear
and nonlinear instability of problem \eqref{henyoupinwei}-\eqref{123456}.
In the last Section \ref{peijian0104}, we provide some proofs of preliminary results in Section \ref{chubu0104}.

\section{Preliminary}\label{chubu0104}
In this section, we introduce
some notations of regions and function spaces, and list some
lemmas which will be proven in Appendix.

\subsection{Notations and function spaces}

\begin{figure}[H]
        \centering
        {\includegraphics[width=3in]{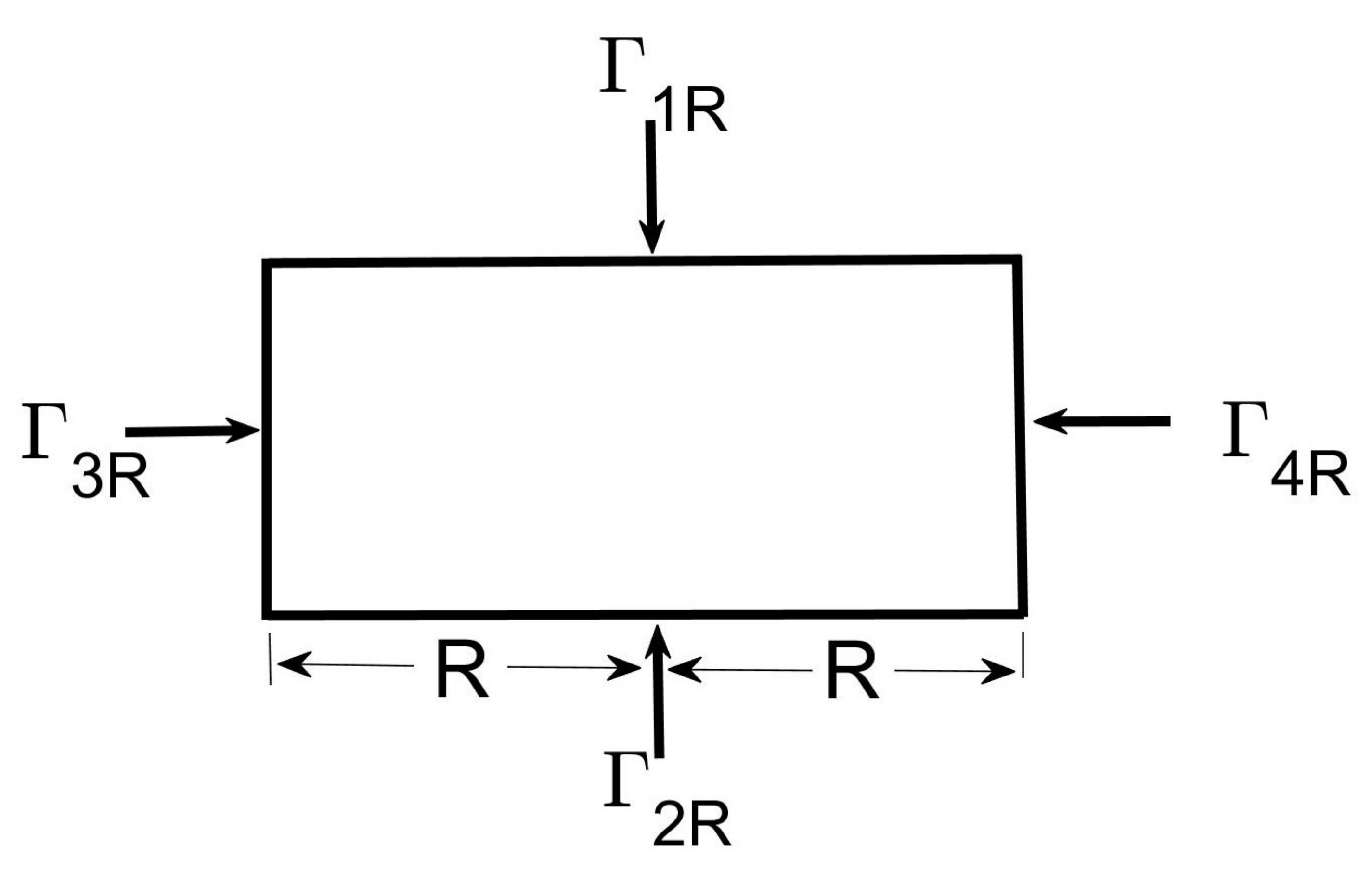}}
    \caption{$\Omega_{R}$}
        \label{quyutu1}
    \end{figure}

As mentioned in the introduction, the unbounded nature of $\Omega_{\infty}$ presents some challenges
in obtaining a smooth orthonormal basis for constructing Galerkin approximate solutions.
To address these difficulties, we utilize the domain expandsion technique. Thus,
we introduce bounded region $\Omega_{R}$ (see \autoref{quyutu1}),
which serves as an approximation to $\Omega_{\infty}$.
\noindent And we define
    \begin{align*}
        \begin{aligned}
           &\Gamma_{1 R}=\left\{(x_{1},x_{2})\in \mathbf{R}^2\big{|} x_{1}\in [-R,R]~\text{and}~x_{2}=1\right\},~
            \Gamma_{2 R}=\left\{(x_{1},x_{2})\in \mathbf{R}^2\big{|} x_{1}\in [-R,R]~\text{and}~x_{2}=0\right\},
            \\
            &\Gamma_{3 R}=\left\{(x_{1},x_{2})\in \mathbf{R}^2\big{|} x_{1}=-R~\text{and}~x_{2}\in [0,1]\right\},
            ~\Gamma_{4 R}=\left\{(x_{1},x_{2})\in \mathbf{R}^2\big{|} x_{1}=R~\text{and}~x_{2}\in [0,1]\right\}.
        \end{aligned}
        \end{align*}
    $\Omega_{R}$ is a rectangular region enclosed by $\Gamma_{j R}~(j=1,2,3,4)$.
  For convenience, let  $q\in(1,+\infty)$ and define the spaces as follows:
\begin{align}\label{kongjiandingyi}
    \begin{aligned}
        &\mathbf{E}^{1,q}(\Omega_{R})=\left\{ \mathbf{u}\in \left[W^{1,q}(\Omega_{R})\right]^2 \big{|}
        \mathbf{u}\cdot\mathbf{n}|_{\Gamma_{1R}\cup\Gamma_{2R}}=0,~
        \mathbf{u}|_{\Gamma_{3 R}\cup\Gamma_{4R}}=\mathbf{0}\right\},
        \\
        &
        \begin{aligned}
            \mathbf{G}^{1,q}(\Omega_{R})=\left\{\mathbf{u}\in \mathbf{E}^{1,q}(\Omega_{R})\big{|}
            \nabla\cdot\mathbf{u}=0\right\},~\widehat{\mathbf{L}}^{q}(\Omega_{R})=\left\{\mathbf{u}\in
            \left[L^{q}(\Omega_{R})\right]^2\big{|}\nabla\cdot\mathbf{u}=0\right\},
        \end{aligned}
        \\
        &\begin{aligned}
        \mathbf{V}^{1,q}(\Omega_{R})=
            \bigg{\{}\mathbf{u}\in &\left[W^{2,q}(\Omega_{R})\right]^2\cap \mathbf{G}^{1,q}(\Omega_{R})\big{|}
       \partial_{{2}}u_{1}(x_{1},1)=\frac{k_{1}}{\mu}u_{1}(x_{1},1),
        \\
        &~\partial_{{2}}u_{1}(x_{1},0)=-\frac{k_{0}}{\mu}u_{1}(x_{1},0),~x_{1}\in [-R,R]
        \bigg{\}},
        \end{aligned}
        \\
        &\mathbf{C}_{0,x_{1}}(\Omega_{\infty})=\bigg{\{}\mathbf{u}\in \big{[}C^{\infty}(\Omega_{\infty})\big{]}^{2}\big{|}
\exists ~r>0,
\text{such~that~}\mathbf{u}=\mathbf{0},~\forall~|x_{1}|\geq r,
\\
&~~~~~~~~~~~~~~~~~~~~~~~~~~\nabla\cdot
\mathbf{u}=0~\text{and}~\mathbf{u}\cdot\mathbf{n}|_{x_{2}=0,1}=0
\bigg{\}},
    \end{aligned}
\end{align}
where $\mathbf{n}$ is the unit outward normal vector,
$W^{1,q}(\Omega_{R})$, $W^{2,q}(\Omega_{R})$ and $L^{q}(\Omega_{R})$ are usual Sobolev and Lebesgue spaces,
$\|\cdot\|_{\mathbf{E}^{1,q}(\Omega_{R})}=\|\cdot\|_{\mathbf{G}^{1,q}(\Omega_{R})}=
\|\cdot\|_{W^{1,q}(\Omega_{R})}$,
 $\|\cdot\|_{\mathbf{V}^{1,q}(\Omega_{R})}=\|\cdot\|_{W^{2,q}(\Omega_{R})}$ and
 $\|\cdot\|_{\widehat{\mathbf{L}}^{q}(\Omega_{R})}=\|\cdot\|_{L^q(\Omega_{R})}$.

In particular, from Lemmas \ref{jidingli} and \ref{sobolveqianrudingli12241},
we can define the inner product and norm of $\mathbf{G}^{1,2}(\Omega_{R})$ as follows,
$\forall~ \mathbf{u}$ and $\mathbf{v}$ $\in \mathbf{G}^{1,2}(\Omega_{R})$,
\begin{align}\label{neiji}
    \begin{aligned}
((\mathbf{u},\mathbf{v}))=\mu\int_{\Omega_{R}}\nabla\mathbf{u}\cdot\nabla\mathbf{v}dx
        &-\int_{-R}^{R}\big{[}k_{1}u_{1}(x_{1},1)v_{1}(x_{1},1)+k_{0}u_{1}(x_{1},0)v_{1}(x_{1},0)\big{]}dx_{1},
        \\
       &\big{\|}\mathbf{u}\big{\|}_{\mathbf{G}^{1,2}(\Omega_{R})}^{2}=((\mathbf{u},\mathbf{u})).
    \end{aligned}
\end{align}
It is not difficult to verify that
$\mathbf{G}^{1,2}(\Omega_{R})$ is a Hilbert space and $\big{\|}\mathbf{u}\big{\|}_{\mathbf{G}^{1,2}(\Omega_{R})}$
is equivalent to $\big{\|}\nabla\mathbf{u}\big{\|}_{L^2(\Omega_{R})}$, see Lemma \ref{sobolveqianrudingli12241}.

\subsection{Lemmas and corollaries}
Generally, the uniform constants in certain inequalities may depend
on the shape or size of domain\cite{Galdi2011,noauthor_partial_nodate}, which poses a challenge in the application of domain expansion
technique. To overcome this challenge, we establish the following some lemmas and corollaries,
where the uniform constants are independent of the horizontal length of $\Omega_{R}$.
The proofs for these inequalities will be provided in the Appendix \ref{11226}-\ref{appendixb1226}.

\begin{lemma}\label{jidingli}Let $u\in W^{1,q}(\Omega_{R})$
 $(1<q<\infty)$. Then we have
    \begin{align}
        \begin{aligned}
        &\big{\|}u(x_{1},0)\big{\|}_{L^{q}([-R,R])}\leq C(q)
       \big{\|}u\big{\|}_{W^{1,q}(\Omega_{R})},~
        \\
        &\big{\|}u(x_{1},1)\big{\|}_{L^{q}([-R,R])}\leq C(q)
        \big{\|}u\big{\|}_{W^{1,q}(\Omega_{R})},
        \end{aligned}
    \end{align}
    where $C(q)$ is a constant only depending on $q$.
\end{lemma}
Thus, let $R\rightarrow+\infty$. And we have the following corollary.
\begin{corollary}\label{jidingli0125}Let $u\in W^{1,q}(\Omega_{\infty})$ $(1<q<\infty)$.
    Then we have
    \begin{align}
        \begin{aligned}
        &\big{\|}u(x_{1},0)\big{\|}_{L^{q}(-\infty,\infty)}\leq C(q)
        \big{\|}u\big{\|}_{W^{1,q}(\Omega_{\infty})},~
        \\
        &\big{\|}u(x_{1},1)\big{\|}_{L^{q}(-\infty,\infty)}\leq C(q)
        \big{\|}u\big{\|}_{W^{1,q}(\Omega_{\infty})},
        \end{aligned}
    \end{align}
    where $C(q)$ is a constant only depending on $q$.
\end{corollary}
\begin{lemma}\label{sobolveqianrudingli12241}
    Let $\mathbf{u}\in \mathbf{G}^{1,2}(\Omega_{R})$. Then for $q\in [2,+\infty)$, we have
    \begin{align*}
        \big{\|}\mathbf{u}\big{\|}_{L^{q}(\Omega_{R})}\leq C(q)\big{\|}\nabla\mathbf{u}\big{\|}_{L^2(\Omega_{R})}.
    \end{align*}
\end{lemma}
Similarly, one has the following corollary.
\begin{corollary}\label{sobolveqianrudingli0125}
    Let $\mathbf{u}\in \mathbf{G}^{1,2}(\Omega_{\infty})$. Then for $q\in [2,+\infty)$, we have
    \begin{align*}
        \big{\|}\mathbf{u}\big{\|}_{L^{q}(\Omega_{\infty})}\leq C(q)\big{\|}\nabla\mathbf{u}\big{\|}_{L^2(\Omega_{\infty})}.
    \end{align*}
\end{corollary}

\begin{lemma}\label{jiushizheyang12242}
    Let $\mathbf{u}\in$ $\mathbf{V}^{1,2}(\Omega_{R})$.
    Then there exist constants independent of $R$ such that
    \begin{align*}
        \begin{aligned}
        &~~~~~~~~~\big{\|}\nabla u_{2}\big{\|}_{L^q(\Omega_{R})}\leq C(q)\left\|\nabla^2 \mathbf{u}\right\|_{L^2(\Omega_{R})},~\text{for~}q\in [2,+\infty),
        \\
        &\big{\|}\partial_{2}u_{1}\big{\|}_{L^4(\Omega_{R})}\leq C\left(\left\|\nabla^2 u_{1}\right\|_{L^2(\Omega_{R})}^{\frac{1}{2}}
        \big{\|}\nabla u_{1}\big{\|}_{L^2(\Omega_{R})}^{\frac{1}{2}}+\big{\|}\nabla u_{1}\big{\|}_{L^2(\Omega_{R})}^{\frac{3}{4}}
        \left\|\nabla^2 u_{1}\right\|_{L^2(\Omega_{R})}^{\frac{1}{4}}\right).
        \end{aligned}
    \end{align*}
\end{lemma}
Recall that $\partial_{1}u_{1}+\partial_{2}u_{2}=0$, we thereby omit the estimate of $\partial_{1}u_{1}$ in Lemma
\ref{jiushizheyang12242}. Passing $R\rightarrow+\infty$, one has the following corollary.
\begin{corollary}\label{jiushizheyang0126}
    Let $\mathbf{u}\in$ $\mathbf{V}^{1,2}(\Omega_{\infty})$.
    Then we have
    \begin{align*}
        \begin{aligned}
        &~~~~~~~~~~~\big{\|}\nabla u_{2}\big{\|}_{L^q(\Omega_{\infty})}\leq C(q)\big{\|}\nabla^2 \mathbf{u}\big{\|}_{L^2(\Omega_{\infty})},
        ~\text{for~}q\in [2,+\infty),
        \\
        & \big{\|}\partial_{2}u_{1}\big{\|}_{L^4(\Omega_{\infty})}\leq C\left(\big{\|}\nabla^2 u_{1}\big{\|}_{L^2(\Omega_{\infty})}^{\frac{1}{2}}
        \big{\|}\nabla u_{1}\big{\|}_{L^2(\Omega_{\infty})}^{\frac{1}{2}}+\big{\|}\nabla u_{1}\big{\|}_{L^2(\Omega_{\infty})}^{\frac{3}{4}}
        \big{\|}\nabla^2 u_{1}\big{\|}_{L^2(\Omega_{\infty})}^{\frac{1}{4}}\right).
        \end{aligned}
    \end{align*}
\end{corollary}

\begin{lemma}\label{jiushizheyang12243}
    Let $u\in W^{1,q}(\Omega_{\infty})$ and $q\in(2,+\infty)$. Then one has
    \begin{align*}
        u\in L^{\infty}(\Omega_{\infty})~\text{and~}\big{\|}u\big{\|}_{L^{\infty}(\Omega_{\infty})}
        \leq C(q)\big{\|}u\big{\|}_{W^{1,q}(\Omega_{\infty})}.
    \end{align*}
\end{lemma}
By the similar method as Lemma  \ref{jiushizheyang12243}, one can obtain the following corollary.
\begin{corollary}\label{ohmygod1228}
    Let $u\in W^{1,q}(\Omega_{R})$ $(q\in(2,+\infty))$. Then we have
    \begin{align*}
        u\in L^{\infty}(\Omega_{R}),~~\big{\|}u\big{\|}_{L^{\infty}(\Omega_{R})}
        \leq C\big{\|}u\big{\|}_{W^{1,q}(\Omega_{R})}.
    \end{align*}
\end{corollary}

\begin{remark}The proofs of Lemmas \ref{jidingli}-\ref{jiushizheyang12243} and Corollary \ref{ohmygod1228}
can be found in Appendix \ref{11226}.
\end{remark}

Generally, we note that the pressure $p$ does not occur in the definition of weak solution. In order to investigate the existence of
$p$ and get the desired Stokes estimate, we shall consider the $q-$generalized solution (the definition can be
found in Appendix \ref{appendix11226}).

\begin{lemma}\label{tishengzhengzexing} If $(\mathbf{u},p)$ is $q-$generalized solution of problem
    \eqref{tuiguanjie}, then we have
$\mathbf{u}\in \mathbf{V}^{1,q}(\Omega_{R})$,~$p\in W^{1,q}(\Omega_{R})$. Moreover, one have the inequality
\begin{align}\label{xiangyaodejieguo}
    \left\|\nabla^2\mathbf{u}\right\|_{L^q(\Omega_{R})}+\big{\|}\nabla p\big{\|}_{L^q(\Omega_{R})}
    \leq C\left(\big{\|}\mathbf{f}\big{\|}_{L^{q}(\Omega_{R})}+\big{\|}\mathbf{u}\big{\|}_{W^{1,q}(\Omega_{R})}\right),
\end{align}
where $C$ is independent of $R$.
\end{lemma}
The proof
of Lemma \ref{tishengzhengzexing} can also be found in Appendix \ref{appendix11226}.
Since the constant $C$ in \eqref{xiangyaodejieguo} is independent of $R$, thus, passing
$R\rightarrow +\infty$ gives the following conclusion.
\begin{corollary}\label{tuilun1216}If $(\mathbf{u},p)$ is $q-$generalized solution of problem
    \eqref{tuiguanjie}, in which $R\rightarrow+\infty$, that is, the region is $\Omega_{\infty}$. Then we have
    \begin{align}\label{xiangyaodejieguo1}
        \left\|\nabla^2\mathbf{u}\right\|_{L^q(\Omega_{\infty})}+\big{\|}\nabla p\big{\|}_{L^q(\Omega_{\infty})}
        \leq C\left(\big{\|}\mathbf{f}\big{\|}_{L^{q}(\Omega_{\infty})}+\big{\|}\mathbf{u}\big{\|}_{W^{1,q}(\Omega_{\infty})}\right).
    \end{align}
\end{corollary}

For the construction of Galerkin approximate solutions, we consider the eigenvalues and
eigenvectors of Stokes operator on $\Omega_{R}$.
\begin{align}\label{eigen1}
    \begin{cases}
        -\mu\Delta \mathbf{u}+\nabla p=\lambda \mathbf{u},
        \\
        u_{2}(x_{1},1)=u_{2}(x_{1},0)=0,~ \mathbf{u}|_{\Gamma_{3 R}\cup\Gamma_{4R}}=\mathbf{0},
        \\
        \frac{\partial u_{1}}{\partial{x_{2}}}(x_{1},1)=\frac{k_{1}}{\mu}u_{1}(x_{1},1),~
        \frac{\partial u_{1}}{\partial{x_{2}}}(x_{1},0)=-\frac{k_{0}}{\mu}u_{1}(x_{1},0),
        \\
        \nabla\cdot \mathbf{u}=0.
    \end{cases}
\end{align}
\begin{lemma}\label{tezhenzhidewenti}
    About problem \eqref{eigen1}, one has the following conclusions,
    \begin{enumerate}[$(1)$]
    \item there exist several sequences $\{\mathbf{e}^n\}\subset [C^\infty(\Omega_{R})]^2$,
    $\{p^n\}\subset C^{\infty}(\Omega_{R})$ and $\{\lambda_{n}\}\subset \mathbf{R}^{+}$, such that
    \begin{align*}
        \begin{cases}
            -\mu\Delta \mathbf{e}^{n}+\nabla p^{n}=\lambda_{n} \mathbf{e}^{n},
            \\
            e_{2}^{n}(x_{1},1)=e_{2}^{n}(x_{1},0)=0,~ \mathbf{e}^{n}|_{\Gamma_{3 R}\cup\Gamma_{4R}}=\mathbf{0},
            \\
            \frac{\partial e_{1}^{n}}{\partial{x_{2}}}(x_{1},1)=\frac{k_{1}}{\mu}e_{1}^{n}(x_{1},1),~
            \frac{\partial e_{1}^{n}}{\partial{x_{2}}}(x_{1},0)=-\frac{k_{0}}{\mu}e_{1}^{n}(x_{1},0),
            \\
            \nabla\cdot \mathbf{e}^{n}=0,~\|\mathbf{e}^{n}\|_{L^{2}(\Omega_{R})}=1;
        \end{cases}
    \end{align*}
    \item $0<\lambda_{1}\leq \lambda_{2}\leq\cdots\leq\lambda_{n}\rightarrow +\infty$ as
    $n\rightarrow\infty$;
    \item $\{\mathbf{e}^{n}\}$ is a set of
    orthogonal basis in $\mathbf{G}^{1,2}(\Omega_{R})$. That is,
    $((\mathbf{e}^{n},\mathbf{e}^m))=\lambda_{n}\delta_{nm}$, where $\delta_{nm}$ is the Kronecker Symbol. And
    $\forall \mathbf{u}\in \mathbf{G}^{1,2}(\Omega_{R})$,
    $\mathbf{u}=\sum\limits_{k=1}^{+\infty}c_{k} \mathbf{e}^{k}$, where
    $c_{k}=\frac{\left(\left(\mathbf{u},\mathbf{e}^{k}\right)\right)}{\lambda_{k}}=\left(\mathbf{u},\mathbf{e}^{k}\right)$.
    \end{enumerate}
    \end{lemma}
    The proof of this lemma is given in Appendix \ref{appendixb1226}.
    Moreover, there are two inequalities which will be used in this paper and listed without proof.
    The reader can refer to see \cite{li_global_2021} for details.

    \begin{lemma}[$L^4-$ estimate]\label{L4guji}
        There exists a constant $C>0$, independent of $R$, such that
        \begin{align*}
            \big{\|}\mathbf{u}\big{\|}_{L^4(\Omega_{R})}^{2}\leq C\big{\|}\mathbf{u}\big{\|}_{L^2(\Omega_{R})}\big{\|}\nabla
            \mathbf{u}\big{\|}_{L^{2}(\Omega_{R})},~\forall~\mathbf{u}\in \mathbf{V}^{1,2}(\Omega_{R})  .      \end{align*}
    \end{lemma}
    \begin{lemma}[$L^\infty-$estimate]\label{Linftyguji}
        There exists a constant $C>0$, independent of $R$, such that
        \begin{align*}
            \big{\|}\mathbf{u}\big{\|}_{L^\infty(\Omega_{R})}^{2}\leq C\big{\|}\mathbf{u}\big{\|}_{L^2(\Omega_{R})}\big{\|}
            \mathbf{u}\big{\|}_{H^{2}(\Omega_{R})},~
        \forall~\mathbf{u}\in \mathbf{V}^{1,2}(\Omega_{R}).
        \end{align*}
    \end{lemma}
\section{Global well-posedness}\label{jiedecunzaixing0104}
This section is focused on the global well-posedness of problem \eqref{model}-\eqref{boundarycondition}. We firstly give the
conclusions on bounded region $\Omega_{R}$, and then derive some uniform estimates independent of $R$, which enable us to obtain the
corresponding conclusions on unbounded domain $\Omega_{\infty}$.
\subsection{Global well-posedness on bounded domain}
In the domain $\Omega_{R}$, the problem \eqref{model}-\eqref{boundarycondition} can be rewritten as follows,
\begin{align}\label{youjiequyu}
    \begin{cases}
     \frac{\partial \mathbf{u}}{\partial t}+(\mathbf{u}\cdot \nabla)\mathbf{u}+\nabla p
     =\mu\Delta\mathbf{u}+\theta
     \mathbf{g},
     \\
     \frac{\partial\theta}{\partial t}+(\mathbf{u}\cdot\nabla)\theta=0,
     \\
     \nabla\cdot\mathbf{u}=0,
     \\
     \mathbf{u}(x,0)=\mathbf{u}_{0}(x),~~\theta(x,0)=\theta_{0}(x),
    \end{cases}
\end{align}
subjected to the boundary conditions,
\begin{align}\label{naiverbouncondition}
    \begin{aligned}
    &u_{2}(x_{1},0)=u_{2}(x_{1},1)=0,~\partial_{2}u_{1}(x_{1},1)=\frac{k_{1}}{\mu}u_{1}(x_{1},1),
    \\
    &\partial_{2}u_{1}(x_{1},0)=-\frac{k_{0}}{\mu}u_{1}(x_{1},0),~x_{1}\in[-R,R],
    ~\mathbf{u}|_{\Gamma_{3R}\cup\Gamma_{4R}}=\mathbf{0}.
    \end{aligned}
\end{align}

We give the weak solution to  above problem as defined below.
\begin{definition}[Weak solution] \label{dingyiruoruo}($\mathbf{u},\theta$)
    is called the
    weak solution of problem \eqref{youjiequyu}-\eqref{naiverbouncondition},
    if
    $\mathbf{u}=(u_{1}(t,x),u_{2}(t,x))\in
    L^{2}\left(0,T;\mathbf{G}^{1,2}(\Omega_{R})\right)$, $\theta\in L^{\infty}([0,T]\times\Omega_{R})$ and satisfies,

   $(i)$  $\forall~ \mathbf{\Phi}(t,x)\in C^{1}\left([0,T];\mathbf{G}^{1,2}
    (\Omega_{R})\right)$ with $\mathbf{\Phi}(T,x)=\mathbf{0}$ and
    $\Psi\in C^{1}\left([0,T];H^{1}(\Omega_{R})\right)$
    with
    $\Psi(T,x)=0$,
    \begin{align}\label{ruojie}
        \begin{aligned}
        &\begin{aligned}
            -\int_{0}^{T}&\int_{\Omega_{R}}\mathbf{u}\cdot\partial_{t}\mathbf{\Phi} dx dt
            -\int_{0}^{T}\int_{\Omega_{R}}(\mathbf{u}\cdot\nabla)\mathbf{\Phi}\cdot \mathbf{u}
            dx dt
         +\mu\int_{0}^{T}\int_{\Omega_{R}}\nabla\mathbf{u}\cdot\nabla\mathbf{\Phi} dx dt
         \\
         &-\int_{0}^{T}\int_{-R}^{R}\left[k_{1}u_{1}(x_{1},1)\mathbf{\Phi}_{1}(x_{1},1)+
         k_{0}u_{1}(x_{1},0)\mathbf{\Phi}_{1}(x_{1},0)\right]dx_{1}dt
         \\
         &~=\int_{\Omega_{R}}\mathbf{u}_{0}\cdot\mathbf{\Phi}(0,x)dx+
         \int_{0}^{T}\int_{\Omega_{R}}\theta\mathbf{g}\cdot\mathbf{\Phi}  dxdt,
        \end{aligned}
        \\
        &\begin{aligned}
            -\int_{0}^{T}\int_{\Omega_{R}}\theta\partial_{t}\Psi dxdt
            -\int_{0}^{T}\int_{\Omega_{R}}(\mathbf{u}\cdot\nabla)\Psi\theta dxdt
         =\int_{\Omega_{R}}\theta_{0}(x)\Psi(0,x)dx;
        \end{aligned}
    \end{aligned}
    \end{align}

  $(ii)$ the initial data $\mathbf{u}_{0}(x)\in \mathbf{G}^{1,2}(\Omega_{R})$
    and $\theta_{0}(x)\in L^{\infty}(\Omega_{R})$.
\end{definition}
Now we state the main results in the following.
\begin{theorem}\label{cunzairuojie} There exists at least one weak solution
    $(\mathbf{u},\theta)$ to problem \eqref{youjiequyu}-\eqref{naiverbouncondition}, when
    $\mathbf{u}_{0}\in \mathbf{G}^{1,2}(\Omega_{R})$ and
     $\theta_{0}\in L^{\infty}(\Omega_{R})$.     Moreover, $\mathbf{u}$ and $\theta$ satisfy the inequalities,
     \begin{align*}
        \begin{aligned}
        & \big{\|}\mathbf{u}\big{\|}_{L^2(\Omega_{R})}^2\leq \left(\big{\|}\theta_{0}\big{\|}_{L^2(\Omega_{R})}^2
        +\big{\|}\mathbf{u}_{0}\big{\|}_{L^2(\Omega_{R})}^{2}\right)e^{|\overrightarrow{
            \mathbf{g}
        }|T},
        ~\big{\|}\theta\big{\|}_{L^{k}(\Omega_{R})}=\big{\|}\theta_{0}\big{\|}_{L^{k}(\Omega_{R})}, (k\in(1,+\infty)),
        \\
        &2\mu\int_{0}^{T}\big{\|}\nabla\mathbf{u}\big{\|}^{2}_{L^2(\Omega_{R})}ds
        \leq
        \big{\|}\mathbf{u}_{0}\big{\|}_{L^2(\Omega_{R})}^2+
       \left(\big{\|}\theta_{0}\big{\|}_{L^2(\Omega_{R})}^2+\big{\|}\mathbf{u}_{0}\big{\|}_{L^2(\Omega_{R})}^2\right)e^{|\mathbf{g}|T}
       +|\mathbf{g}|\big{\|}\theta_{0}\big{\|}_{L^2(\Omega_{R})}^2T.
        \end{aligned}
     \end{align*}
\end{theorem}
Note that the existence of $p$ follows immediately from $\eqref{youjiequyu}_{1}$ and $\eqref{youjiequyu}_{3}$ by
a classical method, see Lemma \ref{dengjia}. For more regular initial value $\mathbf{u}_{0}$ and $\theta_{0}$, we can
improve the regularity of weak solution to obtain the strong solution.
\begin{theorem}\label{bufenqiangjie1227}Under the condition of Theorem \ref{cunzairuojie},
    and $\mathbf{u}_{0}\in \mathbf{V}^{1,2}(\Omega_{R})$, then we have
    \begin{align*}
        \begin{aligned}
    &\partial_{t}\mathbf{u}\in L^{\infty}\left(0,T;\widehat{\mathbf{L}}^{2}(\Omega_{R})\right)
\cap L^2\left(0,T;\mathbf{G}^{1,2}(\Omega_{R})\right),~
\mathbf{u}\in L^{\infty}\left(0,T;\mathbf{V}^{1,2}(\Omega_{R})\right), p\in L^{\infty}\left(0,T;H^{1}(\Omega_{R})\right),
        \\
        &\nabla^2\mathbf{u}\in L^{2}\left(0,T; \left[L^4(\Omega_{R})\right]^6\right),
        ~\nabla p\in L^{2}\left(0,T;\left[L^{4}(\Omega_{R})\right]^2\right).
        \end{aligned}
    \end{align*}
\end{theorem}

\begin{theorem}\label{qiangjie} Under the condition of Theorem \ref{bufenqiangjie1227}
    and $\theta_{0}\in H^{1}(\Omega_{R})$, then we have
    \begin{align*}
        \begin{aligned}
        \theta\in L^{\infty}\left(0,T;H^{1}(\Omega_{R})\right),
\theta_{t}\in L^{\infty}\left(0,T;L^{2}(\Omega_{R})\right).
        \end{aligned}
    \end{align*}
    Moreover, the following inequality holds,
    \begin{align}\label{kanqilaimaodun}
        \begin{aligned}
            \big{\|}\partial_{t}\mathbf{u}\big{\|}_{L^2(\Omega_{R})}^{2}&+
            \big{\|}\nabla p\big{\|}_{L^{2}(\Omega_{R})}^2+
            \big{\|}\nabla\mathbf{u}\big{\|}_{H^{1}(\Omega_{R})}^{2}
            +\big{\|}\theta_{t}\big{\|}_{L^2(\Omega_{R})}^2+
            \int_{0}^{t}\bigg{[}\big{\|}\partial_{s}\mathbf{u}(s)\big{\|}_{L^{2}(\Omega_{R})}^{2}
            \\
            &~~+\big{\|}\nabla\partial_{s}\mathbf{u}\bigg{\|}_{L^2(\Omega_{R})}^{2}
            +\big{\|}\nabla^2\mathbf{u}\big{\|}_{L^{4}(\Omega_{R})}^{2}+\big{\|}\nabla p\big{\|}_{L^{4}(\Omega_{R})}^2\bigg{]}ds
            \leq C,
        \end{aligned}
    \end{align}
    where $C$ is independent of $R$.
\end{theorem}
Furthermore, we can verify that the strong solution is unique.
\begin{theorem}\label{weiyixing1227} Under the condition of Theorem \ref{qiangjie}, the strong
     solution $(\mathbf{u},\theta)$ is unique.
\end{theorem}

\subsection{Global well-posedness on unbounded domain}\label{20240319}
Since all constants of the inequalities in Theorems \ref{cunzairuojie} and \ref{qiangjie}
are independent of $R$, thus one can let $R\rightarrow\infty$ to obtain the similar conclusions to problem
\eqref{model}-\eqref{boundarycondition}.

We firstly give the definition of weak solution on $\Omega_{\infty}$.
\begin{definition}[Weak~solution]\label{dingyiruoruo1}
   ($\mathbf{u},\theta$)
    is called the weak solution of problem \eqref{model}-\eqref{boundarycondition},
    if
    $\mathbf{u}=(u_{1},u_{2})\in
    L^{2}\left(0,T;\mathbf{G}^{1,2}(\Omega_{\infty})\right)$ $\cap L^{\infty}\left(0,T;\left[L^{2}(\Omega_{\infty})\right]^2\right)$
    and $\theta\in L^{\infty}([0,T]\times\Omega_{\infty})$ and satisfies,

    $(i)$
    $\forall~ \mathbf{\Phi}(t,x)\in C^{1}\left([0,T];\mathbf{C}_{0,x_{1}}(\Omega_{\infty})\right)$,
    $\Psi\in C^{1}\left([0,T];H^{1}(\Omega_{\infty})\right)$,
    $\mathbf{\Phi}(T,x)=0$, and
    $\Psi(T,x)=0$,
    \begin{align}\label{ruojie1}
        \begin{aligned}
        &\begin{aligned}
            -\int_{0}^{T}&\int_{\Omega_{\infty}}\mathbf{u}\cdot\partial_{t}\mathbf{\Phi} dx dt
            -\int_{0}^{T}\int_{\Omega_{\infty}}(\mathbf{u}\cdot\nabla)\mathbf{\Phi}\cdot \mathbf{u}
         -\mu\int_{0}^{T}\int_{\Omega_{\infty}}\nabla\mathbf{u}\cdot\nabla\mathbf{\Phi} dx dt
         \\
         &-\int_{0}^{T}\int_{-\infty}^{\infty}\left[k_{1}u_{1}(x_{1},1)\mathbf{\Phi}_{1}(x_{1},1)+
         k_{0}u_{1}(x_{1},0)\mathbf{\Phi}_{1}(x_{1},0)\right]dx_{1}dt
         \\
         &~=\int_{\Omega_{\infty}}\mathbf{u}_{0}\mathbf{\Phi}(0,x)dx+
         \int_{0}^{T}\int_{\Omega_{\infty}}\theta\overrightarrow{g}\cdot\mathbf{\Phi}  dxdt,
        \end{aligned}
        \\
        &\begin{aligned}
            -\int_{0}^{T}\int_{\Omega_{\infty}}\theta\partial_{t}\Psi dxdt
            -\int_{0}^{T}\int_{\Omega_{\infty}}(\mathbf{u}\cdot\nabla)\Psi\theta dxdt
         =\int_{\Omega_{\infty}}\theta_{0}(x)\Psi(0,x)dx.
        \end{aligned}
    \end{aligned}
    \end{align}

     $(ii)$ the initial values $\mathbf{u}_{0}(x)\in \mathbf{V}^{1,2}(\Omega_{\infty})$
    and $\theta_{0}(x)\in L^{\infty}(\Omega_{\infty})\cap L^{1}(\Omega_{\infty})$.
    \begin{remark}
        From $\theta_{0}(x)\in L^{\infty}(\Omega_{\infty})\cap L^{1}(\Omega_{\infty})$, one can conclude that
        $\theta_{0}(x)\in L^{k}(\Omega_{\infty})(~\rm{\forall}k\geq 1)$.
    \end{remark}
\end{definition}
\begin{theorem}\label{cunzairuojie1} Let $\mathbf{u}_{0}\in \mathbf{G}^{1,2}(\Omega_{\infty})$ and
    $\theta_{0}\in L^{\infty}(\Omega_{\infty})\cap L^{1}(\Omega_{\infty})$. Then
   the problem \eqref{model}-\eqref{boundarycondition} has a weak solution $(\mathbf{u},\theta)$.
    Moreover, $\mathbf{u}$ and $\theta$ satisfy the following inequalities,
    \begin{align}\label{threeev}
       \begin{aligned}
       &\big{\|}\mathbf{u}\big{\|}_{L^2}^2\leq \left(\big{\|}\theta_{0}\big{\|}_{L^2}^2
       +\big{\|}\mathbf{u}_{0}\big{\|}_{L^2}\right)e^{|
           \mathbf{g}
       |T},
       ~\big{\|}\theta\big{\|}_{L^{k}}=\big{\|}\theta_{0}\big{\|}_{L^{k}},(k\in(1,+\infty)),
       \\
       &2\mu\int_{0}^{T}\big{\|}\nabla\mathbf{u}\big{\|}^{2}_{L^2}ds
       \leq
       \big{\|}\mathbf{u}_{0}\big{\|}_{L^2}^2+
      \left(\big{\|}\theta_{0}\big{\|}_{L^2}^2+\big{\|}\mathbf{u}_{0}\big{\|}_{L^2}^2\right)e^{|\mathbf{g}|T}
      +|\mathbf{g}|\big{\|}\theta_{0}\big{\|}_{L^2}^2T.
       \end{aligned}
    \end{align}
\end{theorem}
\begin{theorem}\label{qiangjie1} Let $\mathbf{u}_{0}$ $\in \mathbf{V}^{1,2}(\Omega_{\infty})
   $
   and $\theta_{0}\in H^{1}(\Omega_{\infty})\cap L^{1}(\Omega_{\infty})
   \cap L^{\infty}(\Omega_{\infty})$.
   Then the problem \eqref{model}
   -\eqref{boundarycondition} has a unique strong solution $(\mathbf{u},\nabla p,\theta)$ satisfying
   \begin{align}\label{312}
       \begin{aligned}
       &\partial_{t}\mathbf{u}\in L^{\infty}\left(0,T;\widehat{\mathbf{L}}(\Omega_{\infty})\right)
\cap L^2\left(0,T;\mathbf{G}^{1,2}(\Omega_{\infty})\right),~ \theta\in L^{\infty}\left(0,T;H^{1}(\Omega_{\infty})\right),
       \\
       &\theta_{t}\in L^{\infty}\left(0,T;L^{2}(\Omega_{\infty})\right),~\mathbf{u}\in L^{\infty}\left(0,T;
       \mathbf{V}^{1,2}(\Omega_{\infty})\right),
       ~\nabla^2\mathbf{u}\in L^{2}\left(0,T; \left[L^4(\Omega_{\infty})\right]^6\right),
       \\
       &\nabla p\in L^{\infty}\left(0,T;\left[L^2(\Omega_{\infty})\right]^2\right)
       \cap  L^{2}\left(0,T;\left[L^{4}(\Omega_{\infty})\right]^2\right),
       \end{aligned}
   \end{align}
   and
   \begin{align}\label{kanqilaimaodun1}
       \begin{aligned}
           \big{\|}\partial_{t}\mathbf{u}&\big{\|}_{L^2(\Omega_{\infty})}^{2}+
           \big{\|}\nabla p\big{\|}_{L^{2}(\Omega_{\infty})}^2+
           \big{\|}\nabla\mathbf{u}\big{\|}_{H^{1}(\Omega_{\infty})}^{2}
           +\big{\|}\theta_{t}\big{\|}_{L^2(\Omega_{\infty})}^2+
           \int_{0}^{t}\bigg{[}\big{\|}\partial_{s}\mathbf{u}(s)\big{\|}_{L^{2}(\Omega_{\infty})}^{2}
           \\
           &+\big{\|}\nabla\partial_{s}\mathbf{u}\big{\|}_{L^2(\Omega_{\infty})}^{2}
           +\left\|\nabla^2\mathbf{u}\right\|_{L^{4}(\Omega_{\infty})}^{2}+\big{\|}\nabla p\big{\|}_{L^{4}(\Omega_{\infty})}^2\bigg{]}ds
           \leq C.
       \end{aligned}
   \end{align}
\end{theorem}
\begin{remark}Note that the conclusions of Theorem \ref{qiangjie1} hold for any $T>0$, thus the global existence of strong solution is verified.
\end{remark}

\subsection{Proofs of conclusions on bounded domain}
The aim of this section is devoted to the proof of Theorems \ref{cunzairuojie}-\ref{weiyixing1227} for problem \eqref{youjiequyu}-\eqref{naiverbouncondition}. Here, the semi-Galerkin method
\cite{kim_weak_1987} (combination of Galerkin approximation and characteristic methods) is employed to show the existence of weak solution.
Firstly, we show that there exists
a continuous and compact map $\mathbf{F}^m$$(m\in\mathbf{N}^{*})$ (see Lemmas \ref{youjiediyi}-\ref{youjieyinli}) on bounded and convex set; Then by the
Schauder fixed point theorem (see P95 in \cite{daoxingxia2011}), one can prove that there exists
approximate solution
$\{\mathbf{u}^{m},\theta^{m}\}$(see Remark \ref{bijinjie0402}); Furthermore, one derives a prior estimate for $\{\mathbf{u}^{m},\theta^{m}\}$, which implies the existence of limits for approximate solution, and the limits meet the equality \eqref{ruojie}, that
is, these limits are the weak solution. Finally, the standard procedure can be used to improve the regularity of
weak solution, thus one can get the strong solution.

To order to get $\mathbf{F}^{m}$, for a fixed $\mathbf{u}$, we first consider the equation $\eqref{youjiequyu}_{2}$ in forms,
\begin{align}\label{kao1}
    \begin{cases}
        \partial_{t}\theta+(\mathbf{u}\cdot\nabla)\theta=0,
        \\
        \theta(0,x)=\overline{\theta}_{0}(x), ~~x\in\Omega_{R},
    \end{cases}
\end{align}
where $\overline{\theta}_{0}(x)$ is the mollification of $\theta_{0}(x)$.
It is worth mentioning that, by the properties of mollification, one can replace $\theta_{0}$ by $\overline{\theta}_{0}$
involving with integral.
\begin{lemma}\label{youjiediyi}
    For all $(t,x)\in [0,T]\times \overline{\Omega}_{R}$, let
    $\mathbf{u}(t,x)\in C\left(0,T;\left(C^{1}(\Omega_{R})\right)^{2}\cap \mathbf{G}^{1,2}(\Omega_{R})\right)$ and
    $\overline{\theta}_{0}\in C^{1}(\Omega_{R}).$ Then, the equation \eqref{kao1} has a unique solution
    $\theta(t,x)\in C^{1}\left([0,T]\times \overline{\Omega}_{R}\right).$
\end{lemma}
\begin{proof}
    Let $E$ is an open ball in $\mathbf{R}^2$, and $\overline{\Omega}_{R}\subset E$.
    By the continuity of $\mathbf{u}$, we can extend $\mathbf{u}$ from $C\left(0,T;
    \left(C^{1}(\overline{\Omega}_{R})\right)^2\right)$
    to $C\left(0,T; \left(C^{1}(\overline{E})\right)^2\right)$
    denoted by $\mathbf{w}$, that is, for all $(t,x)\in [0,T]\times \overline{\Omega}_{R}$,
    $\mathbf{u}\equiv \mathbf{w}$. Now we shall solve \eqref{kao1} by the classical method of characteristics.
    We first consider the following equations,
    \begin{align}\label{kao3}
        \begin{cases}
            \frac{dX_{1}}{dt}=w_{1}(t,X_{1},X_{2}),
            \\
            \frac{dX_{2}}{dt}=w_{2}(t,X_{2},X_{2}),
            \\
            X(0)=(X_{1}(0),X_{2}(0))=y=(y_{1},y_{2})\in \overline{\Omega}_{R}.
        \end{cases}
    \end{align}
    By Cauchy-Lipschitz theorem (see A.3 in \cite{bedrossian_mathematical_2022}), one can conclude that
    the equation \eqref{kao3} has a unique solution
    $X(t,y)\in C^{1}\left(0,\widetilde{T};
    \left(C^{1}(\overline{\Omega}_{R})\right)^2\right)$, where
    $0<\widetilde{T}\leq T.$

    One hand, if $y\in \Gamma_{3R}\cup\Gamma_{4R}$, then $\mathbf{w}(t,y)=\mathbf{u}(t,y)
    =\mathbf{0}$, which implies that when $t\in[0,\widetilde{T}]$, $X(t,y)=y\in
    \Gamma_{3R}\cup \Gamma_{4R}$. On the other hand,
    if $y\in \Gamma_{1R}\cup \Gamma_{2R}$ which indicates that $w_{2}(t,y)=u_{2}(t,y)=0$, then
     one has $X_{2}(t,y)=y_{2}$ for all $t\in [0,\widetilde{T}]$. In
     conclusion, when $y\in \partial\Omega_{R}$, we have $X(t,y)\in
    \partial\Omega_{R}$,
     for all $t\in [0,\widetilde{T}]$.
    Then by the uniqueness, when initial value $y\in \Omega_{R}$, we have
    $X(t,y)\in \overline{\Omega}_{R}$ for all $t\in[0,\widetilde{T}]$. Thus, for problem \eqref{kao3},
    we can take $\widetilde{T}=T$ and $\mathbf{w}=\mathbf{u}.$

    Moreover, due to $\nabla\cdot\mathbf{u}=0$ for all $(t,y)\in [0,T]\times\overline{\Omega}_{R},$
    one has $\text{det}\left\{\frac{\partial X_{i}}{\partial y_{j}}\right\}=1$. And since $X(t,y)\in
\left(C^1([0,T]\times\overline{\Omega}_{R})\right)^2$, then $X(t,\cdot):$ $\overline{\Omega}_{R}
    \rightarrow \overline{\Omega}_{R}$ is $C^1$-diffeomorphism. Thus, there exists an inverse map
    $A(t,\cdot)$, satisfying $A(0,x)=x$,
    $y=A\left(t,X\left(t,y\right)\right)$ and
    $x=X\left(t,A\left(t,x\right)\right)$,
    where $x\in \overline{\Omega}_{R}$.

    Now, we will verify that $\theta(t,x)=\overline{\theta}_{0}
    (A(t,x))$ is the solution of equation
    \eqref{kao1}.

    From
    $y=A(t,X(t,y))$ and $x=X(t,A(t,x))$,
    one can get
    \begin{align*}
        \begin{aligned}
            &\delta_{ij}=\frac{\partial A_{i}}{\partial y_{j}}=\sum\limits_{k=1}^2
            \frac{\partial A_{i}}{\partial x_{k}}
            \cdot\frac{\partial X_{k}}{\partial y_{j}},
            ~0=\partial_{t}X_{k}+\sum\limits_{i=1}^{2}
            \frac{\partial X_{k}}{\partial A_{i}}\frac{\partial A_{i}}{\partial t},
        \end{aligned}
    \end{align*}
    which implies that
    \begin{align*}
        \frac{\partial A_{i}}{\partial t}+\sum\limits_{k=1}^{2}u_{k}\frac{\partial A_{i}}{\partial x_{k}}=0,
    \end{align*}
    that is,
    \begin{align*}
        \frac{\partial A}{\partial t}+(\mathbf{u}\cdot\nabla)A=0.
    \end{align*}
 A series of calculation gives that
    \begin{align*}
        \partial_{t}\theta=\sum\limits_{i=1}^{2}
        \frac{\partial \overline{\theta}_{0}}{\partial y_{i}}\frac{\partial A_{i}}{\partial t},~
        \partial_{x_{k}}\theta=\sum\limits_{i=1}^{2}\frac{\partial \overline{\theta}_{0}}{\partial y_{i}}
        \frac{\partial A_{i}}{\partial x_{k}}.
    \end{align*}
    Thus, there appears the relation
    \begin{align*}
        \partial_{t}\theta+(\mathbf{u}\cdot\nabla)\theta=\sum\limits_{i=1}^{2}
        \frac{\partial \overline{\theta}_{0}}{\partial y_{i}}\left[\frac{\partial A_{i}}{\partial t}
        +(\mathbf{u}\cdot\nabla)A_{i}\right]=0.
    \end{align*}
    Besides, we know that $\overline{\theta}_{0}(A(0,x))=
    \overline{\theta}_{0}(x)$. Hence, $\overline{\theta}_{0}(A(t,x))$
    is the solution of equation \eqref{kao1}.
    \end{proof}
    \begin{lemma}\label{youjiedier}
        Suppose $\mathbf{u}^{k}(t,x)\in C\left(0,T;\left(C^{1}\left(\overline{\Omega}_{R}\right)\right)^{2}\cap \mathbf{G}^{1,2}(\Omega_{R})\right)$ and
        $\overline{\theta}_{0}\in C^{1}\left(\overline{\Omega}_{R}\right),$ where $k=1,2,\cdots$, and $(t,x)\in [0,T]\times \overline{\Omega}_{R}$.
       Assume further $\mathbf{u}^{k}$ converges to $\mathbf{u}$ in $C\left(0,T;\left(C^{1}\left(\overline{\Omega}_{R}\right)\right)^{2}\right)$.
        If $\theta^{k}(t,x)$ and $\theta(t,x)$ are the solutions of the following
        equations, respectively,
        \begin{align}\label{fangcheng0221}
                \begin{cases}
                    \partial_{t}\theta^{k}+(\mathbf{u}^{k}\cdot\nabla)\theta^{k}=0,
                    \\
                    \theta^{k}(0,x)=\overline{\theta}_{0},
                \end{cases}~ \text{and}~~~
                \begin{cases}
                    \partial_{t}\theta+(\mathbf{u}\cdot\nabla)\theta=0,
                    \\
                    \theta(0,x)=\overline{\theta}_{0},
                \end{cases}
        \end{align}
         then   $\theta^{k}(t,x)$ converges to $\theta(t,x)$ in
         $C\left([0,T]\times\overline{\Omega}_{R}\right)$.
    \end{lemma}
    \begin{proof}
        Let $X^{k}(t,y)$ and $X(t,y)$ solve the following equations, respectively,
        \begin{align*}
            \begin{cases}
                \frac{dX^{k}}{dt}=\mathbf{u}^{k}\left(t,X^{k}\right),
                \\
                X^{k}(0,y)=y\in\overline{\Omega}_{R},
            \end{cases}~~\text{and}~~~
            \begin{cases}
                \frac{dX}{dt}=\mathbf{u}\left(t,X\right),
                \\
                X(0,y)=y\in\overline{\Omega}_{R}.
            \end{cases}
        \end{align*}
    We also take $A^{k}\left(t,X^{k}(t,y)\right)=A\left(t,X(t,y)\right)=y$.

    In order to get the convergence of
    $\{X^{k}\}$, we need to investigate
     the following equations,
    \begin{align*}
        \begin{cases}
            \frac{d\left(X^{k}-X\right)}{dt}=\mathbf{u}^{k}\left(t,X^k\right)-\mathbf{u}\left(t,X\right),
            \\
            \left(X^{k}-X\right)(0,y)=\mathbf{0}.
        \end{cases}
    \end{align*}
    Integrating the above equality from
    $0$ to $t$ gives that
    \begin{align*}
        \left|X^{k}-X\right|\leq \int_{0}^{t}\left|\mathbf{u}^{k}(s,X^{k})-
        \mathbf{u}(s,X^{k})\right|ds+\int_{0}^{t}\left|\mathbf{u}(s,X^{k})-\mathbf{u}(s,X)\right|ds.
    \end{align*}
    Since $\mathbf{u}^{k}$ converges to $\mathbf{u}$ in $C\left(0,T;\left(C^{1}\left(\overline{\Omega}_{R}\right)\right)^{2}\right)$,
    then there exist $c>0$ and $\varepsilon_{k}>0$ such that
    \begin{align}\label{jiuzheyang}
        \left|X^{k}-X\right|\leq \varepsilon_{k}T+c\int_{0}^{t}\left|X^{k}-X\right|ds.
    \end{align}
    From Grownwall's inequality, one has
    \begin{align*}
        \int_{0}^{t}\left|X^{k}-X\right|ds\leq \frac{\varepsilon_{k}T}{c}\left(e^{cT}-1\right)\rightarrow 0,~
        \text{as}~k\rightarrow+\infty.
    \end{align*}
    Thus, from \eqref{jiuzheyang}, one can verify that
    $X^{k}(t,y)$ converges to $X(t,y)$ in $C\left([0,T]\times \overline{\Omega}_{R}\right)$.
As a result, $A^{k}(t,x)$ converges to $A(t,x)$ in $C\left([0,T]\times
    \overline{\Omega}_{R}\right)$, where $A^{k}$ and
    $A$ are the inverse maps of
    $X^{k}$ and $X$, respectively.

    From the proof of Lemma \ref{youjiediyi}, we conclude that $\theta^k(t,x)=\overline{\theta}_{0}\left(A^{k}(t,x)\right)$ and
    $\theta(t,x)=\overline{\theta}_{0}(A(t,x))$ are the solutions of $\eqref{fangcheng0221}_{1}$ and $\eqref{fangcheng0221}_{2}$, respectively. Therefore, one has
    $\theta^k\rightarrow \theta$ in $C\left([0,T]\times \overline{\Omega}_{R}\right)$
    as $k\rightarrow \infty$.
    \end{proof}
    To construct a bounded convex set,  for each fixed $m\in\mathbf{N}^{*}$, we can define the following space,
\begin{align*}
    \mathbf{Y}_{m}=\bigg{(}C([0,T])\bigg{)}^{m},
\end{align*}
endowed with the norm,
\begin{align*}
\big{\|}\mathbf{f^{m}}\big{\|}_{\mathbf{Y}_{m}}=\left(\sum\limits_{i=1}^{m}\max_{t\in [0,T]}
    |f_{mi}(t)|^2\right)^{\frac{1}{2}},
\end{align*}
 where $\mathbf{f^{m}}(t)=\left(f_{m1}(t),\cdots,f_{mm}(t)\right)\in \mathbf{Y}_{m}$.
It is easy to check that $\mathbf{Y}_{m}$ is a Banach space. We choose a convex closed subset of
$\mathbf{Y}_{m}$ in form
\begin{align*}
    B_{\widetilde{R}\mathbf{Y}_{m}}=\left\{\mathbf{f^m}\in\mathbf{Y}_{m}\big{|}~\big{\|}\mathbf{f}^{m}\big{\|}_{\mathbf{Y}_{m}}
    \leq \widetilde{R}\right\},
\end{align*}
where the constant $\widetilde{R}$ is to be determined later.

Now, we shall show how to construct the compact and continuous map $\mathbf{F}^{m}:B_{\widetilde{R}\mathbf{Y}_{m}}\rightarrow B_{\widetilde{R}\mathbf{Y}_{m}}$. Firstly,
 $\forall~\mathbf{f^m}=\left(f_{m1}(t),\cdots,f_{mm}(t)\right)\in B_{\widetilde{R}\mathbf{Y}_{m}}$, we can
introduce $\mathbf{U}^{m}=\left(U_{1}^{m},U_{2}^{m}\right)=\sum\limits_{i=1}^{m}f_{mi}(t)\mathbf{e}^{i},$ and easily find that
 $\mathbf{U}^{m}\in C\left(0,T;\left(C^{1}\left(\overline{\Omega}_{R}\right)\right)^2\cap\mathbf{G}^{1,2}(\Omega_{R})\right)$. Then,  for the following equation
 \begin{align}\label{youjiedisan}
    \begin{cases}
    \partial_{t}\theta^{m}+(\mathbf{U}^{m}\cdot\nabla)\theta^{m}=0,
    \\
    \theta^{m}(0,x)=\overline{\theta}_{0},
    \end{cases}
 \end{align}
we have from Lemma \ref{youjiediyi} that there exists $\theta^{m}
\in C^1\left([0,T]\times \overline{\Omega}_{R}\right)$ satisfying \eqref{youjiedisan}.
 By \eqref{youjiedisan}, a direct calculation shows that $\forall k\geq 2$,
\begin{align*}
    \frac{d}{dt}\big{\|}\theta^{m}\big{\|}_{L^k(\Omega_{R})}^k=0,
\end{align*}
which implies that
\begin{align}\label{xuyaoni0127}
    \big{\|}\theta^{m}\big{\|}_{L^k(\Omega_{R})}=\left\|\overline{\theta}_{0}\right\|_{L^{k}(\Omega_{R})},
    ~~\big{\|}\theta^{m}\big{\|}_{L^\infty(\Omega_{R})}=\left\|\overline{\theta}_{0}\right\|_{L^{\infty}(\Omega_{R})}.
\end{align}

Secondly, based on previous analysis, we further consider equation $\eqref{youjiequyu}_{1}$ under $\theta=\theta^{m}$. From Lemma
\ref{tezhenzhidewenti}, note that $\mathbf{G}^{1,2}(\Omega_{R})$ has an orthogonal basis
$\left\{\mathbf{e}^{i}\right\}_{i=1}^{\infty}\subset [C^{\infty}(\Omega_{R})]^2$. Thus, for fixed
$m\in\mathbf{N}_{+}$, we need to study the following ODE,
\begin{align}\label{changweifen}
    \begin{cases}
       \frac{d \widetilde{f}_{mi}(t)}{dt}
       +\sum\limits_{j,k=1}^{m}\widetilde{f}_{mj}(t)\widetilde{f}_{mk}(t)
       \int_{\Omega_{R}}\left(\mathbf{e}^{j}\cdot\nabla\right)\mathbf{e}^{k}\cdot\mathbf{e}^{i}dx
       +\lambda_{i}\widetilde{f}_{mi}(t)-\int_{\Omega_{R}}\theta^{m}\mathbf{g}
       \cdot\mathbf{e}^{i}dx=0,
       \\
       \widetilde{f}_{mi}(0)=\left(\mathbf{u}_{0},\mathbf{e}^{i}\right),~~i\in\{1,2,\cdots,m\}.
    \end{cases}
\end{align}
From the theory of ODE, there exists a $T_{m}>0$ such that
the above equation has a unique solution $\widetilde{\mathbf{f^{m}}}(t)=
\left(\widetilde{f}_{m1}(t),\cdots, \widetilde{f}_{mm}(t)\right)\in \left(C^{1}([0,T_{m}))\right)^{m},$ where $T_{m}$ is the maximum existence time of $\widetilde{\mathbf{f}}^{m}$. Thus, we can introduce
$\mathbf{u}^{m}=\sum\limits_{i=1}^{m}\widetilde{f}_{mi}(t)\mathbf{e}^{i}$.
Then
$\mathbf{u}^{m}$ satisfies the following equation,
\begin{align}\label{gujiyoujie}
    \begin{cases}
     \int_{\Omega_{R}}\partial_{t}\mathbf{u}^{m}\mathbf{e}^{i}dx
     +\int_{\Omega_{R}}\left(\mathbf{u}^{m}\cdot\nabla\right)\mathbf{u}^{m}\cdot \mathbf{e}^{i}dx
     =\mu\int_{\Omega_{R}}\Delta\mathbf{u}^{m}\cdot \mathbf{e}^{i}dx
     +\int_{\Omega_{R}}\theta^{m}\mathbf{g}\cdot\mathbf{e}^{i}dx,
     \\
     \mathbf{u}^{m}(0,x)=\sum\limits_{i=1}^{m}\left(\mathbf{u}_{0},\mathbf{e}^{i}\right)\mathbf{e}^{i}.
    \end{cases}
\end{align}
Obviously, $T_{m}\leq T$. Note that if $T_{m}<T$, then there exists a $i\in\{1,\cdots,m\}$ such that as $t\rightarrow T_{m}$, $\widetilde{f}_{mi}(t)
\rightarrow \infty$.
Next, we verify that $T_{m}=T$ and determine $\widetilde{R}$ by a prior estimate.
From \eqref{xuyaoni0127} and \eqref{gujiyoujie}, one has
\begin{align}\label{gujiyoujie12}
    \begin{aligned}
    \frac{1}{2}\frac{d}{dt}\big{\|}\mathbf{u}^{m}\big{\|}_{L^2(\Omega_{R})}^2
    &+\mu\big{\|}\nabla\mathbf{u}^{m}\big{\|}_{L^2(\Omega_{R})}^{2}
    -\int_{-R}^{R}\left[k_{1}|u^{m}_{1}(x_{1},1)|^2+k_{0}|u_{1}^{m}(x_{1},0)|^2\right]dx_{1}
    \\
&=\int_{\Omega_{R}}\theta^{m}\mathbf{g}\cdot \mathbf{u}^{m} dx\leq \frac{g}{2}\left(\big{\|}\mathbf{u}^m\big{\|}_{L^2(\Omega_{R})}^2
    +\left\|\overline{\theta}_{0}\right\|_{L^2(\Omega_{R})}^2\right).
    \end{aligned}
\end{align}
From Gronwall's inequality, one can obtain
\begin{align*}
    \sum\limits_{i=1}^{m}\left|\widetilde{f}_{mi}(t)\right|^2=\big{\|}\mathbf{u}^{m}\big{\|}_{L^2(\Omega_{R})}^2\leq e^{gT}\left(\big{\|}\mathbf{u}_{0}\big{\|}_{L^2(\Omega_{R})}^2
    +\left\|\overline{\theta}_{0}\right\|_{L^2(\Omega_{R})}^2\right)<+\infty,
\end{align*}
which implies that  $T_{m}<T$ is impossible.
Thus, $T_{m}=T$.

Moreover, let $e^{gT}\left(\big{\|}\mathbf{u}_{0}\big{\|}_{L^2(\Omega_{R})}^2
+\left\|\overline{\theta}_{0}\right\|_{L^2(\Omega_{R})}^2\right)\leq \widetilde{R}^2.$ We can define the following map,
\begin{align*}
    \mathbf{F}^{m}:~B_{\widetilde{R}\mathbf{Y}_{m}}\rightarrow B_{\widetilde{R}\mathbf{Y}_{m}},
\end{align*}
where $\mathbf{F^{m}}(\mathbf{f^{m}})=\widetilde{\mathbf{f}}^m.$

\begin{lemma}\label{youjieyinli} Let $e^{gT}\left(\big{\|}\mathbf{u}_{0}\big{\|}_{L^2(\Omega_{R})}^2
    +\left\|\overline{\theta}_{0}\right\|_{L^2(\Omega_{R})}^2\right)\leq \widetilde{R}^2$.
  Then  $\mathbf{F^{m}}:~B_{\widetilde{R}\mathbf{Y}_{m}}\rightarrow B_{\widetilde{R}\mathbf{Y}_{m}}$
    is continuous and
    compact.
\end{lemma}
\begin{proof}
    We first prove the continuity of $\mathbf{F}^{m}$.
    Taking $k\rightarrow +\infty$, one has $\mathbf{f_{k}^{m}}(t)=
    \left(f_{m1}^{k}(t),\cdots,f_{mm}^{k}(t)\right)$ $\rightarrow \left(f_{m1}(t),\cdots, f_{mm}(t)\right)=\mathbf{f^{m}}(t)
    $ in $\mathbf{Y_{m}}$. And we define
$\mathbf{U}_{k}^{m}=\sum\limits_{i=1}^{m}f_{mi}^{k}(t)\mathbf{e}^{i}$ and
$\mathbf{U}^{m}=\sum\limits_{i=1}^{m}f_{mi}(t)\mathbf{e}^{i}$. Obviously, $\mathbf{U}_{k}^{m}\rightarrow \mathbf{U}^{m}$ in $C\left(0,T;\left(C^{1}(\Omega_{R})\right)^2\right)$.
    From Lemma \ref{youjiedier}, one has as $k\rightarrow +\infty$,
    $\theta_{k}^{m}\rightarrow \theta^m$ in $C\left([0,T]\times \overline{\Omega}_{R}\right)$,
    where $\theta_{k}^m$ and $\theta^m$ solve the following equations, respectively,
    \begin{align*}
        \begin{cases}
            \partial_{t}\theta_{k}^{m}+(\mathbf{U}_{k}^{m}\cdot\nabla)\theta_{k}^{m}
            =0,
            \\
            \theta^{m}_{k}(0,x)=\overline{\theta}_{0},
        \end{cases}~\text{and}~~~~
        \begin{cases}
            \partial_{t}\theta^{m}+(\mathbf{U}^{m}\cdot\nabla)\theta^{m}=0,
            \\
            \theta(0,x)=\overline{\theta}_{0}.
        \end{cases}
\end{align*}
Since $\theta_{k}^{m}$ and $\theta^{m}$ are given,  one has from \eqref{changweifen} that $\widetilde{\mathbf{f_{k}^{m}}}
=\left(\widetilde{f}_{m1}^{k},\cdots,\widetilde{f}_{mm}^{k}(t)\right)$ and
$\widetilde{\mathbf{f^{m}}}=\left(\widetilde{f}_{m1},\cdots,\widetilde{f}_{mm}(t)\right)$. Let
$\mathbf{u}_{k}^{m}=\sum\limits_{i=1}^{m}\widetilde{f}_{mi}^{k}(t)\mathbf{e}^{i}$ and
$\mathbf{u}^{m}=\sum\limits_{i=1}^{m}\widetilde{f}_{mi}(t)\mathbf{e}^{i}$. We have
\begin{align*}
    \begin{aligned}
        &\int_{\Omega_{R}}\partial_{t}\mathbf{u}_{k}^{m}\mathbf{e}^{i}dx
        +\int_{\Omega_{R}}(\mathbf{u}_{k}^{m}\cdot\nabla)\mathbf{u}_{k}^{m}\cdot\mathbf{e}^i dx
        -\mu\int_{\Omega_{R}}\Delta \mathbf{u}_{k}^{m}\cdot\mathbf{e}^{i}dx
        -\int_{\Omega_{R}}\theta_{k}^{m}\mathbf{g}\cdot\mathbf{e}^{i}dx=0,
        \\
        &\int_{\Omega_{R}}\partial_{t}\mathbf{u}^{m}\mathbf{e}^{i}dx
        +\int_{\Omega_{R}}(\mathbf{u}^{m}\cdot\nabla)\mathbf{u}^{m}\cdot\mathbf{e}^i dx
        -\mu\int_{\Omega_{R}}\Delta \mathbf{u}^{m}\cdot\mathbf{e}^{i}dx
        -\int_{\Omega_{R}}\theta^{m}\mathbf{g}\cdot\mathbf{e}^{i}dx=0,
        \\
        &\mathbf{u}_{k}^{m}(0,x)=\mathbf{u}^{m}(0,x).
    \end{aligned}
\end{align*}
Therefore, a series of calculations yields that
\begin{align*}
    \begin{aligned}
    \frac{1}{2}\frac{d}{dt}\big{\|}\mathbf{u}_{k}^{m}-\mathbf{u}^{m}\big{\|}_{L^2(\Omega_{R})}^2
    &\leq \int_{\Omega_{R}}(\theta_{k}^{m}-\theta^{m})\mathbf{g}\cdot
    (\mathbf{u}_{k}^{m}-\mathbf{u}^{m}) dx
    +\int_{\Omega_{R}}\left[((\mathbf{u}_{k}^{m}-\mathbf{u}^{m})\cdot\nabla)\mathbf{u}_{k}^{m}
    (\mathbf{u}_{k}^{m}-\mathbf{u}^{m})\right]dx
    \\
    ~~~~&\leq \frac{1}{2}|\mathbf{g}|\left(\big{\|}\theta_{k}^{m}-\theta^{m}\big{\|}_{L^2(\Omega_{R})}^2
    +\big{\|}\mathbf{u}_{k}^{m}-\mathbf{u}^{m}\big{\|}_{L^2(\Omega_{R})}^2\right)
    +\big{\|}\nabla\mathbf{u}_{k}^{m}\big{\|}_{L^{\infty(\Omega_{R})}}\big{\|}\mathbf{u}_{k}^{m}-\mathbf{u}^{m}\big{\|}_{L^2(\Omega_{R})}^2.
    \end{aligned}
\end{align*}
Since $\big{\|}\mathbf{u}_{k}^{m}\big{\|}_{L^2(\Omega_{R})}^2=\sum\limits_{i=1}^{m}\left|\widetilde{f}_{mi}^{k}(t)\right|^2$
is uniformly bounded on $[0,T]$
, then there exists a constant $M(m)$ such that
\begin{align*}
    \big{\|}\nabla\mathbf{u}_{k}^{m}\big{\|}_{L^\infty(\Omega_{R})}=\sup\limits_{x\in\overline{\Omega}_{R}}
    \left|\sum\limits_{i}^{m}\widetilde{f}_{mi}^{k}(t)\nabla\mathbf{e}^{i}\right|\leq M(m).
\end{align*}
Thus, one has
\begin{align*}
    \frac{d}{dt}\left[e^{-(2M(m)+|\mathbf{g}|)t}
    \big{\|}\mathbf{u}_{k}^{m}-\mathbf{u}^m\big{\|}_{L^2(\Omega_{R})}^2\right]\leq e^{-(2M(m)+|\mathbf{g}|)t}
    \big{\|}\theta_{k}^{m}-\theta^m\big{\|}_{L^2(\Omega_{R})}^2.
\end{align*}
Then, there appears the relation
\begin{align*}
    \sum\limits_{i=1}^{m}\left|\widetilde{f}_{mi}^{k}(t)-\widetilde{f}_{mi}(t)\right|^2
    =\big{\|}\mathbf{u}_{k}^{m}-\mathbf{u}^{m}\big{\|}_{L^2(\Omega_{R})}^2\leq e^{(2M(m)+\mathbf{g})T}
    \int_{0}^{T}\big{\|}\theta_{k}^{m}-\theta^m\big{\|}_{L^2(\Omega_{R})}^2ds
    \rightarrow 0,~~\text{as~}k\rightarrow+\infty,
\end{align*}
which verifies that $\mathbf{F^{m}}$ is continuous.

Moreover, from \eqref{changweifen} again, one can conclude that
$\left\{\frac{d}{dt}\widetilde{f}_{k}^{m}(t)\right\}_{k=1}^{\infty}$ is uniformly bounded. Thus,
by the Arzela-Acosli theorem, $\left\{\widetilde{f}_{k}^{m}\right\}_{k=1}^{\infty}$
has a convergent subsequence. Thus, $\mathbf{F^{m}}$ is compact.
\end{proof}
\begin{remark}\label{bijinjie0402}By the Schauder fixed point theorem, then for each $m\in N^*$,
    there exists a $\mathbf{f^{m}}(t)
    =(f_{m1}(t),\cdots,$ $f_{mm}(t))\in B_{\widetilde{R}\mathbf{Y_{m}}}$ such that
    $\mathbf{F^m}(\mathbf{f^{m}}(t))=\mathbf{f^{m}}(t).$ Let $\mathbf{u^{m}}=\sum\limits_{i}^{m}
    f_{mi}(t)\mathbf{e}^{i}$. Then there exists a $\theta^{m}$ such that $\mathbf{u}^m$ and $\theta^m$
    satisfying the following equation,
    \begin{align}\label{bijinjie}
        \begin{aligned}
            &\int_{\Omega_{R}}\partial_{t}\mathbf{u}^{m}\cdot\mathbf{e}^{i}dx
            +\int_{\Omega_{R}}(\mathbf{u}^{m}\cdot\nabla)\mathbf{u}^{m}\cdot \mathbf{e}^{i}dx
            =\int_{\Omega_{R}}\mu\Delta\mathbf{u}^{m}\cdot \mathbf{e}^{i}
            +\theta^{m}\mathbf{g}\cdot\mathbf{e}^{i}dx,
            \\
            &\partial_{t}\theta^{m}+(\mathbf{u^{m}}\cdot\nabla)\theta^{m}=0,
            \\
            &\mathbf{u}^{m}(0,x)=\sum\limits_{i=1}^{m}\left(\mathbf{u}_{0},
            \mathbf{e}^{i}\right)\mathbf{e}^{i}
            ,~~\theta^m(0,x)=\overline{\theta}_{0}.
        \end{aligned}
    \end{align}
   We conclude that $\left\{(\theta^m,\mathbf{u}^m)\right\}_{m=1}^{\infty}$ is called the approximate solution sequence.
\end{remark}
\begin{proof}\textbf{Proof of Theorem \ref{cunzairuojie}}.

    To show the existence of weak solution, we shall derive a prior estimate for approximate solution of $\{\mathbf{u}^{m},\theta^{m}\}$.

    \textbf{Step 1.} Estimate of $\theta^m$

    Multiplying $\eqref{bijinjie}_{2}$ by $(\theta^m)^{k-1}$$(k>1)$ and then integrating on $\Omega_{R}$,
     one has from the incompressible condition $\eqref{model}_{3}$ that
    \begin{align}\label{woyaoni2}
        \big{\|}\theta^m\big{\|}_{L^k(\Omega_{R})}=\left\|\overline{\theta}_{0}\right\|_{L^k(\Omega_{R})},
    \end{align}
    in which the right side is independent of $m$. Moreover, $\overline{\theta}_{0}$ is the mollification of
    $\theta_{0}$,
    which means that $\{\theta^m\}$ is bounded in $L^\infty\left([0,T];L^\infty(\Omega_{R})\right)$.

   Hence, there exist a subsequence $\{\theta^{m_{k}}\}$ of $\{\theta^m\}$ and
    $\theta\in L^\infty([0,T];L^\infty(\Omega_{R}))$ such that
    \begin{align}\label{shoulian1}
        \theta^{m_{k}}\rightarrow \theta ~\text{in~} L^\infty([0,T];L^\infty(\Omega_{R})
       )
        ~\text{weakly~star~ as}~k\rightarrow\infty.
    \end{align}

    \textbf{Step 2.} Estimate of $\mathbf{u}^m$

 Replace $\mathbf{e}^{i}$ by $\mathbf{u}^m$ in $\eqref{bijinjie}_{1}$ and
    apply $k_{i}\leq 0$ $(i=0,1)$, then there appears the relation
    \begin{align}\label{jijifen}
        \frac{d}{dt}\big{\|}\mathbf{u}^m\big{\|}_{L^2(\Omega_{R})}^2+2\mu\big{\|}\nabla\mathbf{u}^m\big{\|}_{L^2(\Omega_{R})}^2
        \leq |\mathbf{g}|\left(\big{\|}\mathbf{u}^m\big{\|}_{L^2}^2+\big{\|}\theta^{m}\big{\|}_{L^2(\Omega_{R})}^2\right)
        \leq |\mathbf{g}|\left(\big{\|}\mathbf{u}^m\big{\|}_{L^2(\Omega_{R})}^2+\left\|\overline{\theta}_{0}\right\|_{L^2(\Omega_{R})}^2\right).
    \end{align}
    Since $\overline{\theta}_{0}$ is the mollification of $\theta_{0}$, we have
    \begin{align}
        \frac{d}{dt}\big{\|}\mathbf{u}^m\big{\|}_{L^2(\Omega_{R})}^2+2\mu\big{\|}\nabla\mathbf{u}^m\big{\|}_{L^2(\Omega_{R})}^2\leq
        |\mathbf{g}|\left(\big{\|}\mathbf{u}^m\big{\|}_{L^2(\Omega_{R})}^2+\left\|\theta_{0}\right\|_{L^2(\Omega_{R})}^2\right).
    \end{align}
    From the Gronwall's inequality, one has
\begin{align}\label{yongdedaoha1228}
    \big{\|}\mathbf{u}^m\big{\|}_{L^2(\Omega_{R})}^2\leq \left(\big{\|}\theta_{0}\big{\|}_{L^2(\Omega_{R})}^2+\big{\|}\mathbf{u}_{0}\big{\|}_{L^2(\Omega_{R})}^2\right)e^{|
        \mathbf{g}
    |t}.
\end{align}
Moreover, integrate \eqref{jijifen} on $[0,T]$, one obtains
\begin{align}\label{yongdedaohaha1228}
   2\mu\int_{0}^{T}\big{\|}\nabla\mathbf{u}^m\big{\|}^{2}_{L^2(\Omega_{R})}ds
    \leq
    \big{\|}\mathbf{u}_{0}\big{\|}_{L^2(\Omega_{R})}^2+
   \left(\big{\|}\theta_{0}\big{\|}_{L^2(\Omega_{R})}^2+\big{\|}\mathbf{u}_{0}\big{\|}_{L^2(\Omega_{R})}^2\right)e^{|\mathbf{g}|T}
   +|\mathbf{g}|\big{\|}\theta_{0}\big{\|}_{L^2(\Omega_{R})}^2 T.
\end{align}
By \eqref{yongdedaoha1228} and \eqref{yongdedaohaha1228}, we have
\begin{align}\label{shoulian2}
    \{\mathbf{u}^m\}\text{~is~bounded~in~} L^\infty\left(0,T;\widehat{L}^2(\Omega_{R})\right)
    ~\text{and}~L^{2}\left(0,T;\mathbf{G}^{1,2}(\Omega_{R})\right).
\end{align}
Then there exist a $\mathbf{u}\in L^\infty\left(0,T;\widehat{L}^2(\Omega_{R})\right)
\cap L^{2}\left(0,T;\mathbf{G}^{1,2}(\Omega_{R})\right)$ and a subsequence $\{\mathbf{u}^{m_{k}}\}$ such that as
$k\rightarrow\infty$,
\begin{align}
    \mathbf{u}^{m_{k}}\rightarrow \mathbf{u}~\text{weakly~in~}L^{2}\left(0,T;\mathbf{G}^{1,2}(\Omega_{R})\right)~ \text{and} ~
    \text{weakly~star~in~}L^\infty\left(0,T;\widehat{L}^2(\Omega_{R})\right).
\end{align}
By the same steps in the P287 of \cite{Temam1977}, we can also prove
\begin{align}\label{shoulian3}
    \mathbf{u}^{m_{k}}\rightarrow \mathbf{u}~\text{strongly~in~} L^{2}\left(0,T;\widehat{L}^2(\Omega_{R})\right)
    ~\text{as}~k\rightarrow\infty.
\end{align}
The estimates \eqref{shoulian1}, \eqref{shoulian2} and \eqref{shoulian3} enable us to take the limit in the equality \eqref{bijinjie}.
 Let $\overline{\Phi}(t)=\left(\overline{\Phi}_{1}(t),\overline{\Phi}_{2}(t)\right)$ be continuous and differential on $[0,T]$ and
 $\Psi(t,x)\in C^{1}\left(0,T;H^{1}(\Omega_{R})\right)$. Moreover, we take
$\overline{\Phi}(T)=0$ and $\Psi(T)=0$.
Multiplying $\eqref{bijinjie}_{1}$ and $\eqref{bijinjie}_{2}$ by $\Phi(t)$ and  $\Psi(t)$, then integrating obtained equalities on $[0,T]$ give that
\begin{align*}
    \begin{aligned}
    &\begin{aligned}
        &-\int_{0}^{T}\int_{\Omega_{R}}\mathbf{u}^{m_k}\cdot
        \partial_{t}\overline{\Phi}\mathbf{e^{i}} dxdt
        -\int_{0}^{T}\int_{\Omega_{R}}(\mathbf{u}^{m_k}\cdot\nabla)\overline{\Phi}\mathbf{e^i}\cdot
        \mathbf{u}^{m_k}
        -\mu\nabla\mathbf{u}^{m_k}\cdot \overline{\Phi}\nabla\mathbf{e^i}
        dxdt
        \\
        &~-\int_{0}^{T}\int_{-R}^{R}\overline{\Phi}(t)\left[k_{1}u^{{m}_{k}}_{1}(x_{1},1)
        e^{i}_{1}(x_{1},1)
        +k_{0}u^{m_k}_{1}(x_{1},0)
        e^{i}_{1}(x_{1},0)\right]dx_{1}dt
        \\
        &~=\int_{\Omega_{R}}\mathbf{u}^{m_k}_{0}\overline{\Phi}(0)
        \mathbf{e^{i}} dx+\int_{0}^{T}\int_{\Omega_{R}}\theta^{m_k} \mathbf{g}
        \cdot \overline{\Phi} \mathbf{e^{i}}dxdt,
    \end{aligned}
    \\
    &\begin{aligned}
        -\int_{0}^{T}\int_{\Omega_{R}}\theta^{m_k}\partial_{t}\Psi dxdt
        -\int_{0}^{T}\int_{\Omega_{R}}(\mathbf{u}^{m_k}\cdot\nabla)\Psi\theta^{m_k} dxdt
        =-\int_{\Omega_{R}}\overline{\theta}_{0}\Psi(0,x)dx,
     \end{aligned}
\end{aligned}
\end{align*}
where $\overline{\Phi}\mathbf{e^{i}}:=\left(\overline{\Phi}_{1}(t)e^{i}_{1},
\overline{\Phi}_{2}(t)e^{i}_{2}\right)$.

From \eqref{shoulian3}, one has for $i=1,2,\cdots$,
\begin{align}\label{youxianxiangjia0102}
    \int_{0}^{T}\int_{\Omega_{R}}\left(\mathbf{u}^{m_k}\cdot\nabla\right)\overline{\Phi}\mathbf{e^i}\cdot
        \mathbf{u}^{m_k}
        dxdt\rightarrow
        \int_{0}^{T}\int_{\Omega_{R}}\left(\mathbf{u}\cdot\nabla\right)\overline{\Phi}\mathbf{e^i}\cdot
        \mathbf{u}
        dxdt,
\end{align}
as $m_{k}\rightarrow\infty$. Here,  the above equality holds for any finite linear combination of $\{\mathbf{e}^{i}\}$.
For any $\mathbf{v}\in \mathbf{G}^{1,2}(\Omega_{R})$, we shall show that  $m_{k}\rightarrow\infty$,
\begin{align}\label{zhengdejiushini0102}
    \int_{0}^{T}\int_{\Omega_{R}}(\mathbf{u}^{m_k}\cdot\nabla)\overline{\Phi}\mathbf{v}\cdot
        \mathbf{u}^{m_k}
        dxdt\rightarrow
        \int_{0}^{T}\int_{\Omega_{R}}(\mathbf{u}\cdot\nabla)\overline{\Phi}\mathbf{v}\cdot
        \mathbf{u}
        dxdt.
\end{align}
In fact, since $\mathbf{v}\in \mathbf{G}^{1,2}(\Omega_{R})$,
there exists a $\widetilde{m}\in \mathbf{N}$ and $\forall~\varepsilon>0$
such that
\begin{align}\label{haiyouni0102}
    \left\|\mathbf{v}-\sum\limits_{i=1}^{\widetilde{m}}v_{i}\mathbf{e}^{i}\right\|_{\mathbf{G}^{1,2}(\Omega_{R})}<\varepsilon.
\end{align}
Therefore, we have
\begin{align*}
    \begin{aligned}
        \int_{0}^{T}\int_{\Omega_{R}}&\left[(\mathbf{u}^{m_{k}}\cdot\nabla)\overline{\Phi}\mathbf{v}\mathbf{u}^{m_{k}}
        -(\mathbf{u}\cdot\nabla)\overline{\Phi}\mathbf{v}\mathbf{u}\right]dxdt
        =\int_{0}^{T}\int_{\Omega_{R}}\left[(\mathbf{u}^{m_{k}}\cdot\nabla)\overline{\Phi}
        (\mathbf{v}-\sum\limits_{i=1}^{\widetilde{m}}v_{i}\mathbf{e}^{i})\mathbf{u}^{m_{k}}\right]dxdt
        \\
        &+\int_{0}^{T}\int_{\Omega_{R}}\left[(\mathbf{u}^{m_{k}}\cdot\nabla)\overline{\Phi}
        \sum\limits_{i=1}^{\widetilde{m}}v_{i}\mathbf{e}^{i}\mathbf{u}^{m_{k}}-
        (\mathbf{u}\cdot\nabla)\overline{\Phi}
        \sum\limits_{i=1}^{\widetilde{m}}v_{i}\mathbf{e}^{i}\mathbf{u}\right]dxdt
        \\
        &+\int_{0}^{T}\int_{\Omega_{R}}\left[(\mathbf{u}\cdot\nabla)\overline{\Phi}
        \sum\limits_{i=1}^{\widetilde{m}}v_{i}\mathbf{e}^{i}\mathbf{u}-
        (\mathbf{u}\cdot\nabla)\overline{\Phi}
        \mathbf{v}\mathbf{u}\right]dxdt.
    \end{aligned}
\end{align*}
Thus, from \eqref{shoulian2}, \eqref{youxianxiangjia0102} and \eqref{haiyouni0102},
we can prove that \eqref{zhengdejiushini0102} is valid.

For the nonlinear term $\int_{0}^{T}\int_{\Omega_{R}}(\mathbf{u}^{m_k}\cdot\nabla)\Psi\theta^{m_k} dxdt$, one can easily obtain the convergence due to $\eqref{shoulian1}$.
Thus, take $m_{k}\rightarrow\infty$, one has
\begin{align}\label{dingyiruojie}
    \begin{aligned}
    &\begin{aligned}
        &-\int_{0}^{T}\int_{\Omega_{R}}\mathbf{u}\cdot
        \partial_{t}\overline{\Phi}\mathbf{e^{i}} dxdt
        -\int_{0}^{T}\int_{\Omega_{R}}(\mathbf{u}\cdot\nabla)\overline{\Phi}\mathbf{e^i}\cdot
        \mathbf{u}
        -\mu\nabla\mathbf{u}\cdot \overline{\Phi}\nabla\mathbf{e^i}
        dxdt
        \\
        &~-\int_{0}^{T}\int_{-R}^{R}\overline{\Phi}(t)\left[k_{1}u_{1}(t,x_{1},1)
        e^{i}_{1}(x_{1},1)
        +k_{0}u_{1}(t,x_{1},0)
        e^{i}_{1}(x_{1},0)\right]dx_{1}dt
        \\
        &~=\int_{\Omega_{R}}\mathbf{u}_{0}\overline{\Phi}(0)
        \mathbf{e^{i}} dx+\int_{0}^{T}\int_{\Omega_{R}}\theta \mathbf{g}
        \cdot \overline{\Phi} \mathbf{e^{i}}dxdt,
    \end{aligned}
    \\
    &\begin{aligned}
        -\int_{0}^{T}\int_{\Omega_{R}}\theta\partial_{t}\Psi dxdt
        -\int_{0}^{T}\int_{\Omega_{R}}(\mathbf{u}\cdot\nabla)\Psi\theta dxdt
        =\int_{\Omega_{R}}\theta_{0}\Psi(0,x)dx.
     \end{aligned}
\end{aligned}
\end{align}
Since $\left\{\mathbf{e}^{n}\right\}_{n=1}^{\infty}$ is dense in $\mathbf{G}^{1,2}(\Omega_{R})$, then $(\mathbf{u},
\theta)$ meets \eqref{ruojie}.

The rest part follows the same steps in \cite{Temam1977}, then it follows that  $\mathbf{u}(0,x)=\mathbf{u}_{0}$ and
$\theta(0,x)=\theta_{0}$.
Thus, we end the proof of Theorem \ref{cunzairuojie}.
\end{proof}
\begin{proof}\textbf{Proof of Theorem \ref{bufenqiangjie1227}.}

    \textbf{Step 1.} Estimate of $\partial_{t}\mathbf{u}$

    Since $\mathbf{u}_{0}\in \mathbf{V}^{1,2}(\Omega_{R})$, then
    $\mathbf{u}^{m}(0,x)$ can be taken as the projection of $\mathbf{u}_{0}$ on the space spanned by
    $\mathbf{e}^{1},\cdots,\mathbf{e}^{m}$. Thus, when $m\rightarrow\infty$, one has
    \begin{align}\label{youjie321}
       \mathbf{u}^{m}(0,x)\rightarrow\mathbf{u}_{0}~~\text{strongly~in}~\left[H^2(\Omega_{R})\right]^2,~
       \big{\|}\mathbf{u}^{m}(0,x)\big{\|}_{H^2(\Omega_{R})}\leq\big{\|}\mathbf{u}_{0}\big{\|}_{H^2(\Omega_{R})}.
    \end{align}
    To derive an estimate of $\partial_{t}\mathbf{u}^{m}$,
    we replace $\mathbf{e}^{i}$ by $\partial_{t}\mathbf{u}^{m}$  in $\eqref{bijinjie}_{1}$ and then get
    \begin{align*}
       \begin{aligned}
       \big{\|}\partial_{t}&\mathbf{u}^{m}\big{\|}_{L^2(\Omega_{R})}^{2}
       =\mu\int_{\Omega_{R}}\Delta\mathbf{u}^{m}\cdot\partial_{t}\mathbf{u}^{m}dx
       +\int_{\Omega_{R}}\theta^{m}\mathbf{g}\cdot\partial_{t}\mathbf{u}^m dx
       -\int_{\Omega_{R}}(\mathbf{u}^m\cdot\nabla)\mathbf{u}^m\cdot\partial_{t}\mathbf{u}^m dx
       \\
       &\leq C\big{\|}\Delta\mathbf{u}^{m}\big{\|}_{L^2(\Omega_{R})}^2+
       \big{\|}\theta^{m}\mathbf{g}\big{\|}_{L^2(\Omega_{R})}^2
       +\frac{3}{4}\big{\|}\partial_{t}\mathbf{u}^{m}\big{\|}_{L^2(\Omega_{R})}^2
       +\big{\|}\mathbf{u}^m\big{\|}_{L^{4}(\Omega_{R})}^{2}\big{\|}\nabla\mathbf{u}^{m}\big{\|}_{L^4(\Omega_{R})}^{2}
       \\
       &\leq C\big{\|}\Delta\mathbf{u}^{m}\big{\|}_{L^2(\Omega_{R})}^2+
       \left\|\overline{\theta}_{0}\mathbf{g}\right\|_{L^2(\Omega_{R})}^2
       +\frac{3}{4}\big{\|}\partial_{t}\mathbf{u}^{m}\big{\|}_{L^2(\Omega_{R})}^2
       +C\big{\|}\mathbf{u}^{m}\big{\|}_{H^2(\Omega_{R})}^{4},
       \end{aligned}
    \end{align*}
where the Lemmas \ref{sobolveqianrudingli12241} and \ref{jiushizheyang12242} are used.
Let $t\rightarrow 0^{+}$ in above inequality. We have from \eqref{youjie321} that
       \begin{align}\label{20240225}
           \big{\|}\partial_{t}\mathbf{u}^{m}(0,x)\big{\|}_{L^2(\Omega_{R})}^2
           \leq C(\mu)\big{\|}\mathbf{u}_{0}\big{\|}_{H^2(\Omega_{R})}^2
           +C\big{\|}\mathbf{u}_{0}\big{\|}_{H^2(\Omega_{R})}^4
           +\big{\|}\theta_{0}\mathbf{g}\big{\|}_{L^2(\Omega_{R})}^{2}.
       \end{align}
       Differentiating $\eqref{bijinjie}_{1}$ with respect to $t$, and by virtue of
       $\eqref{bijinjie}_{2}$, one has
       \begin{align}\label{xinxiangshicheng}
           \begin{aligned}
               &\int_{\Omega_{R}}\partial_{tt}\mathbf{u}^{m}\mathbf{e}^{i}dx
               -\mu\int_{\Omega_{R}}\partial_{t}\Delta\mathbf{u}^{m}\cdot \mathbf{e}^{i}dx
               \\
               &=
               -\int_{\Omega_{R}}(\partial_{t}\mathbf{u}^{m}\cdot\nabla)\mathbf{u}^{m}
               \cdot \mathbf{e}^{i}dx-\int_{\Omega_{R}}(\mathbf{u}^{m}\cdot\nabla)
               \partial_{t}\mathbf{u}^{m}
               \cdot \mathbf{e}^{i}
               -\partial_{t}\theta^{m}\mathbf{g}
               \cdot\mathbf{e}^{i}dx,
               \\
               &=-\int_{\Omega_{R}}(\partial_{t}\mathbf{u}^{m}\cdot\nabla)\mathbf{u}^{m}
               \cdot \mathbf{e}^{i}dx-\int_{\Omega_{R}}(\mathbf{u}^{m}\cdot\nabla)
               \partial_{t}\mathbf{u}^{m}
               \cdot \mathbf{e}^{i}+
               (\mathbf{u}^m\cdot\nabla)\theta^{m}\mathbf{g}
               \cdot\mathbf{e}^{i}dx.
           \end{aligned}
       \end{align}
       Multiplying  \eqref{xinxiangshicheng} by $f_{mi}^{'}(t)$ ($i=1,\cdots,m$), and then summing up all  results gives that
\begin{align}\label{fabiao}
   \begin{aligned}
       &\frac{1}{2}\frac{d}{dt}\big{\|}\partial_{t}\mathbf{u}^{m}\big{\|}_{L^2(\Omega_{R})}^{2}
       +\mu\big{\|}\nabla\partial_{t}\mathbf{u}^{m}\big{\|}_{L^2(\Omega_{R})}^2
       \\
       &~\leq \int_{\Omega_{R}}(\mathbf{u}^{m}\cdot\nabla)\partial_{t}\mathbf{u}^{m}\cdot\theta^{m}
       \mathbf{g}dx
       -\int_{\Omega_{R}}(\partial_{t}\mathbf{u}^{m}\cdot\nabla)\mathbf{u}^{m}\cdot
       \partial_{t}\mathbf{u}^{m}dx
       \\
       &~\leq \big{\|}\theta^{m}\big{\|}_{L^{\infty}}\big{\|}\mathbf{u}^{m}\big{\|}_{L^2(\Omega_{R})}
       \big{\|}\nabla\partial_{t}\mathbf{u}^{m}\big{\|}_{L^2(\Omega_{R})}
     +\big{\|}\nabla\mathbf{u}^{m}\big{\|}_{L^2(\Omega_{R})}
       \big{\|}\partial_{t}\mathbf{u}^{m}\big{\|}_{L^4(\Omega_{R})}^{2}.
   \end{aligned}
\end{align}
From Lemma \ref{L4guji}, we know
\begin{align*}
   \big{\|}\partial_{t}\mathbf{u}^{m}\big{\|}_{L^4(\Omega_{R})}^2
   \leq C\big{\|}\nabla\partial_{t}\mathbf{u}^{m}\big{\|}_{L^2(\Omega_{R})}
   \big{\|}\partial_{t}\mathbf{u}^{m}\big{\|}_{L^2(\Omega_{R})}
   ,
\end{align*}
which along with \eqref{fabiao} yields that
\begin{align*}
   \begin{aligned}
       \frac{1}{2}\frac{d}{dt}\big{\|}\partial_{t}\mathbf{u}^m\big{\|}_{L^2(\Omega_{R})}^2
       +\frac{\mu}{2}\big{\|}\nabla\partial_{t}\mathbf{u}^{m}\big{\|}_{L^2(\Omega_{R})}^{2}
       &\leq
       C\big{\|}\nabla\mathbf{u}^{m}\big{\|}_{L^2(\Omega_{R})}^{2}
       \big{\|}\partial_{t}\mathbf{u}^{m}\big{\|}_{L^2(\Omega_{R})}^2
+C\big{\|}\theta^{m}\big{\|}_{L^{\infty}(\Omega_{R})}^{2}\big{\|}\mathbf{u}^{m}\big{\|}_{L^2(\Omega_{R})}^{2}.
   \end{aligned}
\end{align*}
Applying Gronwall's inequality and \eqref{20240225}, one has
\begin{align}\label{woyaoni1}
   \begin{aligned}
       \big{\|}\partial_{t}\mathbf{u}^{m}\big{\|}_{L^2(\Omega_{R})}^{2}<C,~
      \int_{0}^{T} \big{\|}\nabla\partial_{s}\mathbf{u}^{m}\big{\|}_{L^2(\Omega_{R})}^{2}ds<C,
   \end{aligned}
\end{align}
which implies
$$\partial_{t}\mathbf{u}^{m}\text{~is~ bounded~in~}L^{\infty}\left(0,T;\widehat{L}^2(\Omega_{R})\right)
\cap L^2\left(0,T;\mathbf{G}^{1,2}(\Omega_{R})\right).$$
Thus, $\partial_{t}\mathbf{u}\in L^{\infty}\left(0,T;\widehat{L}^2(\Omega_{R})\right)
\cap L^2\left(0,T;\mathbf{G}^{1,2}(\Omega_{R})\right)$.

\textbf{Step 2.} Estimates of $\mathbf{u}$, $\nabla\mathbf{u}$, $\nabla^2\mathbf{u}$ and $\nabla p$

Since $\partial_{t}\mathbf{u}\in L^{2}\left(0,T;\mathbf{G}^{1,2}(\Omega_{R})\right)$ and
$\mathbf{u}\in L^{2}\left(0,T;\mathbf{G}^{1,2}(\Omega_{R})\right)$, thus from Theorem 2 in P286 of \cite{noauthor_partial_nodate}, we can deduce that
\begin{align}\label{3.65}
\mathbf{u}\in C\left(0,T;\mathbf{G}^{1,2}(\Omega_{R})\right)~\text{and~}
\big{\|}\nabla\mathbf{u}\big{\|}_{L^2(\Omega_{R})}^2\leq C\int_{0}^{T}\big{\|}\mathbf{u}\big{\|}_{H^{1}(\Omega_{R})}^2
+\big{\|}\partial_{t}\mathbf{u}\big{\|}_{H^1(\Omega_{R})}^2ds<C.
\end{align}
One can easily check that $-\partial_{t}\mathbf{u}-(\mathbf{u}\cdot\nabla)\mathbf{u}
+\theta\mathbf{g}\in L^{\frac{3}{2}}(\Omega_{R})$ by virtue of holder's inequality. Moreover, the definition of weak solution and the existence of $\partial_{t} \mathbf{u}$ allow us to get that
\begin{align*}
    \begin{aligned}
    \mu\int_{\Omega_{R}}\nabla\mathbf{u}\cdot\nabla \mathbf{e}^{i}dx
    &-\int_{-R}^{R}\left[k_{1}u_{1}(x_{1},1)\mathbf{e}^{i}_{1}(x_{1},1)+
    k_{0}u_{1}(x_{1},0)\mathbf{e}^{i}_{1}(x_{1},0)\right]dx_{1}
    \\
    &=\int_{\Omega_{R}}
    \left(-\partial_{t}\mathbf{u}-(\mathbf{u}\cdot\nabla)\mathbf{u}
    +\theta\mathbf{g}\right)\cdot\mathbf{e}^{i}dx.
    \end{aligned}
\end{align*}
Hence, from Lemma \ref{tishengzhengzexing}, one can obtain $\mathbf{u}\in \left[W^{2,\frac{3}{2}}(\Omega_{R})\right]^2$
 and $\nabla p\in
\left[L^{\frac{3}{2}}(\Omega_{R})\right]^2$, which means that $\mathbf{u}\in \left[L^{\infty}(\Omega_{R})\right]^2$. Thus,
$-\partial_{t}\mathbf{u}-(\mathbf{u}\cdot\nabla)\mathbf{u}
+\theta\mathbf{g}\in L^{2}(\Omega_{R})$. Again
from Lemma \ref{tishengzhengzexing}, one has $\mathbf{u}\in \mathbf{V}^{1,2}(\Omega_{R})$, $\nabla p\in [L^{2}(\Omega_{R})]^2$ and
\begin{align}\label{guijidedao}
    \begin{aligned}
    &\left\|\nabla^2\mathbf{u}\right\|_{L^2(\Omega_{R})}^{2}
    +\big{\|}\nabla p\big{\|}_{L^2(\Omega_{R})}^{2}\leq C \big{\|}-\partial_{t}\mathbf{u}-(\mathbf{u}\cdot\nabla)\mathbf{u}
    +\theta\mathbf{g}\big{\|}_{L^2(\Omega_{R})}^{2}+C\big{\|}\mathbf{u}\big{\|}_{W^{1,2}(\Omega_{R})}^{2}
    \\
    &~~~~
    \leq C\left(\big{\|}\partial_{t}\mathbf{u}\big{\|}_{L^{2}(\Omega_{R})}^{2}
    +\big{\|}(\mathbf{u}\cdot\nabla)\mathbf{u}\big{\|}_{L^2(\Omega_{R})}^{2}
    +\big{\|}\theta_{0}\big{\|}_{L^2(\Omega_{R})}^{2}
   +\big{\|}\mathbf{u}\big{\|}_{W^{1,2}(\Omega_{R})}^{2}\right)
   \\
   &~~~~\leq C\left(\big{\|}\partial_{t}\mathbf{u}\big{\|}_{L^2(\Omega_{R})}^{2}+
   \big{\|}\mathbf{u}\big{\|}_{L^{\infty}(\Omega_{R})}^{2}\big{\|}\nabla\mathbf{u}\big{\|}_{L^{2}(\Omega_{R})}^{2}
   +\big{\|}\theta_{0}\big{\|}_{L^{2}(\Omega_{R})}^{2}+\big{\|}\nabla\mathbf{u}\big{\|}_{L^2(\Omega_{R})}^2\right)
   \\
   &~~~~\leq C\left(\big{\|}\partial_{t}\mathbf{u}\big{\|}_{L^2(\Omega_{R})}^{2}+
   \big{\|}\mathbf{u}\big{\|}_{H^{1}(\Omega_{R})}^{3}+\big{\|}\mathbf{u}\big{\|}_{H^{1}(\Omega_{R})}^3
   \big{\|}\nabla^2\mathbf{u}\big{\|}_{L^{2}(\Omega_{R})}
   +\big{\|}\theta_{0}\big{\|}_{L^{2}(\Omega_{R})}^{2}+\big{\|}\nabla\mathbf{u}\big{\|}_{L^2(\Omega_{R})}^2\right),
    \end{aligned}
\end{align}
where we used the conclusion of Lemma \ref{Linftyguji}. Furthermore, by \eqref{3.65} and Cauchy inequality, and then taking the
coefficient of $\big{\|}\nabla^2\mathbf{u}\big{\|}^2$ small enough gives that
\begin{align}\label{369}
    \left\|\nabla^2\mathbf{u}\right\|_{L^2(\Omega_{R})}^2\leq C\text{~and~}
    \big{\|}\nabla p\big{\|}_{L^2(\Omega_{R})}^2\leq C.
\end{align}
Multiply $\eqref{youjiequyu}_{1}$ by $\partial_{t}\mathbf{u}$ and integrate on $\Omega_{R}$, one has
\begin{align*}
    \begin{aligned}
        \big{\|}\partial_{t}\mathbf{u}\big{\|}_{L^2(\Omega_{R})}^{2}
        &+\frac{\mu}{2}\frac{d}{dt}\big{\|}\nabla\mathbf{u}\big{\|}_{L^{2}(\Omega_{R})}^{2}
        -\frac{1}{2}\frac{d}{dt}\int_{-R}^{R}\left[k_{1}|u_{1}(x_{1},1)|^2+k_{0}|u_{1}(x_{1},0)|^{2}\right]dx_{1}
        \\
        &\leq \big{\|}(\mathbf{u}\cdot\nabla)\mathbf{u}\big{\|}_{L^{2}(\Omega_{R})}^{2}+
        g^{2}\big{\|}\theta_{0}\big{\|}_{L^{2}(\Omega_{R})}^{2}+\frac{1}{2}
        \big{\|}\partial_{t}\mathbf{u}\big{\|}_{L^2(\Omega_{R})}^{2},
    \end{aligned}
\end{align*}
which together with Lemma \ref{Linftyguji} and  \eqref{369} yields that
\begin{align*}
    \begin{aligned}
        \big{\|}\partial_{t}\mathbf{u}\big{\|}_{L^2(\Omega_{R})}^{2}
        &+\mu\frac{d}{dt}\big{\|}\nabla\mathbf{u}\big{\|}_{L^{2}(\Omega_{R})}^{2}
        -\frac{d}{dt}\int_{-R}^{R}\left[k_{1}|u_{1}(x_{1},1)|^2+k_{0}|u_{1}(x_{1},0)|^{2}\right]dx_{1}
        \\
        &
       \leq 2\big{\|}\mathbf{u}\big{\|}_{L^{\infty}(\Omega_{R})}^2\big{\|}\nabla\mathbf{u}\big{\|}_{L^2(\Omega_{R})}^2
       +2|\mathbf{g}|^2\big{\|}\theta_{0}\big{\|}_{L^2(\Omega_{R})}^2\leq C.
    \end{aligned}
\end{align*}
 Integrating above inequality from $0$ to $t$ and using \eqref{369} again, it follows that
\begin{align}\label{woyaoni3}
~\int_{0}^{t}\big{\|}\partial_{s}\mathbf{u}(s)\big{\|}_{L^{2}(\Omega_{R})}^{2}ds\leq C,
~\int_{0}^{t}\big{\|}\nabla^2\mathbf{u}\big{\|}_{L^{2}(\Omega_{R})}^{2}ds \leq C.
\end{align}
By Lemmas \ref{jiushizheyang12242}, \ref{tishengzhengzexing}, \ref{L4guji}
and \ref{Linftyguji}, one can show
\begin{align*}
    \begin{aligned}
        \int_{0}^{t}\big{\|}&\nabla^2\mathbf{u}\big{\|}_{L^{4}(\Omega_{R})}^{2}+\big{\|}\nabla p\big{\|}_{L^{4}(\Omega_{R})}^2 ds
        \leq
        C\int_{0}^{t}\big{\|}\partial_{s}\mathbf{u}+(\mathbf{u}\cdot\nabla)\mathbf{u}
        -\theta\mathbf{g}\big{\|}_{L^4(\Omega_{R})}^2ds+C\int_{0}^{t}\big{\|}\mathbf{u}\big{\|}_{W^{1,4}(\Omega_{R})}^{2}ds
        \\
        &\leq C\int_{0}^{t}\left(\big{\|}\partial_{s}\mathbf{u}\big{\|}_{L^{4}(\Omega_{R})}^2
        +\big{\|}(\mathbf{u}\cdot\nabla)\mathbf{u}\big{\|}_{L^{4}(\Omega_{R})}^{2}
        +\big{\|}\theta_{0}\big{\|}_{L^{4}(\Omega_{R})}^{2}\right) ds+C\int_{0}^{t}\big{\|}\mathbf{u}\big{\|}_{W^{1,4}(\Omega_{R})}^{2}ds
        \\
        &
        \begin{aligned}
            \leq C\int_{0}^{t}
            \left(\big{\|}\partial_{s}\mathbf{u}\big{\|}_{L^2(\Omega_{R})}^{2}+
            \big{\|}\nabla\partial_{s}\mathbf{u}\big{\|}_{L^{2}(\Omega_{R})}^2
            +\big{\|}\mathbf{u}\big{\|}_{L^{\infty}(\Omega_{R})}^{2}\big{\|}\mathbf{u}\big{\|}_{H^2(\Omega_{R})}^2
            +\big{\|}\theta_{0}\big{\|}_{L^{4}(\Omega_{R})}^{2}
        \right)ds
        \end{aligned}
        \\
        &
        \leq C\int_{0}^{t}
        \left(\big{\|}\partial_{s}\mathbf{u}\big{\|}_{L^2(\Omega_{R})}^{2}+
        \big{\|}\nabla\partial_{s}\mathbf{u}\big{\|}_{L^{2}(\Omega_{R})}^2
        +\big{\|}\mathbf{u}\big{\|}_{H^2(\Omega_{R})}^4
        +\big{\|}\theta_{0}\big{\|}_{L^{4}(\Omega_{R})}^{2}
        \right)ds.
    \end{aligned}
\end{align*}
Thus, in view of  \eqref{woyaoni2}, \eqref{woyaoni1}, \eqref{369} and \eqref{woyaoni3}, we conclude that
\begin{align}\label{woyaoni6}
    \int_{0}^{t}\big{\|}&\nabla^2\mathbf{u}\big{\|}_{L^{4}(\Omega_{R})}^{2}+\big{\|}\nabla p\big{\|}_{L^{4}(\Omega_{R})}^2 ds
    \leq C.
\end{align}
This completes the proof of this Theorem.
\end{proof}
\begin{proof}\textbf{Proof of Theorem \ref{qiangjie}.}

    For $\mathbf{u}\in \mathbf{V}^{1,2}(\Omega_{R})$,  one has from Lemma \ref{jiushizheyang12242} that
    \begin{align*}
     \big{\|}\nabla\mathbf{u}\big{\|}_{L^{4}(
         \Omega_{R}
     )}\leq C\big{\|}\nabla^2\mathbf{u}\big{\|}_{L^2(\Omega_{R})}
     \leq C,
    \end{align*}
which together with $\int_{0}^{T}\left\|\nabla^2\mathbf{u}(s)\right\|_{L^4(\Omega_{R})}^2ds\leq C$ gives that
 \begin{align*}
     \int_{0}^{T}\big{\|}\nabla\mathbf{u}(s)\big{\|}_{W^{1,4}(\Omega_{R})}^2ds\leq C.
 \end{align*}
 From Corollary \ref{ohmygod1228}, one obtains
 \begin{align}\label{yongdedao1228}
     \int_{0}^{T}\big{\|}\nabla\mathbf{u}\big{\|}_{L^{\infty}(\Omega_{R})}^{2}ds\leq C\int_{0}^{T}\big{\|}\nabla\mathbf{u}(s)\big{\|}_{W^{1,4}(\Omega_{R})}^2ds\leq C,
 \end{align}
 implying that $\|\nabla\mathbf{u}\|_{L^\infty(\Omega_{R})}$ is integrable.

 We have shown that $(\mathbf{u},\theta)$ satisfies $\partial_{t}\theta+(\mathbf{u}\cdot\nabla)\theta=0$
 in the weak sense.
Thus, to estimate $\theta$, we assume $\theta$ and $\mathbf{u}$ are smooth.
First, $\theta_{x_{i}}$ $(i=1,2)$ satisfy
\begin{align*}
    \partial_{t}\theta_{x_{i}}+(\partial_i\mathbf{u}\cdot\nabla)\theta=-(\mathbf{u}\cdot\nabla)
    \theta_{x_{i}}.
\end{align*}
Multiply the above equality by $\theta_{x_{i}}$, integrate on $\Omega_{R}$ and add up all
result for $i$, then one gets form \eqref{yongdedao1228} that
\begin{align*}
    \frac{d}{dt}\big{\|}\nabla\theta\big{\|}_{L^2(\Omega_{R})}^2
    \leq C\int_{\Omega_{R}}|\nabla\mathbf{u}||\nabla\theta|^{2}dx\leq C \big{\|}\nabla\mathbf{u}\big{\|}_{L^\infty(\Omega_{R})}\big{\|}\nabla\theta\big{\|}_{L^2(\Omega_{R})}^2.
\end{align*}
Applying Gronwall's inequality, we deduce
\begin{align}\label{mingmingbaibai}
 \big{\|}\nabla\theta\big{\|}_{L^2(\Omega_{R})}^2\leq \big{\|}\nabla\theta_{0}\big{\|}_{L^2(\Omega_{R})}^2
 e^{\int_{0}^{t}\|\nabla\mathbf{u}\|_{L^\infty(\Omega_{R})}ds}
 \leq \big{\|}\nabla\theta_{0}\big{\|}_{L^2(\Omega_{R})}^2
e^{C\int_{0}^{t}\|\nabla\mathbf{u}\|_{W^{1,4}(\Omega_{R})}ds} \leq C.
\end{align}
Next, we estimate $\partial_{t}\theta$. Since
    $\theta_{t}=-(\mathbf{u}\cdot\nabla)\theta$, it is found that
\begin{align}\label{yaodejiushi1228}
    \big{\|}\theta_{t}\big{\|}_{L^2(\Omega_{R})}^2\leq \big{\|}\mathbf{u}\big{\|}_{L^{\infty}(\Omega_{R})}^2
    \big{\|}\nabla\theta\big{\|}_{L^2(\Omega_{R})}^2
    \leq C,
\end{align}
which shows that $\theta_{t}\in L^{\infty}\left(0,T;L^{2}(\Omega_{R})\right)$.

Based on estimates \eqref{yongdedaoha1228}, \eqref{yongdedaohaha1228}, \eqref{woyaoni1},
\eqref{369}, \eqref{woyaoni3}, \eqref{woyaoni6},
\eqref{mingmingbaibai} and \eqref{yaodejiushi1228}, we conclude that  \eqref{kanqilaimaodun} holds.
\end{proof}
\begin{proof}\textbf{Proof of Theorem \ref{weiyixing1227}.}

    Suppose that  $(\mathbf{u},p,\theta)$ and
    $\left(\mathbf{\widetilde{u}},\widetilde{p},\widetilde{\theta}\right)$ are the strong solutions under corresponding initial values $\mathbf{u}_{0}\in
\left[H^{2}(\Omega_{R})\right]^2$, $\theta_{0}\in H^{1}(\Omega_{R})$ and $\widetilde{\mathbf{u}}_{0}\in
\left[H^{2}(\Omega_{R})\right]^2$, $\widetilde{\theta}_{0}\in H^{1}(\Omega_{R})$. Thus, one can obtain
\begin{align*}
    \partial_{t}\left(\mathbf{u}-\mathbf{\widetilde{u}}\right)+
    \nabla\left(p-\widetilde{p}\right)-\mu\Delta\left(\mathbf{u}-\mathbf{\widetilde{u}}\right)=
    \left(\theta-\widetilde{\theta}\right)\overrightarrow{g}+\left(\mathbf{u}\cdot\nabla\right)
    \left(\mathbf{u}-\mathbf{\widetilde{u}}\right)
    +\left(\left(\mathbf{u}-\mathbf{\widetilde{u}}\right)\cdot\nabla\right)\mathbf{\widetilde{u}}.
\end{align*}
Multiply the above equality by $\mathbf{u}-\mathbf{\widetilde{u}}$, integrate on
$\Omega_{R}$, and a direct calculation gives that
\begin{align}\label{weiyixing1}
    \begin{aligned}
    \frac{1}{2}&\frac{d}{dt}\left\|\mathbf{u}-\mathbf{\widetilde{u}}\right\|_{L^{2}(\Omega_{R})}^{2}
    +\mu\left\|\nabla(\mathbf{u}-\mathbf{\widetilde{u}})\right\|_{L^2(\Omega_{R})}^2
    \\
    &\leq \int_{\Omega_{R}}\left(\mathbf{u}-\mathbf{\widetilde{u}}
    \right)\cdot\left(\theta-\widetilde{\theta}\right)
    \overrightarrow{g} dx+\int_{\Omega_{R}}\left((\mathbf{u}-\mathbf{\widetilde{u}})\cdot\nabla\right)
    \widetilde{\mathbf{u}}\cdot(\mathbf{u}-\widetilde{\mathbf{u}})dx
    \\
    &\leq \frac{C\varepsilon}{2}\left(\left\|\mathbf{u}-\mathbf{\widetilde{u}}\right\|_{L^2(\Omega_{R})}^{2}
    +\left\|\nabla(\mathbf{u}-\mathbf{\widetilde{u}})\right\|_{L^2(\Omega_{R})}^2\right)
    +\frac{1}{2\varepsilon}\left\|\theta-\widetilde{\theta}\right\|_{L^{\frac{3}{2}}(\Omega_{R})}^2.
    \end{aligned}
\end{align}
On the other hand, from the equality
\begin{align*}
    \partial_{t}\left(\theta-\widetilde{\theta}\right)=(\mathbf{u}\cdot\nabla)\left(\theta-\widetilde{\theta}\right)
    +\left((\mathbf{u}-\widetilde{\mathbf{u}}\right)\cdot\nabla)\widetilde{\theta},
\end{align*}
there appears the relation
\begin{align*}
    \begin{aligned}
    \frac{d}{dt}\left\|\theta-\widetilde{\theta}\right\|_{L^{\frac{3}{2}}(\Omega_{R})}^{\frac{3}{2}}
    &\leq \frac{3}{2}\int_{\Omega_{R}}\left|\mathbf{u}-\mathbf{\widetilde{u}}\right|\left|\nabla\widetilde{\theta}\right|
    \left|\theta-\widetilde{\theta}\right|^{\frac{1}{2}}dx
    \\
    &\leq \frac{3}{2}\left\|\mathbf{u}-\mathbf{\widetilde{u}}\right\|_{L^{6}(\Omega_{R})}
    \left\|\nabla\widetilde{\theta}\right\|_{L^2(\Omega_{R})}
    \left\|\theta-\widetilde{\theta}\right\|_{L^{\frac{3}{2}}(\Omega_{R})}^{\frac{1}{2}}
    \\
    &\leq C\left\|\mathbf{u}-\mathbf{\widetilde{u}}\right\|_{H^{1}(\Omega_{R})}
    \left\|\nabla\widetilde{\theta}\right\|_{L^2(\Omega_{R})}
    \left\|\theta-\widetilde{\theta}\right\|_{L^{\frac{3}{2}}(\Omega_{R})}^{\frac{1}{2}}.
    \end{aligned}
\end{align*}
Thus, one gets
\begin{align*}
    \begin{aligned}
    \frac{d}{dt}\left\|\theta-\widetilde{\theta}\right\|_{L^{\frac{3}{2}}(\Omega_{R})}^2
    =\frac{4}{3}\left\|\theta-\widetilde{\theta}\right\|_{L^{\frac{3}{2}}(\Omega_{R})}^{\frac{1}{2}}
    \frac{d}{dt}\left\|\theta-\widetilde{\theta}\right\|_{L^{\frac{3}{2}}(\Omega_{R})}^{\frac{3}{2}}
    &\leq C\left\|\mathbf{u}-\mathbf{\widetilde{u}}\right\|_{H^1(\Omega_{R})}
    \left\|\nabla\widetilde{\theta}\right\|_{L^2(\Omega_{R})} \left\|\theta-\widetilde{\theta}\right\|_{L^{\frac{3}{2}}(\Omega_{R})}
    \\
    &\leq \frac{\mu}{4}
    \left\|\mathbf{u}-\mathbf{\widetilde{u}}\right\|_{H^1(\Omega_{R})}^2
    +C(t)\left\|\theta-\widetilde{\theta}\right\|_{L^{\frac{3}{2}}(\Omega_{R})}^{2},
    \end{aligned}
\end{align*}
where $C(t)\in L^{1}(0,T)$ is non-negative integrable function. Furthermore, hoosing $\varepsilon>0$ small enough, from the above inequality and \eqref{weiyixing1}, one has
\begin{align*}
    \frac{d}{dt}\left[\left\|\mathbf{u}-\mathbf{\widetilde{u}}\right\|_{L^2(\Omega_{R})}^2+
    \left\|\theta-\widetilde{\theta}\right\|_{L^{\frac{3}{2}}(\Omega_{R})}^2\right]
    \leq \left(\frac{C\varepsilon+2\mu}{4}+\frac{1}{2\varepsilon}+C(t)\right)
   \left[\left\|\mathbf{u}-\mathbf{\widetilde{u}}\right\|_{L^2(\Omega_{R})}^2+
    \left\|\theta-\widetilde{\theta}\right\|_{L^{\frac{3}{2}}(\Omega_{R})}^2\right].
\end{align*}
From Gronwall's inequality, one gets
\begin{align*}
    \left\|\mathbf{u}-\mathbf{\widetilde{u}}\right\|_{L^2(\Omega_{R})}^2+
    \left\|\theta-\widetilde{\theta}\right\|_{L^{\frac{3}{2}}(\Omega_{R})}^2\leq
    \left[\left\|\mathbf{u}_{0}-\mathbf{\widetilde{u}}_{0}\right\|_{L^2(\Omega_{R})}^2+
    \left\|\theta-\widetilde{\theta}_{0}\right\|_{L^{\frac{3}{2}(\Omega_{R})}}^2\right]
    e^{\int_{0}^{T}\left(\frac{C\varepsilon+2\mu}{4}+\frac{1}{2\varepsilon}+C(t)\right)dt}.
\end{align*}
Therefore, the strong solution $(\mathbf{u},\theta)$ continuously depends  on the initial values when the strong solution is considered in the norm of $\left[L^2\right]^2\times L^{\frac{3}{2}}$.
In particular, the uniqueness of strong solution is also obtained.
\end{proof}
\subsection{Proofs of conclusions on unbounded domain}
Based on the conclusions of $\Omega_{R}$, we shall pay attention to the proof of  well-posedness on $\Omega_{\infty}$.
\begin{proof}\textbf{Proof of Theorem \ref{cunzairuojie1}.}

    For each $R>0$, we define
    \begin{align}\label{jieduan}
        \mathbf{u}_{0}^{R}=
        \begin{cases}
            \mathbf{u}_{0},~~x\in\Omega_{R},
            \\
            0,~~~~x\in\Omega_{\infty}/\Omega_{R},
        \end{cases}
        \theta_{0}^{R}=
        \begin{cases}
            \theta_{0},~~x\in\Omega_{R},
            \\
            0,~~~~x\in\Omega_{\infty}/\Omega_{R}.
        \end{cases}
    \end{align}
   Then, $\mathbf{u}_{0}^{R}$ $\rightarrow$ $\mathbf{u}_{0}$
   in $\left[H^{1}(\Omega_{\infty})\right]^2$
   and $\theta_{0}^{R}$ $\rightarrow$ $\theta_{0}$ in
   $L^{\infty}(\Omega_{\infty})\cap L^{1}(\Omega_{\infty})$ as $R\rightarrow+\infty$.

   From Theorem \ref{cunzairuojie}, we know that there exist $\mathbf{u}^{R}$ and $\theta^{R}$
   satisfying the problem \eqref{youjiequyu} under the initial values
   $\mathbf{u}_{0}^{R}$ and $\theta_{0}^{R}$ in the weak sense.
   Moreover, one has
   \begin{align}\label{jixianbaohaoxing}
        \begin{aligned}
        & \left\|\mathbf{u}^R\right\|_{L^2(\Omega_{R})}^2\leq \left(\left\|\theta_{0}^R\right\|_{L^2(\Omega_{R})}^2
        +\left\|\mathbf{u}_{0}^R\right\|_{L^2(\Omega_{R})}^2\right)e^{|\overrightarrow{
            \mathbf{g}
        }|T}
        \leq \left(\big{\|}\theta_{0}\big{\|}_{L^2(\Omega_{R})}^2
        +\big{\|}\mathbf{u}_{0}\big{\|}_{L^2(\Omega_{R})}^2\right)e^{|\overrightarrow{
            \mathbf{g}
        }|T},
        \\
        &
        \begin{aligned}
            2\mu\int_{0}^{T}\left\|\nabla\mathbf{u}^R\right\|^{2}_{L^2(\Omega_{R})}ds
        &\leq
        \left\|\mathbf{u}_{0}^R\right\|_{L^2(\Omega_{R})}^2+
       \left(\left\|\theta_{0}^R\right\|_{L^2(\Omega_{R})}^2+\left\|\mathbf{u}_{0}^R\right\|_{L^2(\Omega_{R})}^2\right)e^{|\mathbf{g}|T}
       +|\mathbf{g}|\left\|\theta_{0}^R\right\|_{L^2(\Omega_{R})}^2T
       \\
       &\leq
       \big{\|}\mathbf{u}_{0}\big{\|}_{L^2(\Omega_{R})}^2+
       \left(\big{\|}\theta_{0}\big{\|}_{L^2(\Omega_{R})}^2+\big{\|}\mathbf{u}_{0}\big{\|}_{L^2(\Omega_{R})}^2\right)e^{|\mathbf{g}|T}
       +|\mathbf{g}|\big{\|}\theta_{0}\big{\|}_{L^2(\Omega_{R})}^2T,
        \end{aligned}
        \\
        &\left\|\theta^R\right\|_{L^{k}(\Omega_{R})}=\left\|\theta_{0}^R\right\|_{L^{k}(\Omega_{R})}\leq \big{\|}\theta_{0}\big{\|}_{L^{k}(\Omega_{R})},
        (k\in (1,+\infty)).
        \end{aligned}
 \end{align}
Thus, we can extend $\left(\mathbf{u}^{R},\theta^{R}\right)$ to $\Omega_{\infty}$ by introducing $\left(\mathbf{u}^{R},\theta^{R}\right)
=(\mathbf{0},0)$ in
$\Omega_{\infty}/\Omega_{R}$. This means that there exists a sequence $\{R_{j}\}$ such that $\left(\mathbf{u}^{R_{j}},\theta^{R_{j}}\right)$
converging to $(\mathbf{u},\theta)$,
\begin{align*}
    \mathbf{u}^{R_{j}}\rightarrow \mathbf{u} ~\text{weakly~star~in~}&L^{\infty}\left(0,T;\left[L^{2}(\Omega_{\infty})\right]^2\right),
    ~\nabla\mathbf{u}^{R_{j}}\rightarrow
    \mathbf{\nabla u} ~\text{weakly~in~} L^{2}\left(0,T;\left[L^{2}(\Omega_{\infty})\right]^{2}\right),
    \\
    &\theta^{R_{j}}\rightarrow \theta~\text{weakly~star~in~} L^{\infty}\left(0,T;L^{1}(\Omega_{\infty})
    \cap L^{\infty}(\Omega_{\infty})\right).
\end{align*}
Moreover, we also have
\begin{align*}
\mathbf{u}^{R_{j}}\rightarrow \mathbf{u} ~\text{strongly~in~} &L^{2}\left(0,T;L^{2}_{loc}(\Omega_{\infty})\right).
\end{align*}
Thus, one gets that $(\mathbf{u},\theta)$ satisfies \eqref{ruojie1}.

From \eqref{jixianbaohaoxing}, it follows that
\begin{align*}
    \mathbf{u}\in L^{\infty}\left(0,T;\left[L^{2}(\Omega_{\infty})\right]^2\right),~
    \mathbf{u}\in L^{2}\left(0,T;\mathbf{G}^{1,2}(\Omega_{\infty})\right),
    ~\theta\in L^{\infty}\left(0,T;L^{1}(\Omega_{\infty})
    \cap L^{\infty}(\Omega_{\infty})\right).
\end{align*}
And, taking $R\rightarrow\infty$ in \eqref{jixianbaohaoxing}, we deduce that $(\mathbf{u},\theta)$ meets the estimate \eqref{threeev}.
\end{proof}
\begin{proof}\textbf{Proof of Theorem \ref{qiangjie1}.}

    First, we consider the problem \ref{model} on bounded $\Omega_{R}$ with initial value
    $\mathbf{u}_{0}^{R}$ and $\theta_{0}^{R}$ defined as \eqref{jieduan}.
   By Theorem \ref{qiangjie}, it is found that the strong solution $\left(\mathbf{u}^{R},p^{R},
    \theta^R\right)$ satisfies
    \begin{align}\label{387}
        \begin{aligned}
            \big{\|}\partial_{t}&\mathbf{u}^R\big{\|}_{L^2(\Omega_{R})}^{2}+
            \left\|\nabla p^R\right\|_{L^{2}(\Omega_{R})}^2+
            \left\|\nabla\mathbf{u}^R\right\|_{H^{1}(\Omega_{R})}^{2}
            +\left\|\theta_{t}^R\right\|_{L^2(\Omega_{R})}^2
            \\&+\int_{0}^{t}\left[\left\|\partial_{s}\mathbf{u}^R(s)\right\|_{L^{2}(\Omega_{R})}^{2}
+\left\|\nabla\partial_{s}\mathbf{u}^R\right\|_{L^2(\Omega_{R})}^{2}
            +\left\|\nabla^2\mathbf{u}^R\right\|_{L^{4}(\Omega_{R})}^{2}+\left\|\nabla p^R\right\|_{L^{4}(\Omega_{R})}^2\right]ds
            \leq C.
        \end{aligned}
    \end{align}
Similarly, we extend $\left(\mathbf{u}^R,p^R,\theta^R\right)$ to $\Omega_{\infty}$ by defining
$\left(\mathbf{u}^R,p^R,\theta^R\right)=(\mathbf{0},0,0)$ in $\Omega_{\infty}/\Omega_{R}$. Hence,
there exists a sequence $\{R_{j}\}$ such that  $R_{j}\rightarrow\infty$,
\begin{align*}
    \begin{aligned}
    \partial_{t}&\mathbf{u}^{R_{j}},~\nabla p^{R_{j}},~\nabla\mathbf{u}^{R_{j}}\text{~and~}\theta_{t}^{R_{j}}
    \rightarrow \partial_{t}\mathbf{u},~\nabla p,~\nabla\mathbf{u}\text{~and~}\theta_{t}\text{~weakly~star~in~}
   L^{\infty}\left(0,T;\left[L^{2}(\Omega_{\infty})\right]^2\right),
   \\
   &L^{\infty}\left(0,T;\left[L^{2}(\Omega_{\infty})\right]^2\right),~
    L^{\infty}\left(0,T;
    \left[H^{1}(\Omega_{\infty})\right]^4\right)\text{~and~}L^{\infty}\left(0,T;L^{2}(\Omega_{\infty})\right),
    \end{aligned}
\end{align*}
and
\begin{align*}
    \begin{aligned}
        \partial_{t}&\mathbf{u}^{R_{j}},~\nabla\partial_{t}\mathbf{u}^{R_{j}},~
        \nabla^2\mathbf{u}^{R_{j}}\text{~and~}\nabla p^{R_{j}}~\rightarrow\partial_{t}\mathbf{u},~
        \nabla\partial_{t}\mathbf{u},~\nabla^2\mathbf{u}\text{~and~}
        \nabla p~
        \text{~weakly~in~}L^{2}\left(0,T;\left[L^2(\Omega_{\infty})\right]^2\right),
        \\
        &L^{2}\left(0,T;\left[L^2(\Omega_{\infty})\right]^2\right),~
        L^{2}\left(0,T;\left[L^4(\Omega_{\infty})\right]^6\right)\text{~and~}L^{2}\left(0,T;\left[L^4(\Omega_{\infty})\right]^2\right).
    \end{aligned}
\end{align*}
From \eqref{387},  one can conclude that \eqref{312} and \eqref{kanqilaimaodun1} are all valid.
The uniqueness of $\nabla p$ arises from the uniqueness of $\mathbf{u}$ and $\theta$. We end the proof of this theorem.
\end{proof}
\section{ Linear instability}\label{RTbuwendingxing0104}
The purpose of this section is to study the linear  instability of a steady state
\eqref{cond1} and \eqref{cond20105}. To investigate this linear instability, we need to seek a growing
solution to equation \eqref{naiverboundarycondition123}-\eqref{xianxingpart}, i.e.,
\begin{align}\label{xuexuezenmeshuo0129}
    \mathbf{v}=\widetilde{\mathbf{v}}(x_{1},x_{2})e^{\lambda t},~~\pi=\widetilde{q}(x_{1},x_{2})e^{\lambda t},
    ~~\Theta=\widetilde{\Theta}(x_{1},x_{2})e^{\lambda t}.
\end{align}
Putting \eqref{xuexuezenmeshuo0129} into \eqref{xianxingpart}, one has
\begin{align}\label{bianhua1}
    \begin{cases}
        \lambda \mathbf{\widetilde{v}}+\nabla \widetilde{q}=\mu\Delta\mathbf{\widetilde{v}}+\widetilde{\Theta}
        \mathbf{g},
        \\
        \lambda\widetilde{\Theta}+\widetilde{v}_{2}D\overline{\theta}=0,
        \\
        \nabla\cdot\widetilde{\mathbf{v}}=0,
    \end{cases}
\end{align}
and the boundary condition is as follows,
\begin{align}\label{naiverboundarycondition1234}
    \begin{aligned}
    &\widetilde{v}_{2}(x_{1},0)=\widetilde{v}_{2}(x_{1},1)=0,
    ~\partial_{x_{2}}\widetilde{v}_{1}(x_{1},1)=\frac{k_{1}}{\mu}\widetilde{v}_{1}(x_{1},1),
    \\
    &\partial_{x_{2}}\widetilde{v}_{1}(x_{1},0)=-\frac{k_{0}}{\mu}
    \widetilde{v}_{1}(x_{1},0),~x_{1}\in\mathbf{R}.
    \end{aligned}
\end{align}
For each spatial frequency $\xi\neq 0$, we can define the following new  functions,
\begin{align}\label{shaodeng0306}
    \begin{aligned}
        \widetilde{\mathbf{v}}=\left(-iU_{1}\left(x_{2}\right),U_{2}(x_{2})\right)e^{i \xi x_{1}},~~
        \widetilde{q}=\varpi (x_{2})e^{i\xi x_{1}},
        ~~\widetilde{\Theta}=h(x_{2})e^{i\xi x_{1}},
    \end{aligned}
\end{align}
where $U_{1},U_{2},\varpi,h:(0,1)\rightarrow\mathbf{R}$ are unknowns . Then putting \eqref{shaodeng0306} into \eqref{bianhua1} gives immediately the following system of ODEs,
\begin{align}\label{tezhengfangcheng2}
    \begin{cases}
        \lambda U_{1}-\xi \varpi =\mu\left(D^2-\xi ^2\right)U_{1},
        \\
        \lambda U_{2}+D\varpi =\mu\left(D^2-\xi ^2\right)U_{2}+g h,
        \\
        \lambda h+U_{2}D\overline{\theta}=0,
        \\
        \xi U_{1}+DU_{2}=0,
    \end{cases}
\end{align}
coupled with
\begin{align}\label{haimeichulai1215}
    U_{2}(0)=U_{2}(1)=0,
    ~DU_{1}(1)=\frac{k_{1}}{\mu}U_{1}(1),
    ~DU_{1}(0)=-\frac{k_{0}}{\mu}
    U_{1}(0).
\end{align}
Eliminating $\varpi$, $U_{1}$ and $h$ in \eqref{tezhengfangcheng2}, one gets a four-order ODE of $U_{2}$,
\begin{align}\label{equ1}
    \begin{aligned}
    -\lambda^2\left( D^2 U_{2}-\xi ^2 U_{2}\right)=-\lambda\mu\left[D^4 U_{2}-2 \xi ^2 D^2 U_{2}+\xi ^4 U_{2}\right]
    -\left(g\xi ^2D\overline{\theta}\right)U_{2},
    \end{aligned}
\end{align}
subjected to boundary conditions:
\begin{align}\label{mofangxietisheng}
    U_{2}(0)=U_{2}(1)=0,~~D^2U_{2}(1)=\frac{k_{1}}{\mu}DU_{2}(1),~~
    D^2U_{2}(0)=-\frac{k_{0}}{\mu}DU_{2}(0).
\end{align}
If \eqref{equ1}-\eqref{mofangxietisheng} has a solution with $\lambda>0$ for some frequency $\xi$,
then the steady state $\left(\mathbf{0},\overline{p},\overline{\theta}\right)$ is regarded as linear instability.
Here, the order of $\lambda$ is 2 which causes a challenge to employ the normal variational method directly. Hence,
we adopt a modified variational method to deal with problem \eqref{equ1}-\eqref{mofangxietisheng}, that is, we consider the following family of ODEs,
\begin{align}\label{equ2}
    \begin{aligned}
    -\lambda^2\left( D^2 U_{2}-\xi ^2 U_{2}\right)=-s\mu\left[D^4 U_{2}-2 \xi ^2 D^2 U_{2}+\xi ^4 U_{2}\right]
    -\left(g\xi ^2D\overline{\theta}\right)U_{2},
    \end{aligned}
\end{align}
where $s>0$ and $\xi\neq 0$.
For convenience, let $W=U_{2}$. Multiplying \eqref{equ2} by $W$ and integrating over $[0,1]$ gives the following
expression,
\begin{align*}
-\lambda^2=\frac{E(W)}{J(W)},
\end{align*}
where
\begin{align}\label{xiangguan0310}
    \begin{aligned}
        &E(W)=E(W,s,\xi):=s E_{0}(W)+s E_{1}(W)+E_{2}(W),
        \\
        &E_{0}(W,\xi ):=E_{0}(W):=\mu\int_{0}^{1}\left[|D^2 W|^2+2 \xi ^2 |DW|^2+\xi ^4 W^2\right]  dx_{2},
        \\
        &E_{1}(W):=-k_{1}\left(DW(1)\right)^2-k_{0}\left(DW(0)\right)^2,
        \\
        &E_{2}(W,\xi ):=E_{2}(W):=g \xi ^2\int_{0}^{1}D\overline{\theta}W^2 dx_{2},
        \\
        &J(W,\xi ):=J(W):=\int_{0}^{1} |DW|^2+\xi ^2 W^2 dx_{2}.
    \end{aligned}
\end{align}
In order to study the existence of $W$
for the above equality, we need to consider the following extreme value problem,
\begin{align}\label{xianxingzuixiaozhi}
\Phi(s):=\Phi(s,\xi):=\inf\limits_{W\in H^2 \cap H_{0}^{1}}\frac{E(W)}{J(W)}.
\end{align}
 We will see that the \eqref{xianxingzuixiaozhi} will be solved in the Proposition \ref{pro1}. And note that
 the above functions in \eqref{xiangguan0310} and \eqref{xianxingzuixiaozhi} are the function of $s$ or $\xi$. Sometimes we will omit the two variables for
convenience if there is no ambiguity.
By the homogeneity, \eqref{xianxingzuixiaozhi} is equivalent to
\begin{align}\label{zhegezhege1}
    \Phi(s):=\Phi(s,\xi):=\inf\limits_{W\in \mathcal{A}_{\xi }}E(W),
\end{align}
where the admissible set is $\mathcal{A}_{\xi }=\left\{W\in H^2\cap H_{0}^{1}\big{|}J(W)=1\right\}$.

Now, we shall state the following conclusion about linear instability.
\begin{theorem}[Linear instability]\label{xianxingbuwending}
    The steady state $\left(\mathbf{0},\overline{p},\overline{\theta}\right)$ is linearly unstable in $H^k-$norm, for
    any $k\in\mathbf{N}$, in the sense that there are exponentially growing mode solutions to the
    linearized perturbed problem \eqref{naiverboundarycondition123}-\eqref{xianxingpart} in
    $H^{k}$.
\end{theorem}
\begin{remark}The existence of a growing normal mode is evident when the equation $\Phi(s)=-s^2$ possesses
    a solution in $(0,+\infty)$ which indicates the linear instability of the steady state
$\left(\mathbf{0},\overline{p},\overline{\theta}\right)$. Through a rigorous analysis for the
    properties of $\Phi(s)$ outlined in Propositions \ref{pro1}-\ref{jiecunzaijiuhao},
    we establish the existence of a solution to $\Phi(s)=-s^2$ in $(0,+\infty)$, confirming the
    existence of a solution to problem \eqref{tezhengfangcheng2}
    for each $\xi\neq 0$ $($see Theorem \ref{jiecunzai1215}$)$. Subsequently, by employing a linear combination of such solutions
    $($see Theorem \ref{gouzaojie}$)$, we can validate Theorem \ref{xianxingbuwending}.
\end{remark}

\subsection{Analysis for the variational problem (4.12)}
In order to prove the instability of Theorem \ref{xianxingbuwending}, we first discuss the properties of $\Phi(s)$ appearing in
variational problem \eqref{zhegezhege1} with a fixed $\xi$ by variational method.
Here, the main task is to verify that the minimum of $E(W)$ can be
achieved over $\mathcal{A}_{\xi }$, see Proposition \ref{pro1}, which implies that $\Phi(s)$ is well-defined. And we also show that
this minimizer solves the Euler-Lagrange equation equivalent to problem \eqref{mofangxietisheng}-\eqref{equ2}
when
$-\lambda^2=\inf\limits_{W\in \mathcal{A}_{\xi }}E(W)$, see Proposition \ref{jiecunzaizhengming}.
\begin{proposition}\label{pro1}
    For any fixed $\xi\neq 0$, $E(W)$ achieves its minimum on $\mathcal{A}_{\xi }$.
\end{proposition}

\begin{proof}
   By a direct calculations, one has
    \begin{align}\label{fangshou}
        \begin{aligned}
            k_{1}(DW(1))^2+k_{0}(DW(0))^2&=\int_{0}^{1}
            \frac{d}{dx_{2}}\left\{\left[(k_{1}+k_{0})x_{2}-k_{0}\right](DW(x_{2}))^2\right\}
             dx_{2}
             \\
             &\leq C(k_{0},k_{1},\varepsilon)\big{\|}DW\big{\|}_{L^2}^2+\varepsilon\big{\|}D^2 W\big{\|}_{L^2}^2.
        \end{aligned}
    \end{align}
     Let $\varepsilon=\mu$ in \eqref{fangshou}. There appears  the relation
    \begin{align}\label{jiushini1216}
        E(W)\geq -C(k_{0},k_{1},\mu)\big{\|}DW\big{\|}_{L^2}^2+g\xi^2 \int_{0}^{1}D\overline{\theta}W^2dx_{2}
        \geq -C(k_{0},k_{1},\mu)-g\big{\|}D\overline{\theta}\big{\|}_{L^\infty},
    \end{align}
    which indicates that there exists the lower bound for $E(W)$.

    Therefore,  a minimizing sequence  $\{W_{n}\}\in \mathcal{A}_{\xi }$ can be found and then satisfies
    \begin{align*}
        \lim\limits_{n\rightarrow \infty}E(W_{n})=\inf\limits_{W\in\mathcal{A}_{\xi }}E(W).
    \end{align*}
Without loss of generality, we assume $E(W_{n})\leq \inf\limits_{W\in\mathcal{A}_{\xi }}E(W)+1$, that is,
    \begin{align}\label{dongde}
        s E_{0}(W_{n})+s E_{1}(W_{n})+E_{2}(W_{n})\leq \inf\limits_{W\in\mathcal{A}_{\xi }}E(W)+1.
    \end{align}

    On the other hand, let $\varepsilon=\frac{\mu}{2}$ in \eqref{fangshou}. By \eqref{dongde}, one has
    that $\{W_{n}\}\in H^2$ is bounded for any fixed $s>0$. Thus, there exists
    a $W_{0}\in H^2$ such that $W_{n}\rightarrow W_{0}$ weakly in $H^2$ and strongly in $H^1_{0}.$
    Since $s E_{0}(W)$ is convex and according to strong convergence, one can conclude that the $E(W)$ is weakly lower
    semi-continuous. Hence,
    \begin{align*}
        E(W_{0})\leq \lim\limits_{\overline{n\rightarrow\infty}}E(W_{n})=
        \inf\limits_{W\in\mathcal{A}_{\xi }}E(W).
    \end{align*}
    And because of $W_{n}\rightarrow W_{0}$ strongly in $H^{1}_{0}$, we have
    \begin{align*}
        W_{0}\in\mathcal{A}_{\xi }.
    \end{align*}
Therefore, the proof is completed.
    \end{proof}
    Next, we will show that the minimizer constructed above satisfies an Euler-Lagrangian equation
equivalent to \eqref{equ2} when $-\lambda^2=\Phi(s)$.

\begin{proposition}\label{jiecunzaizhengming}
    Assume $-\lambda_{0}^2=\Phi(s)=E(W_{0})$, i.e, $W_{0}$ is the minimizer,
    then $W_{0}$ is smooth and $(W_{0},\lambda_{0})$ solves the equation \eqref{mofangxietisheng}-\eqref{equ2}.
\end{proposition}
\begin{proof}
For any $ W\in C_{0}^{\infty}((0,1)),$$~t,r\in\mathbb{R}$, we define the functional
\begin{align*}
    \begin{aligned}
        I(t,r)=J(W_{0}+tW+rW_{0}).
    \end{aligned}
\end{align*}
Then, a series of calculation gives
\begin{align*}
    \begin{aligned}
        &I(0,0)=1,
        \\
        &\partial_{t}I(0,0)=2\int_{0}^{1} \left[DW_{0}DW+\xi ^2 W_{0} W \right]dx_{2},
        \\
        &\partial_{r}I(0,0)=2\int_{0}^{1} \left[|DW_{0}|^2+\xi ^2 W_{0}^2\right]dx_{2}=2\neq 0.
    \end{aligned}
\end{align*}
By the implicit function theorem, there exists a function $r=r(t)$ defined on near $0$, such that
$r(0)=0$, $I(t,r(t))=1$ and $\dot{r}(0)=-\int_{0}^{1} \left[DW_{0}DW+\xi ^2 W_{0}W\right] dx_{2}.$

Furthermore, we consider the functional
\begin{align*}
    \widetilde{I}(t)=E(W_{0}+tW+r(t)W_{0}).
\end{align*}
Since $W_{0}$ is the extreme point, one has
\begin{align*}
    \begin{aligned}
    0=\dot{\widetilde{I}}(0)&=2\bigg\{s\mu \int_{0}^{1}\left[D^2 W_{0}D^2 W
    +2\xi ^2 DW_{0}DW+\xi ^4 W_{0}W\right]
     dx_{2}
     \\
     &
     +g \xi ^2\int_{0}^{1}D\overline{\theta} W_{0}W dx_{2}
     +\bigg{\{}s\mu\int_{0}^{1}\left[|D^2W_{0}|^2+2\xi ^2|DW_{0}|^2+\xi ^4 W_{0}^2\right] dx_{2}
     \\
     &-s k_{1}\left(DW_{0}(1)\right)^2
     -s k_{0}\left(DW_{0}(0)\right)^2+g\xi ^2
     \int_{0}^{1}D\overline{\theta}W_{0}^2dx_{2}\bigg{\}}\dot{r}(0)\bigg\}.
    \end{aligned}
\end{align*}
Moreover, when $-\lambda_{0}^2=E(W_{0})$, it follows that
\begin{align*}
    \begin{aligned}
    &s\mu\int_{0}^{1}\left[D^2 W_{0}D^2 W+2 \xi ^2 DW_{0}DW+\xi ^4 W_{0}W\right]dx_{2}
    \\
    &~~+g \xi ^2\int_{0}^{1}D\overline{\theta}W_{0} W dx_{2}
    +\lambda_{0}^2\int_{0}^{1}\left[ DW_{0}DW+\xi ^2 W_{0}W\right] dx_{2}=0,
    \end{aligned}
\end{align*}
which means that $W_{0}$ solves the equation \eqref{equ2} in the weak sense. Standard bootstrap arguments
(one can refer to the procedure in Appendix \ref{appendix11226})
 show that the solution is smooth.

Next, let $W\in H^2\cap H_{0}^1$. Then according to
$0=\dot{\widetilde{I}}(0)=\frac{d}{dt}\big|_{t=0}E(W_{0}+tW+r(t)W_{0})$, we have
\begin{align}\label{jingjing}
    \begin{aligned}
        s&\mu \int_{0}^{1}\left[
            D^2 W D^2 W_{0}+2\xi^2 DW DW_{0}+\xi^4W W_{0}
        \right]dx_{2}
        +\int_{0}^{1}g\xi^2D\overline{\theta}WW_{0}dx_{2}
        \\
        &+\lambda_{0}^2\int_{0}^1\left[ DW_{0}DW+\xi^2W_{0}W\right]dx_{2}=
        s\left[k_{1}DW(1)DW_{0}(1)+k_{0}DW(0)DW_{0}(0)\right].
    \end{aligned}
\end{align}
Multiplying \eqref{equ2} with $W\in H^2\cap H_{0}^1$ and integrating by part lead to
\begin{align}\label{jingjingjingjing}
    \begin{aligned}
        -s\mu \int_{0}^{1}D^3W_{0}DW dx_{2}&+s\mu
        \int_{0}^{1}\left[
        2\xi^2 DW DW_{0}+\xi^4W W_{0}
        \right]dx_{2}+\int_{0}^{1}g\xi^2D\overline{\theta}WW_{0}dx_{2}
        \\
        &+\lambda_{0}^2\int_{0}^1\left[DW_{0}DW+\xi^2W_{0}W\right]dx_{2}=0.
    \end{aligned}
\end{align}
Then, comparing \eqref{jingjing} with \eqref{jingjingjingjing} gives to
\begin{align*}
    \mu\int_{0}^{1}D\left(D^2W_{0}DW\right) dx_{2}=k_{1}DW(1)DW_{0}(1)+k_{0}DW(0)DW_{0}(0),
\end{align*}
which is equivalent to
\begin{align*}
    \left(\mu D^2W_{0}(1)-k_{1}DW_{0}(1)\right)DW(1)=\left(\mu D^2 W_{0}(0)+k_{0}DW_{0}(0)\right)DW(0).
\end{align*}
Since $W$ is arbitrarily chosen, then
\begin{align*}
    D^2W_{0}(1)=\frac{k_{1}DW_{0}(1)}{\mu},~~D^2W_{0}(0)=-\frac{k_{0}DW_{0}(0)}{\mu}.
\end{align*}
Thus, we end the proof of  Proposition \ref{jiecunzaizhengming}.
\end{proof}
Next, we want to verify that $-s^2=\Phi(s)$ has a solution in $(0,+\infty)$. To this end, we first give some
properties of $\Phi(s)$ as a function defined on $(0,+\infty)$: locally bounded (Proposition
\ref{youjie1223}); locally Lipschitz continuous (Proposition \ref{lianxuxing1215}); monotonous
(Proposition \ref{daodiaoxing1215}) and sign changeable (Propositions \ref{cunzai1215} and
 \ref{zugouxiao1215}).

\begin{proposition}\label{youjie1223} $\Phi(s,\xi)$ is bounded on $I_{1}\times I_{2}$, that is
    there exists a positive constant $M(I_{1},I_{2})$, dependent of $I_{1}$ and $I_{2}$, such that
    \begin{align*}
        |\Phi(s,\xi)|\leq M(I_{1},I_{2}),
    \end{align*}
     where
    $I_{1}=[a,b]\subset(0,+\infty)$ and $I_{2}=[c,d]\subset(0,+\infty)$ are two bounded sets.
\end{proposition}
\begin{proof}
    We first know from \eqref{jiushini1216} that
    \begin{align*}
        \Phi(s,\xi)\geq -C(k_{0},k_{1},\mu)-g\big{\|}D\overline{\theta}\big{\|}_{L^\infty}>-\infty.
    \end{align*}
    Thus, if $\Phi(s,\xi)$ is unbounded on $I_{1}\times I_{2}$, then
    there exist $\widetilde{s}\in I_{1}$ and $\widetilde{\xi}^2\in I_{2}$ such that
    $$\Phi\left(\widetilde{s},\widetilde{\xi}\right)=+\infty.$$
    Then according to the definition of $\Phi(s,\xi)$, one has
    \begin{align*}
        \frac{E\left(W,\widetilde{s},\widetilde{\xi}\right)}{J\left(W,\widetilde{\xi}\right)}=+\infty,~~
        \forall W\in H^2 \cap H_{0}^1.
    \end{align*}
This is impossible case. Hence, we get the conclusion
 of this proposition.
\end{proof}

\begin{proposition}\label{lianxuxing1215}
    $\Phi(s,\xi)\in C_{loc}^{0,1}\left((0,+\infty)\times (0,+\infty)\right)$.
    \end{proposition}
    \begin{proof}
        Let $I_{1}=[a,b]\subset (0,+\infty)$ and $I_{2}=[c,d]\subset(0,+\infty)$ be two bounded
        intervals. We shall divide its proof into two steps.

        (1) ~\textbf{$\Phi(s,\xi)$ is Lipschitz continuous with respect to $s$ on $I_{1}$.}

        $\forall $ $s_{1}\leq s_{2}$ and $s_{1},$ $s_{2}\in I_{1},$ then
        there exist $W_{s_{1}}$ and $W_{s_{2}}\in H^2\cap H_{0}^{1}$ such that
        \begin{align*}
            \begin{aligned}
                \Phi(s_{1})=\frac{E\left(W_{s_{1}},s_{1}\right)}{J\left(W_{s_{1}}\right)},
                ~~\Phi(s_{2})=\frac{E(W_{s_{2}},s_{2})}{J(W_{s_{2}})}.
            \end{aligned}
        \end{align*}
        According to the definition of minimum, we have
        \begin{align}\label{yinci}
            \begin{aligned}
                &\Phi(s_{1})-\Phi(s_{2})\leq
                \frac{E(W_{s_{2}},s_{1})}{J(W_{s_{2}})}-
                \frac{E(W_{s_{2}},s_{2})}{J(W_{s_{2}})}
                =\frac{\left(s_{1}-s_{2}\right)\left(E_{0}(W_{s_{2}})
                +E_{1}(W_{s_{2}})\right)}{J(W_{s_{2}})},
                \\
                &\Phi(s_{2})-\Phi(s_{1})\leq
                \frac{(s_{2}-s_{1})(E_{0}(W_{s_{1}})
                +E_{1}(W_{s_{1}}))}{J(W_{s_{1}})}.
            \end{aligned}
        \end{align}
    Due to the Proposition \ref{youjie1223}, there exists a positive constant
    $M(I_{1},I_{2})$, such that
    \begin{align*}
       |\Phi(s,\xi)|\leq M(I_{1},I_{2}).
    \end{align*}
     It follows that
     \begin{align}\label{20240312}
        \bigg{|}\frac{s_{i}E_{0}(W_{s_{i}})+s_{i}E_{1}(W_{s_{i}})
        +E_{2}(W_{s_{i}})}{J(W_{s_{i}})}\bigg{|}\leq M(I_{1},I_{2}).
     \end{align}
     In addition, it is easy to see that
     \begin{align}\label{202403121}
        \bigg{|}\frac{E_{2}(W_{s_{i}})}{J(W_{s_{i}})}\bigg{|}\leq g\big{\|}D\overline{\theta}\big{\|}_{L^\infty}.
    \end{align}
    In view of \eqref{20240312} and
    \eqref{202403121}, one has
    \begin{align*}
        \bigg{|}\frac{s_{i}E_{0}(W_{s_{i}})+s_{i}E_{1}(W_{s_{i}})
        }{J(W_{s_{i}})}\bigg{|}\leq M(I_{1},I_{2})+g\big{\|}D\overline{\theta}\big{\|}_{L^\infty}.
    \end{align*}
Thus, from \eqref{yinci}, we have
    \begin{align*}
        |\Phi(s_{1})-\Phi(s_{2})|\leq
        \frac{|s_{2}-s_{1}|}{a}\left[M(I_{1},I_{2})+g\big{\|}D\overline{\theta}\big{\|}_{L^\infty}\right].
    \end{align*}

    (2)~\textbf{$\Phi(s,\xi)$ is Lipschitz continuous about $\xi$ on $I_{2}$.}

    $\forall $ $\xi_{1}\leq \xi_{2}$ and $\xi_{1},$ $\xi_{2}\in I_{2},$ then
    there exist $W_{\xi_{1}}$ and $W_{\xi_{2}}\in H^2\cap H_{0}^{1}$ such that
    \begin{align*}
        \begin{aligned}
            \Phi(s,\xi_{1})=\frac{E(W_{\xi_{1}},s,\xi_{1})}{J(W_{\xi_{1}},\xi_{1})},
            ~~\Phi(s,\xi_{2})=\frac{E(W_{\xi_{2}},s,\xi_{2})}{J(W_{\xi_{2}},\xi_{2})}.
        \end{aligned}
    \end{align*}
Similarly, by virtue of the bound of $\Phi(s,\xi)$ on $I_{1}\times I_{2}$ and the definition of
minimum, there is also a positive constant $\widetilde{M}(I_{1},I_{2})$ such that
\begin{align*}
    |\Phi(s,\xi_{1})-\Phi(s,\xi_{2})|\leq \widetilde{M}(I_{1},I_{2})|\xi_{2}-\xi_{1}|.
\end{align*}
Let $M_{I_{1}I_{2}}=\max\left\{\widetilde{M}(I_{1},I_{2}),\frac{M(I_{1},I_{2})+g\left\|D\overline{\theta}\right\|_{L^\infty}}{a}\right\}$. Then,
\begin{align*}
    \begin{aligned}
        &|\Phi(s_{1},\xi)-\Phi(s_{2},\xi)|\leq M_{I_{1}I_{2}}|s_{2}-s_{1}|,
        \\
        &|\Phi(s,\xi_{1})-\Phi(s,\xi_{2})|\leq M_{I_{1}I_{2}}|\xi_{2}-\xi_{1}|.
    \end{aligned}
\end{align*}
Thus, one can easily obtain that
\begin{align*}
|\Phi(s_{1},\xi_{1})-\Phi(s_{2},\xi_{2})|\leq M_{I_{1}I_{2}}(
|s_{1}-s_{2}|+|\xi_{1}-\xi_{2}|
),
\end{align*}
which completes the proof.
    \end{proof}
    \begin{proposition}\label{daodiaoxing1215}

        Under the condition $\max\{k_{0},k_{1}\}\leq 0$,
        $\Phi(s)$ defined on $(0,+\infty)$ is increasing.
    \end{proposition}
    \begin{proof}
    Let $ 0<s_{1}<s_{2}$. According to Proposition \ref{pro1}, there exist
    $W_{s_{1}},W_{s_{2}}\in H^2\cap H_{0}^{1}$ such that
    \begin{align*}
        \Phi(s_{1})=\frac{E(W_{s_{1}},s_{1})}{J(W_{s_{1}})},
        ~~\Phi(s_{2})=\frac{E(W_{s_{2}},s_{2})}{J(W_{s_{2}})}.
    \end{align*}
   From the definition of minimum, one has
    \begin{align}\label{dandiao1}
        \begin{aligned}
        \frac{E(W_{s_{1}},s_{1})}{J(W_{s_{1}})}
        \leq \frac{E(W_{s_{2}},s_{1})}{J(W_{s_{2}})}.
        \end{aligned}
    \end{align}
Utilizing the condition $\max\{k_{0},k_{1}\}\leq 0$, one can deduce that $ E_{0}(W)+ E_{1}(W)>0$, which indicates that
    $\frac{E(W,s)}{J(W)}$ is increasing with respect to $s$.
    Then,
    \begin{align}\label{dandiao2}
        \frac{E(W_{s_{2}},s_{1})}{J(W_{s_{2}})}
        < \frac{E(W_{s_{2}},s_{2})}{J(W_{s_{2}})}=\Phi(s_{2}).
    \end{align}
    Combing \eqref{dandiao1} with \eqref{dandiao2}, we can obtain
    \begin{align*}
        \frac{E(W_{s_{1}},s_{1})}{J(W_{s_{1}})}
        < \frac{E(W_{s_{2}},s_{2})}{J(W_{s_{2}})},
    \end{align*}
    that is, $\Phi(s_{1})< \Phi(s_{2})$ meaning that
    $\Phi(s)$ is increasing on $(0,+\infty)$.
    \end{proof}

     To insure $-s^2=\Phi(s)$ for some  $s>0$,
    it requires that there exists a $W\in H^2\cap H_{0}^{1}$ such that
    \begin{align*}
        s E_{0}(W)+s E_{1}(W)+E_{2}(W)<0,
    \end{align*}
    which guarantees that $\Phi(s)<0.$

    Thus, we will consider the following extreme value problem:
    \begin{align}\label{variational problem1}
        s< \sup\limits_{W\in H^2\cap H_{0}^1}\frac{-E_{2}(W)}{E_{0}(W)+E_{1}(W)}=\sup\limits_{W\in\mathcal{B}}
        -E_{2}(W),
    \end{align}
    where $\mathcal{B}=\left\{W\in H^2\cap H_{0}^{1}\big{|}E_{0}(W)+E_{1}(W)=1\right\}$.
    \begin{proposition}\label{cunzai1215}
        Under the condition \eqref{cond1} and $\forall \xi \neq 0$, then $-E_{2}(W)$ achieves its maximum
        (denoted as $\lambda_{c}(\xi )$)  on $\mathcal{B}$ and $\lambda_{c}(\xi )\in \left(0,
        \frac{g \left\|D\overline{\theta}\right\|_{L^\infty}}{2\mu} \right)$.
    \end{proposition}
    \begin{proof}
        Owe to the condition \eqref{cond1}, there exists a $r>0$ such that $D\overline{\theta}(x_{2})<0$ for
        $x_{2}\in B\left(x_{2}^{0},r\right)$.
Thus, choose a $\overline{W}\in C_{0}^{\infty}\left(B\left(x_{2}^{0},r\right)\right)$, it is easy to get
         $-E_{2}\left(\overline{W}\right)>0,$ which indicates that
        \begin{align}\label{haishiyongdedao}
            \lambda_{c}(\xi )=\sup\limits_{W\in\mathcal{B}}-E_{2}(W)>0.
        \end{align}
        On the other hand, since $W\in\mathcal{B},$ we have
        \begin{align*}
            \xi ^2\int_{0}^{1}W^2 dx_{2}\leq \xi^2 \int_{0}^{1}|DW|^2 dx_{2}
            \leq \frac{1}{2\mu}.
        \end{align*}
        Thus, we conclude that
        \begin{align*}
            |E_{2}(W)|=\bigg{|}\int_{0}^{1}g \xi ^2 D\overline{\theta}W^2 dx_{2}\bigg{|}\leq
            \frac{g \left\|D\overline{\theta}\right\|_{L^\infty}}{2\mu}.
        \end{align*}

 Subsequently, we outline the key points to prove that $-E_{2}(W)$ can achieve its maximum.
        First, $-E_{2}(W)$ is bounded on $\mathcal{B}$; Second, $-E_{2}(W)$ is weakly upper semi-continuous;
        Finally, by the method of contradiction, we can prove that  $-E_{2}(W)$ can achieve its maximum on
        $\mathcal{B}$.
\end{proof}
\begin{remark}\label{zongshiruci} With respect to $\lambda_{c}(\xi)$, we have the estimate:
            when $W\in\mathcal{B}$,
            \begin{align}
                 |E_{2}(W)|=\bigg{|}\int_{0}^{1}g \xi ^2 D\overline{\theta}W^2 dx_{2}\bigg{|}\leq
                 \frac{\left\|gD\overline{\theta}\right\|_{L^\infty}}{\mu \xi^2},
            \end{align}
        which implies that
        \begin{align}\label{lingling}
            \lim\limits_{\xi\rightarrow \infty}\lambda_{c}(\xi)=0.
        \end{align}
    \end{remark}
    \begin{proposition}\label{zugouxiao1215}
        Under the condition \eqref{cond1} and $\forall \xi \neq 0$,
        then when $0<s<\lambda_{c}(\xi )$, $\Phi(s)<0$;
        when  $s \geq \lambda_{c}(\xi )$, $\Phi(s)\geq 0$.
    \end{proposition}
    \begin{proof}
    According to Proposition \ref{cunzai1215}, the existence of $\lambda_{c}(\xi )$ is verified.

    When $0<s<\lambda_{c}(\xi )$, by the definition of $\lambda_{c}(\xi )$, there exists a
    $\widehat{W}$ such that
    $s<-\frac{E_{2}\left(\widehat{W}\right)}{E_{1}\left(\widehat{W}\right)+E_{0}\left(\widehat{W}\right)}$, which shows that
    \begin{align*}
        s E_{0}\left(\widehat{W}\right)+s E_{1}\left(\widehat{W}\right)+E_{2}\left(\widehat{W}\right)<0.
    \end{align*}
    Then, we deduce that $\Phi(s)<0$.

    When $s>\lambda_{c}(\xi )$ and $\forall~W\in\mathcal{B}$, we have
    \begin{align*}
        s E_{0}(W)+s E_{1}(W)+E_{2}(W)>0,
    \end{align*}
    which proves that $\Phi(s)>0$.

    Thus, according to the continuity of $\Phi(s)$ with respect to $s$ (see  Proposition
    \ref{lianxuxing1215}), we can find that $\Phi(\lambda_{c}(\xi ),\xi )=0.$
    \end{proof}

    Next, based on the above propositions, we can show that the equation
    $-s^2=\Phi(s)$ has a solution on $(0,+\infty)$.

    \begin{proposition}\label{jiecunzaijiuhao}
        Under the condition \eqref{cond1} and $\forall \xi \neq 0$, there exists a
        $\lambda_{0}(\xi )\in (0, \lambda_{c}(\xi ))$ satisfying
        $-\lambda_{0}^2(\xi )=\Phi(\lambda_{0}(\xi ),\xi ).$
    \end{proposition}
    \begin{proof}

    Let $f(s)=\Phi(s)+s^2$. By the virtue of \eqref{haishiyongdedao}, we have
    \begin{align*}
        \begin{aligned}
            \Phi(0)=\inf\limits_{W\in H^2\cap H_{0}^1}\frac{E_{2}(W)}{J(W)}<0,
        \end{aligned}
    \end{align*}
    which implies that $f(0)<0$. Besides, when $s=\lambda_{c}(\xi )$, it follows that $f(s)=
    \lambda_{c}^2(\xi )>0$. Then, by the intermediate value theorem of continuous function,
    there exists a $\lambda_{0}\in (0,\lambda_{c}(\xi ))$ satisfying $f(\lambda_{0})=0$, that is,
    $-\lambda_{0}^2=\Phi(\lambda_{0}).$
\end{proof}

    Consequently, in view of above Propositions \ref{pro1}-\ref{jiecunzaijiuhao}, we conclude  the
existence of solutions for problem \eqref{equ1}-\eqref{mofangxietisheng}.
\begin{theorem}\label{jiecunzai1215} For each $|\xi|\neq 0$, there exist
    $U_{2}=U_{2}(|\xi|,x_{2})$ and $\lambda(|\xi|)>0$ satisfying \eqref{equ1}-
    \eqref{mofangxietisheng}. Moreover, $U_{2}\in H^{k}([0,1])$ for any positive integer $k$.
\end{theorem}

Since the solution $\lambda_{0}(\xi)$ to equality $-s^2=\Phi(s)$ depends on $\xi$, in the
following proposition, we give some properties
of $\lambda_{0}(\xi)$ as a function defined on $\xi\in (0,+\infty)$ which will be used later.
\begin{proposition}\label{youjiexing4}
    $\lambda_{0}(\xi)$ defined on $(0,+\infty)$ is continuous. And $\lambda_{0}(\xi)$
    is bounded,
    \begin{align}\label{zuidazhi1217}
        \Lambda:=\sup\limits_{\xi\in(0,+\infty)}\lambda_{0}(\xi)>0.
    \end{align}
\end{proposition}
\begin{proof}
 Let $f(s,\xi)=\Phi(s,\xi)+s^2$. Then
 $f(s,\xi)$ is increasing with respect to $s$ (see Proposition \ref{daodiaoxing1215}) and continuous on $(0,+\infty)
 \times (0,+\infty)$. Moreover, from Proposition \ref{jiecunzaijiuhao}, we have $f(\lambda_{0}(\xi),\xi)=0.$

 $\forall~ \xi_{0}\in (0,+\infty)$ and $\forall~ \varepsilon>0$. Since $f(\lambda_{0}(\xi_{0}),
 \xi_{0})=0$ and $f$ is increasing about $s$, we have
 \begin{align*}
    f(\lambda_{0}(\xi_{0})-\varepsilon,\xi_{0})<0,~f(\lambda_{0}(\xi_{0})+\varepsilon,\xi_{0})>0.
 \end{align*}
Then due to the continuity of $f(s,\xi)$ (see Proposition \ref{lianxuxing1215}),
there exists a $\delta>0$ such that
when $\xi\in O(\xi_{0},\delta)$,
\begin{align*}
    f(\lambda_{0}(\xi_{0})-\varepsilon,\xi)<0,~f(\lambda_{0}(\xi_{0})+\varepsilon,\xi)>0.
 \end{align*}
Then according to the monotonicity of $f$ and $f(\lambda_{0}(\xi),\xi)=0$, we have
\begin{align*}
    |\lambda_{0}(\xi)-\lambda_{0}(\xi_{0})|<\varepsilon,
\end{align*}
which indicates that $\lambda_{0}(\xi)$ is continuous.

In view of Proposition \ref{cunzai1215} and \eqref{lingling}, one can obtain that $\lambda_{0}(\xi)$ is bounded and
$
\lim\limits_{\xi\rightarrow\infty}\lambda_{0}(\xi)=0.
$ We complete the proof of this proposition.
\end{proof}

\subsection{Construction of a solution to system (4.5)}

Grounded on the above discussion, we get that the equation \eqref{equ1}-\eqref{mofangxietisheng} is solved.
Thus, for the solution of equation \eqref{tezhengfangcheng2}, we have the following conclusion.
\begin{theorem}\label{buheshijiuladao}
    Let $\xi\neq 0$ be fixed. Then the equation \eqref{tezhengfangcheng2} has a solution
    $(U_{1},U_{2},\varpi,h)=(U_{1}(\xi,x_{2}),$
    $U_{2}(\xi,x_{2}),\varpi(\xi,x_{2}),h(\xi,x_{2}))$,  which
    belongs to $H^{k}([0,1])$ for any positive integer $k$, with $\lambda_{0}=\lambda_{0}(\xi)>0$.
\end{theorem}

\begin{proof}
    With the help of Theorem \ref{jiecunzai1215}, we first prove that there exists a solution $(U_{2},\lambda_{0})$
    satisfying \eqref{equ1}-\eqref{mofangxietisheng}.
Afterwards, we introduce
    \begin{align}\label{hanshu}
        U_{1}(\xi,x_{2})=-\frac{DU_{2}}{\xi},~ h(\xi,x_{2})=-\frac{DU_{2}}{\lambda_{0}},~
        ~\varpi(\xi,x_{2})=-\frac{\mu\left(D^2-\xi^2\right)U_{1}-\lambda_{0}U_{1}}{\xi}.
    \end{align}
It is easy to verify that $(U_{1},U_{2},h,\varpi)$ satisfies the equations \eqref{tezhengfangcheng2}-\eqref{haimeichulai1215}.

 Due to the smoothness of $U_{2}$, the functions in \eqref{hanshu} are also smooth
with respect to $x_{2}$.
\end{proof}

\begin{remark}\label{henhaohenhao}
    By observing the coefficient of  equation \eqref{equ1}, we can deduce that
     the solution constructed in Theorem \ref{buheshijiuladao} has the following
     conclusions:

   $(1)$ $U_{2}(\xi,x_{2})$, $h(\xi,x_{2})$ and $\varpi(\xi,x_{2})$ are even on $\xi$,
       when $x_{2}$ is fixed;

 $(2)$ $U_{1}(\xi,x_{2})$ is odd on $\xi$, when $x_{2}$ is fixed.
\end{remark}
Next, we consider the boundedness of the solutions constructed in above theorem.
\begin{lemma}\label{henzhongyaodeyinli}
    Let $0<R_{1}<|\xi|<R_{2}$, $U_{1}(x_{2},\xi)$, $U_{2}(x_{2},\xi)$,
    $\varpi(x_{2},\xi)$ and $h(x_{2},\xi)$ be constructed
    as in Theorem \ref{buheshijiuladao}. Then there exist positive
    constants $C_{k}$ depending on $k$, $R_{1}$, $R_{2}$, $\overline{\theta}$,
     such that,
    \begin{align}\label{juanqilai}
        \big{\|}(U_{1},U_{2},\varpi,h)(x_{2},\xi)\big{\|}_{H^k([0,1])}\leq C_{k},~\text{for~any~integer}~k\geq 0.
    \end{align}
    In addition, we also conclude that
    \begin{align}\label{xianeryijian}
        \big{\|}U_{2}(x_{2},\xi)\big{\|}_{L^2([0,1])}^2>0.
    \end{align}
\end{lemma}
\begin{proof}
Since $U_{2}(x_{2},\xi)\in \mathcal{A}_{\xi},$ it is easy to see that  \eqref{xianeryijian} holds and there exists a constant
$C(R_{1},R_{2})$ depending on $R_{1}$ and $R_{2}$ such that
\begin{align}
    \big{\|}U_{2}(x_{2},\xi)\big{\|}_{H^1([0,1])}\leq C(R_{1},R_{2}).
\end{align}
On the other hand, by virtue of Proposition \ref{youjiexing4} and $-\lambda_{0}^2=\Phi(\lambda_{0},\xi)$, we have
\begin{align*}
    0< C(R_{1},R_{2})\leq\lambda_{0}(\xi)\leq \Lambda,~~\text{for~any~}|\xi|\in(R_{1},R_{2}).
\end{align*}

Moreover, using the following equality,
$$-\lambda_{0}^2(\xi)=\lambda_{0}(\xi)E_{0}(U_{2}(x_{2},\xi))+\lambda_{0}(\xi)E_{1}(U_{2}(
    x_{2},\xi))+E_{2}(U_{2}(x_{2},\xi)),~~\text{for}~
U_{2}(x_{2},\xi)\in\mathcal{A}_{\xi},$$
we can obtain that
\begin{align*}
    \big{\|}U_{2}(x_{2},\xi)\big{\|}_{H^2([0,1])}\leq C\left(R_{1},R_{2},\overline{\theta}\right).
\end{align*}
According to \eqref{hanshu}, a direct calculation gives that
\begin{align*}
    \begin{aligned}
    \big{\|}U_{1}(x_{2},\xi)\big{\|}_{H^1([0,1])}\leq &C\left(R_{1},R_{2},\overline{\theta}\right),~
    \big{\|}h(x_{2},\xi)\big{\|}_{H^2([0,1])}\leq C\left(R_{1},R_{2},\overline{\theta}\right),
    \\
    &\big{\|}\varpi(x_{2},\xi)\big{\|}_{L^2([0,1])}\leq  C\left(R_{1},R_{2},\overline{\theta}\right).
    \end{aligned}
\end{align*}
Finally, by iterative application of \eqref{equ1} and  \eqref{hanshu},
we derive that \eqref{juanqilai} holds.
\end{proof}
\subsection{Exponential growth rate}
Now, we use the Fourier synthesis to obtain growing solutions to \eqref{naiverboundarycondition123}-\eqref{xianxingpart} by virtue of the solutions constructed in Theorem \ref{buheshijiuladao} for any fixed spatial frequency
$\xi\in\mathbf{R}$ with $|\xi|>0$. The constructed solutions will grow in-time under the norm of $H^k$ for any positive integer $k$.
\begin{theorem}\label{gouzaojie}
    Let $f\in C_{0}^{\infty}(0,+\infty)$ be a real-valued function, $U_{2}=W_{0}$ in Proposition
    \ref{jiecunzaizhengming}, $U_{1}$, $\varpi$ and $h$ in \eqref{hanshu}.
    Define the functions
    \begin{align}\label{xianxingjie}
        \begin{aligned}
            &v_{1}(x_{1},x_{2},t)=-\frac{1}{2\pi}\int_{R}
            f(|\xi|)iU_{1}(\xi,x_{2})e^{i\xi x_{1}}e^{\lambda_{0}(\xi)t} d\xi,
            \\
            &v_{2}(x_{1},x_{2},t)=\frac{1}{2\pi}\int_{R}
            f(|\xi|)U_{2}(\xi,x_{2})e^{i\xi x_{1}}e^{\lambda_{0}(\xi)t} d\xi,
            \\
            &\pi(x_{1},x_{2},t)=\frac{1}{2\pi}\int_{R}
            f(|\xi|)\varpi(\xi,x_{2})e^{i\xi x_{1}}e^{\lambda_{0}(\xi)t} d\xi,
            \\
            &\Theta(x_{1},x_{2},t)=\frac{1}{2\pi}\int_{R}
            f(|\xi|)h(\xi,x_{2})e^{i\xi x_{1}}e^{\lambda_{0}(\xi)t} d\xi,
        \end{aligned}
    \end{align}
    then we have the following conclusions,

        $(1)$~
     $(\mathbf{v}=(v_{1},v_{2}),\pi,\Theta)$ is a solution to linearized problem
    \eqref{naiverboundarycondition123}-\eqref{xianxingpart};

    $(2)$~
    due to the smoothness of functions
    $U_{1}$, $U_{2}$, $\varpi$ and $h$, the following inequality holds
    \begin{align}\label{yinguoxunhuan}
        \big{\|}(\mathbf{v},\pi,\Theta)(0)\big{\|}_{H^k(\Omega_{\infty})}\leq \widehat{C}
        \left(\int_{\mathbf{R}}\left(1+\xi^2\right)^k |f(|\xi|)|^2d\xi\right)^{\frac{1}{2}}<+\infty,~k\in\mathbf{N},
    \end{align}
    where $\widehat{C}$ is a positive constant depending on $k_{0}$, $k_{1}$, $\mu$,
    $\overline{\theta}$ and $k$;

    $(3)$~ for every $t>0,$ the boundedness of $\lambda_{0}(\xi)$ over $(0,+\infty)$ implies
    that the solution \eqref{xianxingjie} $\in H^k$ and satisfies
    \begin{align}\label{xianranyoujie}
        e^{\lambda_{f}t}\big{\|}(\mathbf{v},\pi,\Theta)(0)\big{\|}_{H^k(\Omega_{\infty})}\leq
        \big{\|}(\mathbf{v},\pi,\Theta)(t)\big{\|}_{H^k(\Omega_{\infty})}
        \leq  e^{\Lambda t}\big{\|}(\mathbf{v},\pi,\Theta)(0)\big{\|}_{H^k(\Omega_{\infty})},
    \end{align}
    where
$\lambda_{f}:=\inf\limits_{|\xi|\in supp(f)}\lambda_{0}(\xi)>0,$
    and $\Lambda$ is a positive number defined in Proposition \ref{youjiexing4};

$(4)$~ if~$f\nequiv 0$, we have
    \begin{align}\label{jinghao}
       \big{\|}v_{2}(0)\big{\|}_{H^k(\Omega_{\infty})}>0;
    \end{align}

    $(5)$~ we can choose some $f$ such that
    \begin{align}\label{yesruci}
        \lambda_{f}=\frac{\Lambda}{2}.
    \end{align}
\end{theorem}
\begin{proof}
    For each fixed $\xi\neq 0$, the following functions
    \begin{align*}
        \begin{aligned}
            &\widetilde{v}_{1}=-f(|\xi|)iU_{1}(\xi,x_{2})e^{i\xi x_{1}}e^{\lambda_{0}(\xi)t},
            \\
            &\widetilde{v}_{2}=f(|\xi|)U_{2}(\xi,x_{2})e^{i\xi x_{1}}e^{\lambda_{0}(\xi)t},
            \\
            &\widetilde{q}=f(|\xi|)\varpi(\xi,x_{2})e^{i\xi x_{1}}e^{\lambda_{0}(\xi)t},
            \\
            &\widetilde{\Theta}=f(|\xi|)h(\xi,x_{2})e^{i\xi x_{1}}e^{\lambda_{0}(\xi)t},
        \end{aligned}
    \end{align*}
solve the linear equation \eqref{xianxingpart}. Since $f\in C_{0}^{\infty}(0,+\infty)$,
Lemma \ref{henzhongyaodeyinli} implies that
\begin{align}\label{zhegezhege}
    \sup\limits_{\xi\in \text{supp(f)}}\left\|\partial_{x_{2}}^{k}(\widetilde{v}_{1},\widetilde{v}_{2}
    ,\widetilde{q},\widetilde{\Theta})
    (\xi,x_{2})\right\|_{L^{2}(\Omega_{\infty})}<\infty,~~
    \text{for~all}~k\in \mathbf{N},
\end{align}
which together with the Sobolev embedding inequality indicates
\begin{align*}
\sup\limits_{\xi\in \text{supp}(f)}
\left|(\widetilde{v}_{1},\widetilde{v}_{2}
    ,\widetilde{q},\widetilde{\Theta})\right|\leq C_{S}\sup\limits_{\xi\in \text{supp(f)}}\left\|\partial_{x_{2}}^{k}(\widetilde{v}_{1},\widetilde{v}_{2}
    ,\widetilde{q},\widetilde{\Theta})
    (\xi,x_{2})\right\|_{L^{2}(\Omega_{\infty})}<\infty,
\end{align*}
where $C_{S}>0$ is the embedding constant.
Then the Fourier synthesis of the solution given by \eqref{xianxingjie} is also a solution to
\eqref{xianxingpart}. Moreover,  since $f(|\xi|)$ is real-valued and radial, we have from Remark \ref{henhaohenhao}  that the solution in \eqref{xianxingjie}
is real-valued.

A direct calculation gives
\begin{align*}
    \begin{aligned}
        \frac{\partial^{k} v_{1}}{\partial x_{1}^{j}\partial_{x_{2}}^{k-j}}
        =-\frac{1}{2\pi}
        \int_{\mathbf{R}}f(|\xi|)\partial^{k-j}_{x_{2}}(iU_{1}(\xi,x_{2}))(i\xi)^j e^{i\xi x_{1}}
        e^{\lambda_{0}(\xi)t} d\xi.
    \end{aligned}
\end{align*}
According to Parseval equality and \eqref{zhegezhege}, we have
\begin{align*}
    \begin{aligned}
\bigg{\|}\frac{\partial^{k} v_{1}}{\partial x_{1}^{j}\partial x_{2}^{k-j}}\bigg{\|}_{L^2(\Omega_{\infty})}^2&=
2\pi \left\|f(|\xi|)\partial^{k-j}_{x_{2}}(iU_{1}(\xi,x_{2}))(i\xi)^j
e^{\lambda_{0}(\xi)t}\right\|_{L^2(\Omega_{\infty})}^2
\\
&\leq \widehat{C}\int_{\mathbf{R}}\left(1+|\xi|^2\right)^j f^2(|\xi|)e^{2\lambda_{0}(\xi)t} d\xi.
    \end{aligned}
\end{align*}
Thus, \eqref{yinguoxunhuan} is valid.

According to Proposition \ref{youjiexing4}, it follows that \eqref{xianranyoujie} holds.
From \eqref{xianeryijian}, it is easy to see that \eqref{jinghao} is obviously valid.
 And \eqref{yesruci} comes from the continuity of $\lambda_{0}$
with respect to $\xi$.
Thus, we end the proof of this  theorem.
\end{proof}
\begin{remark}
    The constant $\Lambda$ is called maximal linear growth rate. And from Remark \ref{zongshiruci}, we can see that
    $\Lambda<\infty$ and $\Lambda\rightarrow 0$ as $\mu\rightarrow\infty$.
    This shows that the viscosity plays a stabilizing role in the linear instability.
\end{remark}

\section{Nonlinear instability}\label{feixian0331}
The goal of this section is to  establish the RT instability for nonlinear perturbed equation \eqref{henyoupinwei}-
\eqref{naiverboundarycondition123}, which can be stated in following theorem.
\begin{theorem}[Nonlinear instability]\label{feixianxingbuwending1217}
    Assume that the steady temperature profile $\overline{\theta}$ satisfies \eqref{cond1}, then
    the steady state $\left(\mathbf{0},\overline{p},\overline{\theta}\right)$ of
    \eqref{henyoupinwei}-\eqref{naiverboundarycondition123} is unstable in the following sense. That is,
    for any $s\geq 2$, $\delta>0$, $K>0$ and $F$ satisfying
    \begin{align}\label{quzheyangdehanshu1217}
        F(y)\leq Ky,~~\text{for~any~}y\in[0,\infty),
    \end{align}
then the unique strong solution $(\mathbf{v},\Theta)$ to nonlinear problem \eqref{henyoupinwei}-\eqref{naiverboundarycondition123}
    with some smooth initial value $ (\mathbf{v}_{0},\Theta_{0})\in \left(H^{\infty}(\Omega_{\infty})\right)^3$ and $\big{\|}(\mathbf{v}_{0},\Theta_{0})\big{\|}_{H^{s}(\Omega_{\infty})}
        :=\sqrt{\big{\|}\mathbf{v}_{0}\big{\|}_{H^s(\Omega_{\infty})}^2
        +\big{\|}\Theta_{0}\big{\|}_{H^{s-1}(\Omega_{\infty})}^2}<\delta$  satisfies
\begin{align}\label{feixianxinginstability1217}
        \big{\|}v_{2}(t_{K})\big{\|}_{L^2(\Omega_{\infty})}>F\left(\big{\|}(\mathbf{v}_{0},\Theta_{0})\big{\|}_{H^s(\Omega_{\infty})}\right),
    ~\text{for}~t_{K}\in \left(0,\frac{2}{\Lambda}\ln{\frac{2K}{i_{0}}}\right],
    \end{align}
    where $i_{0}$ is a constant depending on $s$,  $\Lambda$ is given by \eqref{zuidazhi1217}
    and $H^{\infty}=\cap_{k=1}^{\infty}H^{k}$.
\end{theorem}
\begin{remark}
The main steps of the proof to Theorem \ref{feixianxingbuwending1217} are outlined as follows:

$(1)$ Verify the uniqueness of the strong solution to linearized perturbed problem
\eqref{naiverboundarycondition123}-\eqref{xianxingpart}, as detailed in Theorem \ref{weiyi}.

$(2)$ Establish the estimate of the strong solution to  nonlinear problem
\eqref{henyoupinwei}-\eqref{naiverboundarycondition123}, as presented in Proposition \ref{feichangzhongyaode}.

$(3)$ Utilize above estimate to construct a family of solutions for
nonlinear problem \eqref{henyoupinwei}-\eqref{naiverboundarycondition123}. Verify that the limits of this family of solutions correspond to the solution
of linearized problem \eqref{naiverboundarycondition123}-\eqref{xianxingpart}, leveraging the uniqueness Theorem \ref{weiyi}.

$(4)$ Exploit the exponential growth rate $($see \eqref{duideqiziji}$)$ of solution to
linearized problem \eqref{naiverboundarycondition123}-\eqref{xianxingpart} and the estimate $($see Proposition \ref{feichangzhongyaode}$)$ to derive contradiction, as showed in
Lemma \ref{maodunchengli0105}.
\end{remark}
\subsection{Uniqueness of linearized equations (1.13)}
Here, we shall show the uniqueness of solution to linearized problem \eqref{xianxingpart}, which will
be used in the proof of nonlinear instability. From the above linear instability,
one can conclude the existence of solution to linearized problem
\eqref{xianxingpart} under cretain initial value.
For convenience, we define the function space of strong solutions to linearized problem as follows,
\begin{align*}
    \begin{aligned}
    \mathcal{Q}(T)=\bigg{\{}(\mathbf{v},\pi,\Theta)\big{|}&\Theta\in C\left([0,T],L^2(\Omega_{\infty})\right),
    \nabla \pi\in L^\infty\left(0,T;(L^2(\Omega_{\infty}))^2\right),
    \\
    &\mathbf{v}\in C\left([0,T],\widehat{L}^2(\Omega_{\infty})\right)
    \cap L^\infty\left(0,T;\mathbf{V}^{1,2}(\Omega_{\infty})\right),\partial_{t}\mathbf{v}\in \left(L^2((0,T)\times\Omega_{\infty})\right)^2\bigg{\}}.
    \end{aligned}
\end{align*}

\begin{theorem}[Uniqueness]\label{weiyi} Assume that $(\mathbf{v},\pi,\Theta)\in\mathcal{Q}(T)$
    is a strong solution to \eqref{xianxingpart} with $(\mathbf{v},\Theta)(0)=((0,0),0)$.
    Then, $(\mathbf{v},\nabla \pi,\Theta)=((0,0),(0,0),0).$
\end{theorem}
\begin{proof}
   Multiplying $\eqref{xianxingpart}_{1}$ and $\eqref{xianxingpart}_{2}$ with $\mathbf{v}$ and  $\theta$, respectively, then integrating them over $\Omega_{\infty}$ gives that
    \begin{align*}
        \begin{aligned}
            &\frac{1}{2}\frac{d}{dt}\big{\|}\mathbf{v}\big{\|}_{L^2(\Omega_{\infty})}^2+\mu\big{\|}\nabla\mathbf{v}\big{\|}_{L^2(\Omega_{\infty})}^2
            =k_{1}\int_{\mathbf{R}}(v_{1}(x_{1},1))^2dx_{1}+
            k_{0}\int_{\mathbf{R}}(v_{1}(x_{1},0))^2dx_{1}+\int_{\Omega}g\Theta v_{2}dx_{1}dx_{2},
            \\
            &\frac{1}{2}\frac{d}{dt}\big{\|}\Theta\big{\|}_{L^2(\Omega_{\infty})}^2=-\int_{\Omega}D\overline{\theta}
            \Theta v_{2}dx_{1}dx_{2}.
        \end{aligned}
    \end{align*}
Since $\max\{k_{1},k_{0}\}\leq 0$, then we have
\begin{align*}
    \begin{aligned}
        \frac{1}{2}\frac{d}{dt}\big{\|}\mathbf{v}\big{\|}_{L^2(\Omega_{\infty})}^2\leq \int_{\Omega}g\Theta v_{2}dx_{1}dx_{2},
        ~~~~
        \frac{1}{2}\frac{d}{dt}\big{\|}\Theta\big{\|}_{L^2(\Omega_{\infty})}^2=-\int_{\Omega}D\overline{\theta}
        \Theta v_{2}dx_{1}dx_{2},
    \end{aligned}
\end{align*}
which implies that
\begin{align*}
    \frac{d}{dt}\left[\big{\|}\mathbf{v}\big{\|}_{L^2(\Omega_{\infty})}^2+\big{\|}\Theta\big{\|}_{L^2(\Omega_{\infty})}^2\right]\leq C\left(\overline{\theta}\right)
    \left[\big{\|}\mathbf{v}\big{\|}_{L^2(\Omega_{\infty})}^2+\big{\|}\Theta\big{\|}_{L^2(\Omega_{\infty})}^2\right].
\end{align*}
By the Gronwall's inequality, we have
\begin{align*}
\big{\|}\mathbf{v}\big{\|}_{L^2(\Omega_{\infty})}^2+\big{\|}\Theta\big{\|}_{L^2(\Omega_{\infty})}^2=0,~~\text{for~any~}t\in[0,T].
\end{align*}
Thus, the proof is completed.
\end{proof}
\subsection{Energy estimates of perturbed problem (1.10)-(1.11)}
This subsection investigates some energy estimates for the perturbed problem, which
will be useful in the proof of nonlinear instability. For this purpose, we assume that
$(\mathbf{v},\pi,\Theta)$ is a strong solution of nonlinear problem \eqref{henyoupinwei}-\eqref{naiverboundarycondition123}. In previous subsection \ref{20240319}, we have established Theorem \ref{qiangjie1}
about the existence of strong solution to problem \eqref{model}-\eqref{boundarycondition}. Thus,
$(\mathbf{v},\pi,\Theta
)=\left(\mathbf{u},p-\overline{p},\theta-\overline{\theta}\right)$ should be the strong solution to
\eqref{henyoupinwei}-\eqref{naiverboundarycondition123} if $(\mathbf{u},p,\theta)$ is a strong solution of
\eqref{model}-\eqref{boundarycondition}.

\begin{proposition}\label{feichangzhongyaode}
    There exists a $\delta\in (0,1)$ such that
    if $\big{\|}\Theta_{0}\big{\|}_{H^1(\Omega_{\infty})}^2+\big{\|}\mathbf{v}_{0}\big{\|}_{H^2(\Omega_{\infty})}^2=\delta_{0}^2<\delta^2$, where
    $\Theta_{0}$ and $\mathbf{v}_{0}$ are the initial values to problem \eqref{henyoupinwei}-\eqref{naiverboundarycondition123}, then any strong solution
    $(\mathbf{v},\pi,\Theta)$ to problem \eqref{henyoupinwei}-\eqref{naiverboundarycondition123},
    emanating from $\Theta_{0}$ and $\mathbf{v}_{0}$, satisfies
    \begin{align*}
        \begin{aligned}
            \sup\limits_{t\in[0,T]}&\left[\big{\|}\mathbf{v}\big{\|}_{H^2(\Omega_{\infty})}^2+\big{\|}\Theta\big{\|}_{H^1(\Omega_{\infty})}^2+\big{\|}\nabla \pi\big{\|}_{L^2(\Omega_{\infty})}^2
+\big{\|}\mathbf{v}_{t}\big{\|}_{L^2(\Omega_{\infty})}^2+\big{\|}\Theta_{t}\big{\|}_{L^2(\Omega_{\infty})}^2\right]
\\
&+\int_{0}^{t}\left[\big{\|}\mathbf{v}_{s}\big{\|}_{H^1(\Omega_{\infty})}^2+\big{\|}\nabla\mathbf{v}\big{\|}_{H^1(\Omega_{\infty})}^2\right] ds
            \leq C(T)\delta_{0}^2,
        \end{aligned}
    \end{align*}
    where $t\in [0,T]$ and $T>0$ is a arbitrarily fixed constant.
    \end{proposition}
    \begin{proof}
        $(1)$~ Estimations of $\big{\|}\Theta\big{\|}_{L^2(\Omega_{\infty})}$, $\big{\|}\mathbf{v}\big{\|}_{L^2(\Omega_{\infty})}$ and
        $\int_{0}^{t}\big{\|}\nabla\mathbf{v}\big{\|}_{L^2(\Omega_{\infty})}^2ds$.

        Multiplying $\eqref{henyoupinwei}_{2}$ by
        $\Theta$, using $\eqref{henyoupinwei}_{3}$ and integrating over $\Omega_{\infty}$ gives that
        \begin{align*}
            \frac{1}{2}\frac{d}{dt}\big{\|}\Theta\big{\|}_{L^2(\Omega_{\infty})}^2=-\int_{\Omega_{\infty}}
            v_{2}D\overline{\theta}\Theta dx_{1}dx_{2}
            \leq\frac{\big{\|}D\overline{\theta}\big{\|}_{L^\infty(\Omega_{\infty})}}{2}\left(\big{\|}\mathbf{v}\big{\|}_{L^2(\Omega_{\infty})}^2+\big{\|}\Theta\big{\|}_{L^2(\Omega_{\infty})}^2\right).
        \end{align*}
        And multiplying $\eqref{henyoupinwei}_{1}$ by $\mathbf{v}$ and integrating over $\Omega_{\infty}$,
        we obtain
        \begin{align*}
            \begin{aligned}
            \frac{1}{2}\frac{d}{dt}\big{\|}\mathbf{v}\big{\|}_{L^2(\Omega_{\infty})}^2&+\mu\big{\|}\nabla\mathbf{v}\big{\|}_{L^2(\Omega_{\infty})}^2
            -\int_{\mathbf{R}}\left[k_{1}(v_{1}(x_{1},1))^2
            +k_{0}(v_{1}(x_{1},0))^2\right]dx_{1}
            \\
            &=\int_{\Omega_{\infty}}\Theta \mathbf{g}\cdot \mathbf{v} dx_{1}dx_{2}
            \leq\frac{g}{2}\left(\big{\|}\Theta\big{\|}_{L^2(\Omega_{\infty})}^2+\big{\|}\mathbf{v}\big{\|}_{L^2(\Omega_{\infty})}^2\right).
            \end{aligned}
        \end{align*}
        Since $\max\{k_{1},k_{0}\}\leq 0$, we have
        \begin{align}\label{zhijie}
            \frac{1}{2}\frac{d}{dt}\left[\big{\|}\Theta\big{\|}_{L^2(\Omega_{\infty})}^2+\big{\|}\mathbf{v}\big{\|}_{L^2(\Omega_{\infty})}^2\right]+\mu\big{\|}\nabla\mathbf{v}\big{\|}_{L^2(\Omega_{\infty})}
            \leq\frac{g+\big{\|}D\overline{\theta}\big{\|}_{L^\infty(\Omega_{\infty})}}{2}\left(\big{\|}\Theta\big{\|}_{L^2(\Omega_{\infty})}^2+\big{\|}\mathbf{v}\big{\|}_{L^2(\Omega_{\infty})}^2\right).
        \end{align}
    Then, according to the Gronwall's inequality and taking $C_{1}=g+\|D\overline{\theta}\|_{L^\infty}$, we have
        \begin{align}\label{guji1}
            \frac{1}{2}\left[\big{\|}\Theta\big{\|}_{L^2(\Omega_{\infty})}^2+\big{\|}\mathbf{v}\big{\|}_{L^2(\Omega_{\infty})}^2\right]+\mu\int_{0}^{t}\big{\|}\nabla\mathbf{v}\big{\|}_{L^2(\Omega_{\infty})}^2ds
            \leq\frac{1}{2}e^{C_{1}t}\left(\big{\|}\Theta_{0}\big{\|}_{L^2(\Omega_{\infty})}^2
            +\big{\|}\mathbf{v}_{0}\big{\|}_{L^2(\Omega_{\infty})}^2\right).
        \end{align}
        Thus, we also have
        \begin{align}\label{houmian1}
            \begin{aligned}
            \big{\|}\Theta\big{\|}_{L^2(\Omega_{\infty})}^2+\big{\|}\mathbf{v}\big{\|}_{L^2(\Omega_{\infty})}^2&+2\mu\int_{0}^t \big{\|}\nabla\mathbf{v}\big{\|}_{L^2(\Omega_{\infty})}^2 ds
           \leq C(T)\left[\big{\|}\Theta_{0}\big{\|}_{L^2(\Omega_{\infty})}^2+\big{\|}\mathbf{v}_{0}\big{\|}_{L^2(\Omega_{\infty})}^2\right]
            <C(T)\delta_{0}^2.
            \end{aligned}
        \end{align}

        $(2)$~Estimations of $\big{\|}\mathbf{v}_{t}\big{\|}_{L^2(\Omega_{\infty})}^2$ and $\int_{0}^{t}\big{\|}\nabla\mathbf{v}_{s}\big{\|}_{L^2(\Omega_{\infty})}^2 ds.$

        To control $\mathbf{v}_{t}$, multiplying  $\eqref{henyoupinwei}_{1}$ by $\mathbf{v}_{t}$ in
        $L^2$-inner product, then we obtain
        \begin{align*}
            \begin{aligned}
               \big{\|}\mathbf{v}_{t}\big{\|}_{L^2(\Omega_{\infty})}^2&=-\int_{\Omega_{\infty}}
               (\mathbf{v}\cdot\nabla)\mathbf{v}\cdot\mathbf{v}_{t}
               dx_{1}dx_{2}+\mu\int_{\Omega_{\infty}}\Delta\mathbf{v}\cdot\mathbf{v}_{t} dx_{1}dx_{2}
               +\int_{\Omega_{\infty}}\Theta\mathbf{g}\cdot\mathbf{v}_{t}dx_{1}dx_{2}
               \\
               &\leq \big{\|}(\mathbf{v}\cdot\nabla)\mathbf{v}\big{\|}_{L^2}\big{\|}\mathbf{v}_{t}\big{\|}_{L^2(\Omega_{\infty})}
               +\mu\big{\|}\Delta\mathbf{v}\big{\|}_{L^2(\Omega_{\infty})}\big{\|}\mathbf{v}_{t}\big{\|}_{L^2(\Omega_{\infty})}
               +g\big{\|}\Theta\big{\|}_{L^2(\Omega_{\infty})}\big{\|}\mathbf{v}_{t}\big{\|}_{L^2(\Omega_{\infty})}
               \\
               &\leq \big{\|}(\mathbf{v}\cdot\nabla)\mathbf{v}\big{\|}_{L^2(\Omega_{\infty})}^2+
               \frac{3}{4}\big{\|}\mathbf{v}_{t}\big{\|}_{L^2(\Omega_{\infty})}^2+\mu^2\big{\|}\Delta\mathbf{v}\big{\|}_{L^2(\Omega_{\infty})}^2
               +g^2\big{\|}\Theta\big{\|}_{L^2(\Omega_{\infty})}^2
               \\
               &\leq \big{\|}\mathbf{v}\big{\|}_{L^\infty}^2\big{\|}\nabla\mathbf{v}\big{\|}_{L^2(\Omega_{\infty})}^2
               +\frac{3}{4}\big{\|}\mathbf{v}_{t}\big{\|}_{L^2(\Omega_{\infty})}^2+\mu^2\big{\|}\Delta\mathbf{v}\big{\|}_{L^2(\Omega_{\infty})}^2
               +g^2\big{\|}\Theta\big{\|}_{L^2(\Omega_{\infty})}^2
               \\
               &\leq C\big{\|}\mathbf{v}\big{\|}_{L^2(\Omega_{\infty})}\big{\|}\mathbf{v}\big{\|}_{H^2(\Omega_{\infty})}\big{\|}\nabla\mathbf{v}\big{\|}_{L^2(\Omega_{\infty})}^2
               +\frac{3}{4}\big{\|}\mathbf{v}_{t}\big{\|}_{L^2(\Omega_{\infty})}^2+\mu^2\big{\|}\Delta\mathbf{v}\big{\|}_{L^2(\Omega_{\infty})}^2
               +g^2\big{\|}\Theta\big{\|}_{L^2(\Omega_{\infty})}^2,
            \end{aligned}
        \end{align*}
        where the Cauchy-Schwarz's inequality and Lemma \ref{Linftyguji} are used.
        Thus, one can obtain
        \begin{align}\label{chushizhi}
            \begin{aligned}
            \big{\|}\mathbf{v}_{t}(0)\big{\|}_{L^2(\Omega_{\infty})}^2&\leq C\left(\big{\|}\mathbf{v}_{0}\big{\|}_{H^2(\Omega_{\infty})}^4+
            \big{\|}\mathbf{v}_{0}\big{\|}_{H^2(\Omega_{\infty})}^2\right)+g^2\big{\|}\Theta_{0}\big{\|}_{L^2(\Omega_{\infty})}^2
            \leq C(T)\delta_{0}^2.
            \end{aligned}
        \end{align}
    Then, differentiating $\eqref{henyoupinwei}_{1}$ with respect to $t$, multiplying
    the result by $\mathbf{v}_{t}$ and integrating over $\Omega_{\infty}$, one has
        \begin{align}\label{zongshizheyang}
            \begin{aligned}
                \frac{1}{2}&\frac{d}{dt}\big{\|}\mathbf{v}_{t}\big{\|}_{L^2(\Omega_{\infty})}^2+\mu\big{\|}\nabla\mathbf{v}_{t}\big{\|}_{L^2(\Omega_{\infty})}^2
               -\int_{\mathbf{R}}\left[k_{1}(v_{1t}(x_{1},1))^2+k_{0}(v_{1t}(x_{1},0))^2\right]dx_{1}
                \\
                &=-\int_{\Omega_{\infty}}(\mathbf{v}_{t}\cdot\nabla)\mathbf{v}\cdot\mathbf{v}_{t} dx
                +\int_{\Omega_{\infty}}\Theta_{t}\mathbf{g}\cdot\mathbf{v}_{t}dx
                \\
                &=\int_{\Omega_{\infty}}(\mathbf{v}_{t}\cdot\nabla)\mathbf{v}_{t}\cdot\mathbf{v} dx
                +\int_{\Omega_{\infty}}(\mathbf{v}\cdot\nabla)\mathbf{v}_{t}\cdot \Theta \mathbf{g} dx
                -\int_{\Omega_{\infty}}v_{2}D\overline{\theta}\mathbf{g}\cdot\mathbf{v}_{t}dx
                \\
                &\leq \big{\|}\mathbf{v}_{t}\big{\|}_{L^4(\Omega_{\infty})}\big{\|}\mathbf{v}\big{\|}_{L^4(\Omega_{\infty})}\big{\|}\nabla\mathbf{v}_{t}\big{\|}_{L^2(\Omega_{\infty})}
                +g\big{\|}\mathbf{v}\big{\|}_{L^4(\Omega_{\infty})}\big{\|}\Theta\big{\|}_{L^4(\Omega_{\infty})}\big{\|}\nabla\mathbf{v}_{t}\big{\|}_{L^2(\Omega_{\infty})}
               \\ &~~~~+g\big{\|}D\overline{\theta}\big{\|}_{L^\infty(\Omega_{\infty})}\big{\|}\mathbf{v}\big{\|}_{L^2(\Omega_{\infty})}\big{\|}\mathbf{v}_{t}\big{\|}_{L^2(\Omega_{\infty})}.
            \end{aligned}
        \end{align}
    Applying Cauchy-Schwarz's inequality and Lemma \ref{L4guji},  one
    can obtain from \eqref{zongshizheyang} that
    \begin{align}\label{xuyaoni1216}
        \begin{aligned}
        \frac{d}{dt}\big{\|}\mathbf{v}_{t}\big{\|}_{L^2(\Omega_{\infty})}^2&+\mu\big{\|}\nabla\mathbf{v}_{t}\big{\|}_{L^2(\Omega_{\infty})}^2
        \leq C\left(\big{\|}\mathbf{v}\big{\|}_{L^2(\Omega_{\infty})}^2\big{\|}\nabla\mathbf{v}\big{\|}_{L^2(\Omega_{\infty})}^2
        +1\right)\big{\|}\mathbf{v}_{t}\big{\|}_{L^2(\Omega_{\infty})}^2
        \\
&+C\big{\|}\mathbf{v}\big{\|}_{L^2(\Omega_{\infty})}\big{\|}\nabla\mathbf{v}v_{L^2(\Omega_{\infty})}
        \big{\|}\Theta\big{\|}_{L^4(\Omega_{\infty})}^2+C\big{\|}\mathbf{v}\big{\|}_{L^2(\Omega_{\infty})}^2.
        \end{aligned}
    \end{align}
    Since
            $\frac{1}{4}\frac{d}{dt}\big{\|}\Theta\big{\|}_{L^4(\Omega_{\infty})}^4=-\int_{\Omega}v_{2}D\overline{\theta}\Theta^3dx_{1}dx_{2}
            \leq \big{\|}D\overline{\theta}\big{\|}_{L^\infty(\Omega_{\infty})}\big{\|}\mathbf{v}\big{\|}_{L^4(\Omega_{\infty})}\big{\|}\Theta\big{\|}_{L^4(\Omega_{\infty})}^3,
         $
         then, by \eqref{houmian1}, one has
    \begin{align}\label{woxuyaode1216}
        \|\Theta\|_{L^4(\Omega_{\infty})}\leq C(T)\delta_{0}.
    \end{align}
    Combing \eqref{houmian1} with \eqref{woxuyaode1216}, the inequality \eqref{xuyaoni1216} can be
    rewritten as follows
    \begin{align}\label{wobuxuyaokelian1216}
        \frac{d}{dt}\big{\|}\mathbf{v}_{t}\big{\|}_{L^2(\Omega_{\infty})}^2&+\mu\big{\|}\nabla\mathbf{v}_{t}\big{\|}_{L^2(\Omega_{\infty})}^2
        \leq C\left(\delta_{0}^2\big{\|}\nabla\mathbf{v}\big{\|}_{L^2(\Omega_{\infty})}^2+1\right)\big{\|}\mathbf{v}_{t}\big{\|}_{L^2(\Omega_{\infty})}^2
        +C\big{\|}\nabla\mathbf{v}\big{\|}_{L^2(\Omega_{\infty})}^2+C\delta_{0}^2.
    \end{align}
    Then by Gronwall's inequality and \eqref{chushizhi}, one obtains
    \begin{align*}
        \big{\|}\mathbf{v}_{t}\big{\|}_{L^2(\Omega_{\infty})}^2\leq C\big{\|}\mathbf{v}_{t}(0)\big{\|}_{L^2(\Omega_{\infty})}^2+C\delta_{0}^2
    \leq C(T) \delta_{0}^2.
    \end{align*}
    Thus, we conclude the estimate
        \begin{align}\label{kujin}
            \big{\|}\mathbf{v}_{t}\big{\|}_{L^2(\Omega_{\infty})}^2+\mu\int_{0}^{t}\big{\|}\nabla\mathbf{v}_{s}\big{\|}_{L^2(\Omega_{\infty})}^2 ds
             \leq C\delta_{0}^2.
        \end{align}

        $(3)$~Estimations of $\big{\|}\nabla\mathbf{v}\big{\|}_{L^2(\Omega_{\infty})}$, $\big{\|}\nabla^2\mathbf{v}\big{\|}_{L^2(\Omega_{\infty})}$
        and $\big{\|}\nabla \pi\big{\|}_{L^2(\Omega_{\infty})}$.

        From \eqref{houmian1} and \eqref{kujin}, one can derive that
        \begin{align*}
           \mathbf{v}_{t}\in L^{2}\left(0,T;\left(H^1(\Omega_{\infty})\right)^2\right),~
           \mathbf{v}\in L^2\left(0,T;\left(H^1(\Omega_{\infty})\right)^2\right).
        \end{align*}
    Thus, from Theorem 2 in the P286 of \cite{Temam1977}, \eqref{houmian1} and \eqref{kujin}, one can conclude that
    \begin{align}\label{lairi1216}
        \big{\|}\nabla\mathbf{v}\big{\|}_{L^2(\Omega_{\infty})}^2<C(T)\delta_{0}^2.
    \end{align}

   Furthermore, using the following equation
        \begin{align*}
            -\mu\Delta \mathbf{v}+\nabla \pi=-\mathbf{v}_{t}-(\mathbf{v}\cdot \nabla )
            \mathbf{v}+\Theta\mathbf{g},
        \end{align*}
   and Corollary \ref{tuilun1216}, one has
        \begin{align}\label{stokesguji1216}
            \begin{aligned}
                \left\|\nabla^2\mathbf{v}\right\|_{L^2(\Omega_{\infty})}^2+\big{\|}\nabla \pi\big{\|}_{L^2(\Omega_{\infty})}^2\leq
                C\left(\big{\|}\mathbf{v}_{t}\big{\|}_{L^2(\Omega_{\infty})}^2+\big{\|}(\mathbf{v}\cdot\nabla)\mathbf{v}\big{\|}_{L^2(\Omega_{\infty})}^2
                +\big{\|}\Theta\big{\|}_{L^2(\Omega_{\infty})}^2+\big{\|}\nabla\mathbf{v}\big{\|}_{L^2(\Omega_{\infty})}^2\right).
            \end{aligned}
        \end{align}
    By virtue of Lemma \ref{Linftyguji}, one gets
        \begin{align}\label{qujianma1216}
            \begin{aligned}
                \big{\|}(\mathbf{v}\cdot\nabla)\mathbf{v}\big{\|}_{L^2(\Omega_{\infty})}^2&\leq  \big{\|}
                \mathbf{v}\big{\|}_{L^\infty(\Omega_{\infty})}^2\big{\|}\nabla\mathbf{v}\big{\|}_{L^2(\Omega_{\infty})}^2
                \leq C\big{\|}\mathbf{v}\big{\|}_{L^2(\Omega_{\infty})}\big{\|}\mathbf{v}\big{\|}_{H^2(\Omega_{\infty})}\big{\|}\nabla\mathbf{v}\big{\|}_{L^2(\Omega_{\infty})}^2
                \\
                &\leq C\big{\|}\mathbf{v}\big{\|}_{L^2(\Omega_{\infty})}^2\big{\|}\nabla\mathbf{v}\big{\|}_{L^2(\Omega_{\infty})}^4
                +\varepsilon\left(\big{\|}\mathbf{v}\big{\|}_{H^1(\Omega_{\infty})}^2+\left\|\nabla^2\mathbf{v}\right\|_{L^2(\Omega_{\infty})}^2\right),
            \end{aligned}
        \end{align}
    which together with  \eqref{stokesguji1216}  concludes that
\begin{align}\label{ruguoyouyitian1216}
        \left\|\nabla^2\mathbf{v}\right\|_{L^2(\Omega_{\infty})}^2+\big{\|}\nabla \pi\big{\|}_{L^2(\Omega_{\infty})}^2\leq C(T)\delta_{0}^2.
    \end{align}

        $(4)$~Estimation of $\|\nabla\Theta\|_{L^2(\Omega_{\infty})}^2.$

        By the Stokes' estimate (see  Corollary \ref{tuilun1216}), we have
        \begin{align}\label{anning}
            \begin{aligned}
                \left\|\nabla^2\mathbf{v}\right\|_{L^4(\Omega_{\infty})}+\big{\|}\nabla \pi\big{\|}_{L^4(\Omega_{\infty})}
                &\leq C\left(\big{\|}\mathbf{v}_{t}\big{\|}_{L^4(\Omega_{\infty})}+\big{\|}(\mathbf{v}\cdot
                \nabla)\mathbf{v}\big{\|}_{L^4(\Omega_{\infty})}+\big{\|}\Theta\big{\|}_{L^4(\Omega_{\infty})}+\big{\|}\nabla\mathbf{v}\big{\|}_{L^4(\Omega_{\infty})}\right).
            \end{aligned}
        \end{align}
Then, from \eqref{houmian1}, \eqref{woxuyaode1216},\eqref{kujin},  \eqref{lairi1216} and \eqref{ruguoyouyitian1216}, it follows that  \begin{align}\label{jiuzheyang1216}
        \begin{aligned}
            \int_{0}^{t}\left\|\nabla^2\mathbf{v}\right\|_{L^4(\Omega_{\infty})}+\big{\|}\nabla\pi\big{\|}_{L^4(\Omega_{\infty})}ds
            &\leq C\int_{0}^{t}\left(\big{\|}\mathbf{v}_{s}\big{\|}_{L^4(\Omega_{\infty})}
            +\big{\|}\mathbf{v}\|_{L^\infty(\Omega_{\infty})}\big{\|}\nabla\mathbf{v}\big{\|}_{L^4(\Omega_{\infty})}
            \big{\|}\Theta\big{\|}_{L^4(\Omega_{\infty})}\right)ds
            \\
            &\leq
            C\int_{0}^{t}\left(\big{\|}\mathbf{v}_{s}\big{\|}_{H^1(\Omega_{\infty})}+\big{\|}\mathbf{v}\big{\|}_{H^2(\Omega_{\infty})}^2
            +\big{\|}\Theta\big{\|}_{L^4(\Omega_{\infty})}\right) ds\leq C\delta_{0}.
        \end{aligned}
    \end{align}
    Differentiating equation $\eqref{henyoupinwei}_{2}$ in $x_{i}$ $(i=1,2)$, then we have
    \begin{align}\label{xiangcheng1216}
        \begin{aligned}
            &\partial_{t}\Theta_{x_{1}}+(\mathbf{v}_{x_{1}}\cdot\nabla)\Theta+
            (\mathbf{v}\cdot\nabla)\Theta_{x_{1}}
            +(v_{2})_{x_{1}}D\overline{\theta}=0,
            \\
            &\partial_{t}\Theta_{x_{2}}+(\mathbf{v}_{x_{2}}\cdot\nabla)\Theta
            +(\mathbf{v}\cdot\nabla)\Theta_{x_{2}}
            +(v_{2})_{x_{2}}D\overline{\theta}+v_{2}D^2\overline{\theta}=0.
        \end{aligned}
    \end{align}
    Then, adding $\eqref{xiangcheng1216}_{1}\times \Theta_{x_{1}}$ to
    $\eqref{xiangcheng1216}_{2}\times \Theta_{x_{2}}$, and integrating the result in $L^2$, one has
    \begin{align}\label{guanyushijian1216}
        \begin{aligned}
        \frac{1}{2}\frac{d}{dt}\big{\|}\nabla\Theta\big{\|}_{L^2(\Omega_{\infty})}^2&=-
        \int_{\Omega_{\infty}}\sum\limits_{i=1}^2
       \left[ (\mathbf{v}_{x_{i}}\cdot\nabla)\Theta\Theta_{x_{i}}+
       (v_{2})_{x_{i}}D\overline{\theta}\Theta_{x_{i}}\right]+v_{2}D^2\overline{\theta}\Theta_{x_{2}} dx
       \\
       &\leq C\left(\big{\|}\nabla
       \mathbf{v}\big{\|}_{L^\infty(\Omega_{\infty})}\big{\|}\nabla\Theta\big{\|}_{L^2(\Omega_{\infty})}^2+\big{\|}\mathbf{v}\big{\|}_{H^1(\Omega_{\infty})}^2+
       \big{\|}\nabla\Theta\big{\|}_{L^2(\Omega_{\infty})}^2\right).
    \end{aligned}
    \end{align}

        Then, by Gronwall's inequality and \eqref{jiuzheyang1216}, we have
        \begin{align}\label{jiuzheyangba}
            \begin{aligned}
                \big{\|}\nabla\Theta\big{\|}_{L^2(\Omega_{\infty})}^2\leq C(T)\delta_{0}^2.
            \end{aligned}
        \end{align}

        $(5)$ Estimation of $\big{\|}\Theta_{t}\big{\|}_{L^2(\Omega_{\infty})}^2$.

    A series of calculation gives that
        \begin{align*}
            \begin{aligned}
                \big{\|}\Theta_{t}\big{\|}_{L^2(\Omega_{\infty})}^2&=-\int_{\Omega}(\mathbf{v}\cdot\nabla)\Theta\Theta_{t} dx_{1}dx_{2}
                -\int_{\Omega}v_{2}D\overline{\theta}\Theta_{t} dx_{1}dx_{2}
                \\
                &\leq \big{\|}(\mathbf{v}\cdot\nabla)\Theta\big{\|}_{L^2(\Omega_{\infty})}\big{\|}\Theta_{t}\big{\|}_{L^2(\Omega_{\infty})}+
                \big{\|}D\overline{\theta}\big{\|}_{L^\infty(\Omega_{\infty})}\big{\|}\mathbf{v}\big{\|}_{L^2(\Omega_{\infty})}\big{\|}\Theta_{t}\big{\|}_{L^2(\Omega_{\infty})}
                \\
                &\leq \big{\|}(\mathbf{v}\cdot\nabla)\Theta\big{\|}_{L^2(\Omega_{\infty})}^2
                +\big{\|}D\overline{\theta}\big{\|}_{L^\infty(\Omega_{\infty})}^2\big{\|}\mathbf{v}\big{\|}_{L^2(\Omega_{\infty})}^2
                +\frac{1}{2}\big{\|}\Theta_{t}\big{\|}_{L^2(\Omega_{\infty})}^2.
            \end{aligned}
        \end{align*}
    Then form the previous estimations, one has
\begin{align}\label{zilv}
            \big{\|}\Theta_{t}\big{\|}_{L^2(\Omega_{\infty})}^2\leq 2\big{\|}\mathbf{v}\big{\|}_{L^\infty(\Omega_{\infty})}^2\big{\|}\nabla\Theta\big{\|}_{L^2(\Omega_{\infty})}^2
            +2\big{\|}D\overline{\theta}\big{\|}_{L^\infty(\Omega_{\infty})}^2\big{\|}\mathbf{v}\big{\|}_{L^2(\Omega_{\infty})}^2
            \leq C\delta_{0}^2.
        \end{align}
     Therefore, In view of \eqref{houmian1}, \eqref{kujin}, \eqref{lairi1216}, \eqref{ruguoyouyitian1216},
        \eqref{jiuzheyangba} and \eqref{zilv}, there appears the relation
        \begin{align*}
        \begin{aligned}
            \sup\limits_{t\in[0,T]}&\left[\big{\|}\mathbf{v}\big{\|}_{H^2(\Omega_{\infty})}^2+\big{\|}\Theta\big{\|}_{H^1(\Omega_{\infty})}^2+\big{\|}\nabla \pi\big{\|}_{L^2(\Omega_{\infty})}^2
+\big{\|}\mathbf{v}_{t}\big{\|}_{L^2(\Omega_{\infty})}^2+\big{\|}\Theta_{t}\big{\|}_{L^2(\Omega_{\infty})}^2\right]
\\
&+\int_{0}^{t}\left[\big{\|}\mathbf{v}_{s}\big{\|}_{H^1(\Omega_{\infty})}^2+\big{\|}\nabla\mathbf{v}\big{\|}_{H^1(\Omega_{\infty})}^2\right] ds
            \leq C(T)\delta_{0}^2.
        \end{aligned}
    \end{align*}
    Thus, we complete the proof of this proposition.
        \end{proof}
        \subsection{Construction of the solution to linearized problem}
    By virtue of Theorem \ref{gouzaojie}, we can construct a classical solution $(\mathbf{v},\pi,\Theta)$ to
\eqref{naiverboundarycondition123}-\eqref{xianxingpart} satisfying the properties \eqref{xianranyoujie}-\eqref{yesruci}. Noticing that the solution $(\mathbf{v},\Theta)$ is smooth
 and $\|v_{2}(0)\|_{H^s(\Omega_{\infty})}>0$.
Here, when we define
 \begin{align*}
    \begin{aligned}
        \left(\widetilde{\mathbf{v}},\widetilde{\pi},\widetilde{\Theta}\right)=\frac{\delta_{0}(\mathbf{v},\pi,\Theta)}
        {\big{\|}(\mathbf{v},\Theta)(0)\big{\|}_{H^s(\Omega_{\infty})}}.
    \end{aligned}
 \end{align*}
It is not difficult to verify that $\left(\widetilde{\mathbf{v}},\widetilde{\pi},\widetilde{\Theta}\right)$ is still a
classical solution to \eqref{naiverboundarycondition123}-\eqref{xianxingpart} and satisfy all the
properties of $(\mathbf{v},\pi,\Theta)$. Moreover, it is easy to see that
$$\left\|\left(\widetilde{\mathbf{v}},\widetilde{\Theta}\right)(0)\right\|_{H^s(\Omega_{\infty})}=\delta_{0}.$$
For convenience, without loss of generality, we introduce
 \begin{align}\label{jiushizheyang}
    \begin{aligned}
        \big{\|}(\mathbf{v},\Theta)(0)\big{\|}_{H^s(\Omega_{\infty})}:=\sqrt{\big{\|}\mathbf{v}_{0}\big{\|}_{H^s(\Omega_{\infty})}^2+
        \big{\|}\Theta_{0}\big{\|}_{H^{s-1}(\Omega_{\infty})}^2}=\delta_{0}.
    \end{aligned}
 \end{align}
Furthermore, we denote
\begin{align*}
    \tau_{0}:=\tau(s):=\frac{\big{\|}v_{2}(0)\big{\|}_{L^2(\Omega_{\infty})}}{\big{\|}(\mathbf{v},\Theta)(0)\big{\|}_{H^s(\Omega_{\infty})}}\leq 1.
\end{align*}
Obviously, from \eqref{jinghao}, one has
$\tau_{0}>0$. Take $t_{K}=\frac{2}{\Lambda}\ln{\frac{2K}{\tau_{0}}}$ and recall
\eqref{xianranyoujie}, we have
\begin{align}\label{duideqiziji}
    \big{\|}v_{2}(t_{K})\big{\|}_{L^2(\Omega_{\infty})}\geq e^{\frac{t_{K}\Lambda}{2}}\tau_{0}\delta_{0}\geq 2 K\delta_{0}.
\end{align}
\subsection{Construction of the solution to nonlinear problem}
Subsequently, we construct a family of solutions to nonlinear problem
\eqref{henyoupinwei}-\eqref{naiverboundarycondition123}
based on the initial value \eqref{jiushizheyang}. We define an initial values as follows
\begin{align*}
    \left(\mathbf{v}_{0}^{\varepsilon},\Theta_{0}^{\varepsilon}\right):=
    \varepsilon(\mathbf{v},\Theta)(0),~~\text{for}~\varepsilon\in(0,1).
\end{align*}
Noting that
\begin{align*}
    \left(\mathbf{v}_{0}^{\varepsilon},\Theta_{0}^{\varepsilon}\right)\in \left(H^\infty(\Omega_{\infty})\right)^3,~~
    \text{and}~~\left\|\left(\mathbf{v}_{0}^{\varepsilon},\Theta_{0}^{\varepsilon}\right)\right\|_{H^s(\Omega_{\infty})}=\delta_{0}
    \varepsilon<\delta_{0}.
\end{align*}
By Proposition \ref{feichangzhongyaode}, there exists a strong solution
$(\mathbf{v}^{\varepsilon},\pi^{\varepsilon},\Theta^{\varepsilon})$ to nonlinear
problem \eqref{henyoupinwei}-\eqref{naiverboundarycondition123} with initial $(\mathbf{v}_{0}^{\varepsilon},\Theta_{0}^{\varepsilon})$ on $(0,T)
\times \Omega$ with $T>t_{K}$, which also satisfies
\begin{align}\label{2024032011}
    \begin{aligned}
        \sup\limits_{t\in[0,t_{K}]}\bigg{[}\big{\|}\mathbf{v}^{\varepsilon}(t)\big{\|}_{H^2(\Omega_{\infty})}^2
        &+\big{\|}\Theta^{\varepsilon}(t)\big{\|}_{H^1(\Omega_{\infty})}^2+\big{\|}\nabla \pi^{\varepsilon}\big{\|}_{L^2(\Omega_{\infty})}^2
+\big{\|}\mathbf{v}_{t}^{\varepsilon}\big{\|}_{L^2(\Omega_{\infty})}^2+\big{\|}\Theta_{t}^{\varepsilon}\big{\|}_{L^2(\Omega_{\infty})}^2\bigg{]}
        \\
        &+\int_{0}^{t}\left[\big{\|}\mathbf{v}_{s}^{\varepsilon}\big{\|}_{H^1(\Omega_{\infty})}^2+\big{\|}\nabla\mathbf{v}^{\varepsilon}\big{\|}_{H^1(\Omega_{\infty})}^2\right] d\tau
        \leq C(t_{K})\delta_{0}^2\varepsilon^2,
    \end{aligned}
\end{align}
where $C(t_{K})$ is independent of $\varepsilon$.
\begin{lemma}\label{maodunchengli0105}
For any $\varepsilon_{0}$ small enough, the strong solution
$(\mathbf{v}^{\varepsilon_{0}},\Theta^{\varepsilon_{0}})$ to  nonlinear problem
    \eqref{henyoupinwei}-\eqref{naiverboundarycondition123}, emanating from the initial
    value $(\mathbf{v}^{\varepsilon_{0}},\Theta^{\varepsilon_{0}})(0)$, satisfies
    $$
   \left\|v_{2}^{\varepsilon_{0}}(t_{K})\right\|_{L^2((\Omega_{\infty}))}>F\left(\big{\|}(\mathbf{v}^{\varepsilon_{0}},
    \Theta^{\varepsilon_{0}})(0)\big{\|}_{H^s(\Omega_{\infty})}\right),~~
    \text{for~~some~~}t_{K}\in \left(0,\frac{2}{\Lambda}\ln{\frac{2K}{\tau_{0}}}\right]\subset
    (0,T).$$
\end{lemma}
\begin{proof}
This lemma is proved by the method of contradiction. Suppose that for initial value $\left(\mathbf{v}^{\varepsilon_{0}},\Theta^{\varepsilon_{0}}\right)(0)$,
 the strong solution $\left(\mathbf{v}^{\varepsilon_{0}},\pi^{\varepsilon_{0}},\Theta^{\varepsilon_{0}}\right)$ of nonlinear problem  \eqref{henyoupinwei}-\eqref{naiverboundarycondition123}  satisfies
$$
\left\|v_{2}^{\varepsilon_{0}}(t)\right\|_{L^2(\Omega_{\infty})}\leq F\left(\big{\|}(\mathbf{v}^{\varepsilon_{0}},
\Theta^{\varepsilon_{0}}\right)(0)\big{\|}_{H^s(\Omega_{\infty})}),~~\forall~t\in \left(0,\frac{2}{\Lambda}\ln{\frac{2K}{\tau_{0}}}\right].
$$
Thus, according to the property of $F$, we have the inequality
\begin{align}\label{enenen}
    \left\|v_{2}^{\varepsilon_{0}}(t)\right\|_{L^2(\Omega_{\infty})}\leq K\big{\|}\left(\mathbf{v}^{\varepsilon_{0}},
    \Theta^{\varepsilon_{0}}\right)(0)\big{\|}_{H^s(\Omega_{\infty})}\leq K\delta_{0}\varepsilon_{0},~\forall
    ~t\in(0,t_{k}].
\end{align}
We set $\left(\overline{\mathbf{v}}^{\varepsilon_{0}},\overline{\pi}^{
    \varepsilon_{0}},\overline{\Theta}^{\varepsilon_{0}}\right):=
    \frac{\left(\mathbf{v}^{\varepsilon_{0}},\pi^{\varepsilon_{0}},\Theta^{\varepsilon_{0}}\right)}{\varepsilon_{0}}$.
Then it is easy to verify that $\left(\overline{\mathbf{v}}^{\varepsilon_{0}},\overline{\pi}^{
    \varepsilon_{0}},\overline{\Theta}^{\varepsilon_{0}}\right)$ satisfies
    \begin{align}\label{tianzhidao}
        \begin{cases}
            \partial_{t}\overline{\mathbf{v}}^{\varepsilon_{0}}
            +\varepsilon_{0}\left(\overline{\mathbf{v}}^{\varepsilon_{0}}\cdot\nabla\right)
            \overline{\mathbf{v}}^{\varepsilon_{0}}
            +\nabla \overline{\pi}^{\varepsilon_{0}}=\mu\Delta\overline{\mathbf{v}}^{\varepsilon_{0}}+
            \overline{\Theta}^{\varepsilon_{0}} \mathbf{g},
            \\
            \partial_{t}\overline{\Theta}^{\varepsilon_{0}}
            +\varepsilon_{0}\left(\overline{\mathbf{v}}^{\varepsilon_{0}}\cdot\nabla\right)
            \overline{\Theta}^{\varepsilon_{0}}
            +\overline{v_{2}}^{\varepsilon_{0}}\frac{d\overline{\theta}}{dx_{2}}=0,
            \\
            \nabla\cdot\overline{\mathbf{v}}^{\varepsilon_{0}}=0,
        \end{cases}
    \end{align}
with the initial value
\begin{align*}
    \left(\overline{\mathbf{v}}^{\varepsilon_{0}},\overline{\Theta}^{\varepsilon_{0}}\right)
    =(\mathbf{v},\Theta)(0).
\end{align*}
Moreover, by the inequality \eqref{2024032011}, the following estimate holds
\begin{align*}
    \begin{aligned}
    \sup\limits_{t\in[0,t_{k}]}\bigg{[}\big{\|}\overline{\mathbf{v}}^{\varepsilon_{0}}(t)\big{\|}_{H^2(\Omega_{\infty})}^2
    &+\big{\|}\overline{\Theta}^{\varepsilon_{0}}(t)\big{\|}_{H^1(\Omega_{\infty})}^2+\big{\|}\nabla \overline{\pi}^{\varepsilon_{0}}\big{\|}_{L^2(\Omega_{\infty})}^2
    +\big{\|}\overline{\mathbf{v}}_{t}^{\varepsilon_{0}}\big{\|}_{L^2(\Omega_{\infty})}^2
    +\big{\|}\overline{\Theta}_{t}^{\varepsilon_{0}}\big{\|}_{L^2(\Omega_{\infty})}^2\bigg{]}
   \\
   &+\int_{0}^{t}\left[\big{\|}\overline{\mathbf{v}}_{s}^{\varepsilon_{0}}\big{\|}_{H^1(\Omega_{\infty})}^2
    +\big{\|}\nabla\overline{\mathbf{v}}^{\varepsilon_{0}}\big{\|}_{H^1(\Omega_{\infty})}^2\right] d\tau
    \leq C(t_{k})\delta_{0}^2.
    \end{aligned}
\end{align*}
Thus, by virtue of the Banach-Alaoglu theorem, there exists a subsequence (not relabeled)
of $\left\{\left(\overline{\mathbf{v}}^{\varepsilon_{0}},\overline{\pi}^{\varepsilon_{0}},
\overline{\Theta}^{\varepsilon_{0}}\right)\right\}$ such that for $\varepsilon_{0}\rightarrow 0$,
\begin{align*}
    \begin{aligned}
        &\left(\overline{\mathbf{v}}^{\varepsilon_{0}}_{t},
        \nabla\overline{\pi}^{\varepsilon_{0}},\overline{\Theta}_{t}^{\varepsilon_{0}}\right)
        \rightarrow \left(\overline{\mathbf{v}}_{t},
        \nabla\overline{\pi},\overline{\Theta}_{t}\right)~~\text{weakly~star~in~}
        L^{\infty}\left(0,t_{k};\left(L^2(\Omega_{\infty})\right)^5\right),
        \\
        &\left(\overline{\mathbf{v}}^{\varepsilon_{0}},
        \overline{\Theta}^{\varepsilon_{0}}\right)
        \rightarrow \left(\overline{\mathbf{v}},
        \overline{\Theta}\right)~~\text{weakly~star~in~}
        L^{\infty}\left(0,t_{k};\left(H^2(\Omega_{\infty})\right)^2\times H^1(\Omega_{\infty})\right),
        \\
        &\left(\overline{\mathbf{v}}^{\varepsilon_{0}},
        \overline{\Theta}^{\varepsilon_{0}}\right)
        \rightarrow \left(\overline{\mathbf{v}},
        \overline{\Theta}\right)~~\text{strongly~in~}
        C\left([0,t_{k}],\left(L_{\text{loc}}^2(\Omega_{\infty})\right)^3\right),
    \end{aligned}
\end{align*}
and
\begin{align}\label{enenenen}
    \sup\limits_{0<t\leq t_{K}}\big{\|}\overline{v}_{2}(t)\big{\|}_{L^2(\Omega_{\infty})}\leq K\delta_{0},
    ~~\left(\overline{\mathbf{v}},
    \overline{\Theta}\right)\in C\left([0,t_{K}],\left(L^2(\Omega_{\infty})\right)^3\right).
\end{align}

If one takes the limit  $\varepsilon_{0}\rightarrow 0$ in \eqref{tianzhidao}, then we have
\begin{align*}
    \begin{cases}
        \partial_{t}\overline{\mathbf{v}}
        +\nabla \overline{\pi}=\mu\Delta\overline{\mathbf{v}}+
        \overline{\Theta}\mathbf{g},
        \\
        \partial_{t}\overline{\Theta}
        +\overline{v}_{2}\frac{d\overline{\theta}}{dx_{2}}=0,
        \\
        \nabla\cdot\overline{\mathbf{v}}=0.
    \end{cases}
\end{align*}
Therefore, $\left(\overline{\mathbf{v}},\overline{\Theta}\right)$ is just a strong solution to linearized
problem \eqref{naiverboundarycondition123}-\eqref{xianxingpart}. Then, according to Theorem \ref{weiyi}, we discover
$$
\overline{\mathbf{v}}=\mathbf{v}~~\text{on~}[0,t_{k}]\times\Omega.
$$
Thus, from \eqref{duideqiziji} and \eqref{enenenen}, we deduce
\begin{align*}
    2 K\delta_{0}\leq \big{\|}v_{2}(t_{K})\big{\|}_{L^2(\Omega_{\infty})}=\big{\|}\overline{v}_{2}(t_{K})\big{\|}_{L^2(\Omega_{\infty})}\leq K\delta_{0},
\end{align*}
which is a contradiction. Thus, we end the proof of this lemma.
\end{proof}
\section{Appendix}\label{peijian0104}
In this section, we give the proofs of some lemmas listed in the previous sections.

\subsection{Appendix A}\label{11226}

\begin{proof}\textbf{Proof of Lemma \ref{jidingli}.}

    Due to density, one can assume $u$ is smooth. We just prove the conclusion
    of $u(x_{1},0)$.
    Let $\Gamma^{dR}=\left\{x=(x_{1},x_{2})\big{|}x_{1}\in[-R,R],x_{2}\in[0,d]\right\},$ where $d\in (0,1]$.
    One can choose $\phi(x_{2})\in C^{\infty}\left(\Gamma^{dR}\right)$ such that (1)
    $\phi(x_{2})\equiv 1$ when $x_{2}\in \left[0,\frac{d}{2}\right]$; (2)$\phi(x_{2})\equiv 0$ when
    $x_{2}>\frac{3d}{4}$; (3) $0\leq\phi(x_{2})\leq 1$ when $x\in\Gamma^{dR}$.
Thus, a direct calculations gives
    \begin{align*}
        \begin{aligned}
        \int_{-R}^{R}|u(x_{1},0)|^{q}&dx_{1}
        =\int_{-R}^{R}\big{|}\int_{0}^{d}(\phi(x_{2})u(x_{1},x_{2}))_{{2}} dx_{2}\big{|}^q dx_{1}
        \\
        &=\int_{-R}^{R}\big{|}\int_{0}^{d}\phi'(x_{2})u(x_{1},x_{2})+
        \phi(x_{2})\partial_{{2}}u(x_{1},x_{2}) dx_{2}\big{|}^q dx_{1}
        \\
        &\leq 2^q\left(\int_{-R}^{R}\big{|}\int_{0}^{d}\phi'(x_{2})u(x_{1},x_{2}) dx_{2}\big{|}^q dx_{1}+
        \int_{-R}^{R}\big{|}\int_{0}^{d}
        \phi(x_{2})\partial_{{2}}u(x_{1},x_{2}) dx_{2}\big{|}^q dx_{1}\right).
        \end{aligned}
    \end{align*}
    By H\"{o}lder inequality, one can obtain
    \begin{align*}
        \begin{aligned}
        \int_{-R}^{R}\big{|}u(x_{1},0)|^{q}&dx_{1}\leq
        2^{q}\left(\int_{-R}^{R}\big{|}\int_{0}^{1}|\phi(x_{2})|^{\frac{q}{q-1}}dx_{2}\big{|}^{q-1}
\int_{0}^{1}\big{|}\partial_{2}u(x_{1},x_{2})\big{|}^{q} dx_{2} dx_{1}\right)
        \\
        &~~+2^q\left(\int_{-R}^{R}\big{|}\int_{0}^{1}\big{|}\phi'(x_{2})\big{|}^{\frac{q}{q-1}}dx_{2}\big{|}^{q-1}
\int_{0}^{1}\big{|}u(x_{1},x_{2})\big{|}^{q} dx_{2} dx_{1}\right).
        \end{aligned}
        \end{align*}
        Hence, we deduce that 
    \begin{align*}
        \int_{-R}^{R}\big{|}u(x_{1},0)\big{|}^{q}dx_{1}\leq C(q)\big{\|}u\big{\|}_{W^{1,q}(\Omega_{R})}^{q}.
    \end{align*}
    The proof of this lemma is completed.
\end{proof}

\begin{proof}\textbf{Proof of Lemma \ref{sobolveqianrudingli12241}.}

    Assume that $\mathbf{u}=(u_{1},u_{2})$ is smooth.
We first prove the conclusion for $u_{2}$. Since $u_{2}(x_{1},0)=u_{2}(x_{1},x_{2})|_{\Gamma_{3R}\cup
    \Gamma_{4R}}=0$, we have
    \begin{align*}
        u_{2}(x_{1},x_{2})
        = \int_{0}^{x_{2}}\partial_{{2}}u_{2}(x_{1},s_{2})ds_{2},~
            x_{1}\in [-R,R].
    \end{align*}
Thus, one obtains
\begin{align*}
    |u_{2}(x_{1},x_{2})|\leq \int_{0}^{1}\left|\partial_{2}u_{2}(x_{1},x_{2})\right|dx_{2},
\end{align*}
which implies by H\"{o}lder inequality that
\begin{align*}
    \big{\|}u_{2}\big{\|}_{L^2(\Omega_{R})}\leq \big{\|}\nabla u_{2}\big{\|}_{L^2(\Omega_{R})}.
\end{align*}

Next, we show by induction that when $q=2k~(k\in\mathbf{N}^{*})$,
$\|u_{2}\|_{L^q(\Omega_{R})}\leq C\|\nabla u_{2}\|_{L^2(\Omega_{R})}$.

We have proven that the conclusion is valid for $k=1$. Now, we assume
$\|u_{2}\|_{L^{2k}(\Omega_{R})}\leq C\|\nabla u_{2}\|_{L^{2}(\Omega_{R})}$.
Since the following equalities hold,
\begin{align*}
    \begin{aligned}
    &\left(u_{2}(x_{1},x_{2})\right)^{{k+1}}=(k+1)
        \int_{0}^{x_{2}}\left(u_{2}(x_{1},s_{2})\right)^{k}\partial_{2}u_{2}(x_{1},s_{2})ds_{2},~
        x_{1}\in [-R,R],
    \\
    &\left(u_{2}(x_{1},x_{2})\right)^{{k+1}}=(k+1)
        \int_{-R}^{x_{1}}\left(u_{2}(s_{1},x_{2})\right)^{k}\partial_{1}u_{2}(s_{1},x_{2})ds_{1}, ~
        x_{2}\in [0,1],
\end{aligned}
\end{align*}
then we have
\begin{align*}
    \begin{aligned}
    &(u_{2}(x_{1},x_{2}))^{{k+1}}\leq(k+1)
    \int_{0}^{1}\left|(u_{2}(x_{1},x_{2}))^{k}\partial_{2}u_{2}(x_{1},x_{2})\right|dx_{2},
    \\
    &(u_{2}(x_{1},x_{2}))^{{k+1}}\leq(k+1)
    \int_{-R}^{R}\left|(u_{2}(x_{1},x_{2}))^{k}\partial_{1}u_{2}(x_{1},x_{2})\right|dx_{1}.
    \end{aligned}
\end{align*}
By H\"{o}lder inequality and the conclusion of $q=2k$, it follows that
\begin{align*}
    \begin{aligned}
        \big{\|}u_{2}\big{\|}_{L^{2k+2}(\Omega_{R})}^{2k+2}
    \leq C\big{\|}\nabla u_{2}\big{\|}_{L^2(\Omega_{R})}^2
    \big{\|}u_{2}\big{\|}_{L^{2k}(\Omega_{R})}^{2k}
    \leq C\big{\|}\nabla u_{2}\big{\|}_{L^2(\Omega_{R})}^{2(k+1)},
    \end{aligned}
\end{align*}
which implies that for $q=2k+2~ (k\in\mathbf{N}^{*})$,
\begin{align*}
    \big{\|}u_{2}\big{\|}_{L^{2k+2}(\Omega_{R})}\leq C\big{\|}\nabla u_{2}\big{\|}_{L^2(\Omega_{R})}.
\end{align*}
Then, by the interpolation inequalities, one can conclude that $\forall q\in [2,+\infty)$,
\begin{align*}
    \big{\|}u_{2}\big{\|}_{L^{q}(\Omega_{R})}\leq C(q)\big{\|}\nabla u_{2}\big{\|}_{L^2(\Omega_{R})}.
\end{align*}

Let us turn to prove the conclusion of $u_{1}$. Since $\partial_{1}\int_{0}^{1}u_{1}(x_{1},x_{2})dx_{2}=-\int_{0}^{1}\partial_{2}u_{2}(x_{1},x_{2})dx_{2}=0$ for
$x_{1}\in [-R,R]$, then $\int_{0}^{1}u_{1}(x_{1},x_{2})dx_{2}\equiv constant$.
We will show that
\begin{align}\label{yaozheng0207}
\int_{0}^{1}
u_{1}(x_{1},x_{2})dx_{2}\equiv 0.
\end{align}
In fact, we have
\begin{align}\label{yongyong1225}
    0=\int_{\partial\Omega_{R}}x_{1}\mathbf{u}\cdot\mathbf{n}ds=\int_{\Omega_{R}}\nabla\cdot(x_{1}\mathbf{u}) dx=\int_{\Omega_{R}}u_{1}(x_{1},x_{2})dx_{1}dx_{2},
\end{align}
which implies \eqref{yaozheng0207}.

As a result, there exists a point $(x_{1},\xi)$ for each $x_{1}\in[-R,R]$ such that
$u_{1}(x_{1},\xi)=0$. Moreover, we know that $u_{1}|_{\Gamma_{3R}\cup\Gamma_{4R}}=0$. Thus,
a similar argument can tell us that for $q\in [2,+\infty)$,
\begin{align*}
    \big{\|}u_{1}\big{\|}_{L^q(\Omega_{R})}\leq C\big{\|}\nabla u_{1}\big{\|}_{L^2(\Omega_{R})}.
\end{align*}
Therefore, we end the proof of this lemma.
\end{proof}
\begin{remark}
It is not difficult to find that if for each $x_{1}$ and $x_{2}$, there always exist
$(x_{1},\xi_{2})$ and $(\xi_{1},x_{2})$ such that $u(x_{1},\xi_{2})=u(\xi_{1},x_{2})=0$. So
the New-Leibniz formula with zero boundary-value for vertical and horizontal direction can
be employed. Then, by
the mathematical induction and interpolation inequalities we can
verify conclusions. Hence, in the proof of next lemma, we just prove that in each of two directions
there exists a point where the function takes zero, i.e., the New-Leibniz formula with
zero boundary-value term can be utilized.
\end{remark}

\begin{proof}\textbf{Proof of Lemma \ref{jiushizheyang12242}}.

    We first prove the conclusion for $\partial_{2}u_{2}$. Since $u_{2}(-R,x_{2})=u_{2}(R,x_{2})=0$ for each
    $x_{2}\in[0,1]$, thus we have  $\partial_{2}u_{2}|_{\Gamma_{3R}\cup\Gamma_{4R}}=0$; Since
    $u_{2}(x_{1},0)=u_{2}(x_{1},1)=0$ for each  $x_{1}\in[-R,R]$, thus from the Lagrange's theorem there exists a point $(x_{1},\xi_{2})$ such
    that $\partial_{2}u_{2}(x_{1},\xi_{2})=0$. Then, one can verify the conclusions of $\partial_{2}u_{2}$.
    By the incompressible condition $\partial_{1}u_{1}+\partial_{2}u_{2}=0$, one can conclude that the conclusion of $\partial_{1}u_{1}$ is valid.

    Regarding $\partial_{1}u_{2}$,
    since $u_{2}(-R,x_{2})=u_{2}(R,x_{2})=0$, then
    there exists a point $(\xi_{1},x_{2})$ for each
    $x_{2}\in[0,1]$, such that $\partial_{1}u_{2}(\xi_{1},x_{2})$;
     Since $u_{2}(x_{1},0)=u_{2}(x_{1},1)=0$, then we know $\partial_{1}u_{2}(x_{1},0)=
    \partial_{1}u_{2}(x_{1},1)=0$. A similar argument as Lemma \ref{sobolveqianrudingli12241} gives the proof of $\partial_{1}u_{2}$.

     With respect to $\partial_{2}u_{1}$,  since $u_{1}(-R,x_{2})=u_{1}(R,x_{2})=0$, then
    $\partial_{2}u_{1}(-R,x_{2})=\partial_{2}u_{1}(R,x_{2})=0$;
    We can not find the point we need for all $x_{1}\in(-R,R)$ fixed. However, one still has
\begin{align}\label{womenxuyaode0126}
        \begin{aligned}
        &\begin{aligned}
            (\partial_{2}u_{1}(x_{1},x_{2}))^2-(\partial_{2}u_{1}(x_{1},0))^2
            &=2\int_{0}^{x_{2}}\partial_{22}^2 u_{1}(x_{1},s_{2})\partial_{2}u_{1}(x_{1},s_{2})ds_{2}
            \\
            &\leq 2\int_{0}^{1}|\partial_{22}^2 u_{1}(x_{1},x_{2})||\partial_{2}u_{1}(x_{1},x_{2})|dx_{2},
        \end{aligned}
        \\
        &\begin{aligned}
            (\partial_{2}u_{1}(x_{1},x_{2}))^2-(\partial_{2}u_{1}(-R,x_{2}))^2
            &=2\int_{-R}^{x_{1}}\partial_{12}^2 u_{1}(s_{1},x_{2})\partial_{2}u_{1}(s_{1},x_{2})ds_{1}
            \\
            &\leq 2\int_{-R}^{R}|\partial_{12}^2 u_{1}(x_{1},x_{2})||\partial_{2}u_{1}(x_{1},x_{2})|dx_{1}.
        \end{aligned}
    \end{aligned}
    \end{align}
    Since $\partial_{2}u_{1}(-R,x_{2})=0$ and $\partial_{2}u_{1}(x_{1},0)=-\frac{k_{0}}{\mu}u_{1}(x_{1},0)$,
    then one has from \eqref{womenxuyaode0126} that
    \begin{align}\label{buyaobeishang0126}
        \begin{aligned}
            (\partial_{2}u_{1}(x_{1},x_{2}))^{4}&\leq 4\int_{0}^{1}|\partial_{22}^{2}u_{1}||\partial_{2}u_{1}|dx_{2}
            \int_{-R}^{R}|\partial_{12}^{2}u_{1}||\partial_{2}u_{1}|dx_{1}
            \\
            &~~+\frac{2k_{0}^{2}}{\mu^2}|u_{1}(x_{1},0)|^2\int_{-R}^{R}|
            \partial_{12}^{2}u_{1}||\partial_{2}u_{1}|dx_{1}:=O_{1}+O_{2}.
        \end{aligned}
    \end{align}
    By H\"{o}lder inequality, Lemma \ref{jidingli} and Lemma \ref{sobolveqianrudingli12241}, one can get that
    \begin{align*}
        \begin{aligned}
            \int_{\Omega_{R}}O_{1}dx\leq \big{\|}\nabla^2 u_{1}\big{\|}_{L^2(\Omega_{R})}^2\big{\|}\nabla u_{1}\big{\|}_{L^2(\Omega_{R})}^2,~
            \int_{\Omega_{R}}O_{2}dx\leq C\big{\|}\nabla u_{1}\big{\|}_{L^2(\Omega_{R})}^3\big{\|}\nabla^2 u_{1}\big{\|}_{L^2(\Omega_{R})},
        \end{aligned}
    \end{align*}
    where $C$ is independent of $R$,
which together with \eqref{buyaobeishang0126} yields that
    \begin{align*}
\big{\|}\partial_{2}u_{1}\big{\|}_{L^4(\Omega_{R})}\leq C\left(\big{\|}\nabla^2 u_{1}\big{\|}_{L^2(\Omega_{R})}^{\frac{1}{2}}
       \big{\|}\nabla u_{1}\big{\|}_{L^2(\Omega_{R})}^{\frac{1}{2}}+\big{\|}\nabla u_{1}\big{\|}_{L^2(\Omega_{R})}^{\frac{3}{4}}
        \big{\|}\nabla^2 u_{1}\big{\|}_{L^2(\Omega_{R})}^{\frac{1}{4}}\right).
    \end{align*}
    We end the proof of this lemma.
\end{proof}

\begin{proof}\textbf{Proof of Lemma \ref{jiushizheyang12243}}.
    Assume $u\in W^{1,q}(\Omega_{\infty})\cap C^{\infty}(\Omega_{\infty})$. First, we  extend $u$ as follows
\begin{align*}
    \widetilde{u}(x)=
    \begin{cases}
        -3 u(x_{1},-x_{2})+4 u\left(x_{1},-\frac{x_{2}}{2}\right),~~~~~~x_{1}\in(-\infty,+\infty),~
        x_{2}\in(-1,0],
        \\
        u(x_{1},x_{2}),~~~~~~~~~~~~~~~~~~~~~~~~~~~~~~~~~~~x_{1}\in(-\infty,+\infty),~
        x_{2}\in[0,1],
        \\
        -3u(x_{1},2-x_{2})+4u\left(x_{1},\frac{3-x_{2}}{2}\right),
        ~~~~x_{1}\in(-\infty,+\infty),~
        x_{2}\in [1,2).
    \end{cases}
\end{align*}
We shall show  $\widetilde{u}(x)\in C^{1}(\Omega_{\infty,3})$, where $
\Omega_{\infty,3}=\left\{x=(x_{1},x_{2})\big{|}x_{1}\in(-\infty,+\infty)~\text{and}~x_{2}\in (-1,2)\right\}$. In fact, it is easy to see
\begin{align*}
    u(x_{1},0)=-3u(x_{1},0)+4u(x_{1},0),~u(x_{1},1)=-3u(x_{1},2-1)+4u\left(x_{1},\frac{3-1}{2}\right),
\end{align*}
which implies that $
\widetilde{u}$ is continuous. And a
series of calculation gives that
\begin{align*}
    u_{x_{2}}(x_{1},0)=\left[-3u\big{(}x_{1},-x_{2}\big{)}+4u\left(x_{1},-\frac{x_{2}}{2}\right)\right]_{x_{2}}\bigg{|}_{x_{2}=0},
    ~u_{x_{2}}(x_{1},1)=\left[-3u\big{(}x_{1},-x_{2}\big{)}+4u\left(x_{1},-\frac{x_{2}}{2}\right)\right]_{x_{2}}\bigg{|}_{x_{2}=1},
\end{align*}
indicating that $\partial_{2}\widetilde{u}$ is
continuous. Thus, $\widetilde{u}\in C^{1}(\Omega_{\infty,3})$.

Moreover, note that there exists a constant $C>0$ such that
\begin{align}\label{wuzhideshihou1225}
    \big{\|}\widetilde{u}\big{\|}_{W^{1,q}(\Omega_{\infty,3})}\leq C\big{\|}u\big{\|}_{W^{1,q}(\Omega_{\infty})}.
\end{align}
Let $x\in\Omega_{\infty}$ and $B(x,r_{x})\subset\Omega_{\infty,3}$, where $B(x,r_{x})$ is a ball
whose center is $x$ and radius is $r_{x}$.
Next, we will show that there exists a constant $C>0$ such that
\begin{align}\label{haoxiang1225}
    \frac{1}{r_{x}^2}\int_{B(x,r_{x})}\left|\widetilde{u}(x)-\widetilde{u}(y)\right|dy\leq C\int_{B(x,r_{x})}
    \frac{|D\widetilde{u}(y)|}{|x-y|}dy.
\end{align}
Indeed,  $\forall z\in \partial B(0,1)$ and $0<s<r_{x}$, one has
\begin{align*}
    \begin{aligned}
   \left|\widetilde{u}(x+sz)-\widetilde{u}(x)\right|=\left|\int_{0}^{s}
    \frac{d}{dt}\widetilde{u}(x+tz)dt
    \right|\leq \int_{0}^{s}\left|D\widetilde{u}(x+tz)\right|dt.
    \end{aligned}
\end{align*}
Taking $y=x+tz$, then one gets
\begin{align*}
    \begin{aligned}
        \int_{\partial B(0,1)}|\widetilde{u}(x+sz)-\widetilde{u}(x)|dS(z)
        &\leq \int_{0}^{s}\int_{\partial B(0,1)}|D\widetilde{u}(x+tz)|dS(z)dt
        \\
        &=\int_{B(x,s)}\frac{|D\widetilde{u}(y)|}{|x-y|}dy
        \\
        &\leq \int_{B(x,r_{x})}
        \frac{|D\widetilde{u}(y)|}{|x-y|}dy.
    \end{aligned}
\end{align*}
Multiplying the above inequality by $s$ and integrating from $0$ to $r_{x}$ with respect to $s$, one has
\begin{align*}
    \int_{B(x,r_{x})}|\widetilde{u}(y)-\widetilde{u}(x)|dy \leq \frac{r_{x}^2}{2}
    \int_{B(x,r_{x})}\frac{|D\widetilde{u}(y)|}{|x-y|}dy.
\end{align*}
Thus, we deduce that \eqref{haoxiang1225} is valid.

Since $x\in\Omega_{\infty}$, then $\widetilde{u}(x)=u(x)$. When $r_{x}=1$ in \eqref{haoxiang1225}, we have
\begin{align*}
    \begin{aligned}
    |u(x)|=\frac{1}{\pi}\int_{B(x,1)}|\widetilde{u}(x)|dy
    &\leq \frac{1}{\pi}\int_{B(x,1)}|\widetilde{u}(y)-\widetilde{u}(x)|dy +
    \frac{1}{\pi}\int_{B(x,1)}|\widetilde{u}(y)|dy
    \\
    &\leq C\int_{B(x,1)}\frac{|D\widetilde{u}(y)|}{|x-y|}dy +
    \frac{1}{\pi}\int_{B(x,1)}|\widetilde{u}(y)|dy.
    \end{aligned}
\end{align*}
By H\"{o}lder inequality, it follows that
\begin{align*}
    |u(x)|\leq C\big{\|}D\widetilde{u}\big{\|}_{L^q(\Omega_{\infty,3})}
    \left[\int_{B(x,1)}\frac{1}{|x-y|^{\frac{q}{q-1}}}dy\right]^{1-\frac{1}{q}}
    +C\big{\|}\widetilde{u}\big{\|}_{L^{q}(\Omega_{\infty,3})}.
\end{align*}
In addition, the condition $q>2$ implies $\frac{q}{q-1}<2$, then we discover
$\int_{B(x,1)}\frac{1}{|x-y|^{\frac{q}{q-1}}}dy<+\infty$. And from \eqref{wuzhideshihou1225},
one obtains
\begin{align*}
    |u(x)|\leq C\big{\|}\widetilde{u}\big{\|}_{W^{1,q}(\Omega_{\infty,3})}\leq
    C\big{\|}u\big{\|}_{W^{1,q}(\Omega_{\infty})},
\end{align*}
which means
\begin{align*}
    \sup\limits_{x\in\Omega_{\infty}}|u(x)|\leq C\big{\|}u\big{\|}_{W^{1,q}(\Omega_{\infty})}.
\end{align*}
This completes the proof.
\end{proof}

\begin{proof}\textbf{Proof of Corollary \ref{ohmygod1228}}.
The same argument as the proof in the Lemma \eqref{jiushizheyang12243}, we extend $u$ as follows,
    \begin{align*}
        \widetilde{u}(x)=
        \begin{cases}
            -3 u(x_{1},-x_{2})+4 u\left(x_{1},-\frac{x_{2}}{2}\right),~~~~~~x_{1}\in [-R,R],~
            x_{2}\in[-1,0],
            \\
            u(x_{1},x_{2}),~~~~~~~~~~~~~~~~~~~~~~~~~~~~~~~~~~~x_{1}\in [-R,R],~
            x_{2}\in[0,1],
            \\
            -3u(x_{1},2-x_{2})+4u\left(x_{1},\frac{3-x_{2}}{2}\right),
            ~~~~x_{1}\in[-R,R],~
            x_{2}\in [1,2].
        \end{cases}
    \end{align*}
    And we further need to extend $\widetilde{u}(x)$ as follows
    \begin{align*}
    \widehat{u}(x)=
    \begin{cases}
    -3\widetilde{u}(-2R-x_{1},x_{2})+4\widetilde{u}\left(-\frac{x_{1}}{2}-\frac{3R}{2},x_{2}\right),~~~~x_{1}\in[-3R,-R],~x_{2}\in[-1,2],
    \\
    \widetilde{u}(x_{1},x_{2}),~~~~~~~~~~~~~~~~~~~~~~~~~~~~~~~~~~~~~~~~~~~~~x_{1}\in[-R,R],~x_{2}\in[-1,2],
    \\
    -3\widetilde{u}(2R-x_{1},x_{2})+4\widetilde{u}\left(\frac{3R}{2}-\frac{x_{1}}{2},x_{2}\right),~~x_{1}\in[R,3R],~x_{2}\in[-1,2].
    \end{cases}
    \end{align*}
    Obviously, one can verify that  $\big{\|}\widehat{u}\big{\|}_{L^{\infty}(\Omega_{R})}=\big{\|}u\big{\|}_{L^{\infty}(\Omega_{R})}$,
    and $\big{\|}\widehat{u}\big{\|}_{W^{1,q}([-3R,3R]\times[-1,2])}\leq C\big{\|}u\big{\|}_{W^{1,q}(\Omega_{R})}$.
    Besides, $\forall$ $x\in\Omega_{R}$, it is easy to see that $B(x,1)\subset[-3R,3R]\times[-1,2]$.
    Then the rest proof can be completed by the analogous analysis of Lemma \ref{jiushizheyang12243}.
\end{proof}

\subsection{Appendix B}\label{appendix11226}
We give the definition of $q-$generalized solution of
the following equation:
\begin{align}\label{tuiguanjie}
    \begin{cases}
        -\mu\Delta \mathbf{u}+\nabla p=\mathbf{f}(x),~~x\in\Omega_{R},
        \\
        \mathbf{u}\cdot\mathbf{n}|_{\Gamma_{1R}\cup\Gamma_{2R}}=0,
        ~ \mathbf{u}|_{\Gamma_{3 R}\cup\Gamma_{4R}}=\mathbf{0},
        \\
        \partial_{2} u_{1}(x_{1},1)=\frac{k_{1}}{\mu}u_{1}(x_{1},1),~
        \partial_{2} u_{1}(x_{1},0)=-\frac{k_{0}}{\mu}u_{1}(x_{1},0),
        \\
        \nabla\cdot \mathbf{u}=0,
    \end{cases}
\end{align}
where $\mathbf{f}\in \left[L^q(\Omega_{R})\right]^2$.
\begin{definition}\label{tuiguanjiedingyi}
    $(\mathbf{u},p)$ is called the $q-$generalized solution of problem \eqref{tuiguanjie}, if
    \begin{enumerate}[$(1)$]
        \item  $\mathbf{u}\in \mathbf{G}^{1,q}(\Omega_{R})$ and $p\in L^{q}(\Omega_{R})$;
        \item $\forall \mathbf{v}\in \mathbf{E}^{1,q'}(\Omega_{R})$, where $\frac{1}{q}+\frac{1}{q'}=1$, one has
        \begin{align}\label{dengshichengli}
            \begin{aligned}
            &\mu\int_{\Omega_{R}}\nabla\mathbf{u}\cdot\nabla\mathbf{v}dx
        -\int_{-R}^{R}\left[k_{1}u_{1}(x_{1},1)v_{1}(x_{1},1)+k_{0}u_{1}(x_{1},0)v_{1}(x_{1},0)\right]dx_{1}
        \\
        &~~~~~~~~~~~~~~-\int_{\Omega_{R}}p\nabla\cdot\mathbf{v}dx
        =\int_{\Omega_{R}}\mathbf{f}\cdot\mathbf{v}dx.
            \end{aligned}
        \end{align}
    \end{enumerate}
    In particular, when $q=2$, $(\mathbf{u},p)$ is called the weak solution.
\end{definition}

Next, we give an equivalent definition of $q-$ generalized as follows.
\begin{lemma}\label{dengjia}
    Let $\mathbf{u}\in \mathbf{G}^{1,q}(\Omega_{R})$. And $\forall~ \mathbf{v}\in \mathbf{G}^{1,q'}(\Omega_{R})$, one has
    \begin{align}\label{dengshichenglidengjia}
        \begin{aligned}
\mu\int_{\Omega_{R}}\nabla\mathbf{u}\cdot\nabla\mathbf{v}dx
    -\int_{-R}^{R}\left[k_{1}u_{1}(x_{1},1)v_{1}(x_{1},1)
    +k_{0}u_{1}(x_{1},0)v_{1}(x_{1},0)\right]dx_{1}
    =\int_{\Omega_{R}}\mathbf{f}\cdot\mathbf{v}dx.
        \end{aligned}
    \end{align}
    Then, we can conclude that there exists a unique $p\in \widetilde{L}^{q}(\Omega_{R})=\left\{p\in L^{q}(\Omega_{R})\big|\int_{\Omega_{R}}p
    dx=0\right\} $ such that $(\mathbf{u},p)$ is the $q-$generalized solution of
    problem \eqref{tuiguanjie} and satisfies the inequality
    \begin{align}\label{tidup0117}
        \big{\|}p\big{\|}_{L^{q}}\leq C(R)
        \left(\big{\|}\mathbf{f}\big{\|}_{L^q(\Omega_{R})}+\big{\|}\mathbf{u}\big{\|}_{W^{1,q}(\Omega_{R})}\right).
    \end{align}
\end{lemma}
\begin{remark} Note that $C$ depends on $R$ in \eqref{tidup0117}. However, one can still prove the existence of strong solution by
    domain expanding technique since only $\nabla p$ is needed.
\end{remark}
Before giving the proof of Lemma \ref{dengjia}, we first give the following proposition.
\begin{proposition}[Riesz represent theorem on
    $\widetilde{L}^{q}(\Omega_{R})$]\label{Rieszbiaoshiyinli0110}
Let $1<q<+\infty$, and$\mathcal{L}\in\left[\widetilde{L}^{q}(\Omega_{R})\right]'$,
i.e., $\mathcal{L}$ is a linear and bounded functional on $\widetilde{L}^{q}(\Omega_{R})$.
Then there exists a unique $v\in \widetilde{L}^{q'}(\Omega_{R})$ $\left(q'=\frac{q}{q-1}\right)$,
such that $\forall$ $u\in \widetilde{L}^{q}(\Omega_{R})$,
\begin{align*}
\mathcal{L}(u)=\int_{\Omega_{R}}u(x)v(x)dx,
\end{align*}
and $\big{\|}v\big{\|}_{q'}\leq 2\big{\|}\mathcal{L}\big{\|}$.
\end{proposition}
\begin{remark}\label{yongyongpinglun0127} From above proposition, one can see
    \begin{enumerate}[$1$]
    \item Note that $\left[\widetilde{L}^{q}(\Omega_{R})\right]'$ is not isometric
    isomorphism to $\widetilde{L}^{q'}(\Omega_{R})$;
    \item With respect to the proof of Proposition \ref{Rieszbiaoshiyinli0110}, we follow the proof of
    the theorem saying that $\left[L^{q}(\Omega_{R})\right]'$ is isometric isomorphism to $L^{q'}(\Omega_{R})$
    in \cite{Adams2003}. Thus, we list the several key steps: (1) one can easily prove
    $\widetilde{L}^{q}(\Omega_{R})$ is a Banach space under the norm of $\|\cdot\|_{L^{q}(\Omega_{R})}$;
    (2) by Clarkson's inequalities (see P37 in \cite{Adams2003}),
    one can prove that $\widetilde{L}^{q}(\Omega_{R})$ is uniform convex;
    (3) by the similar argument of Lemma 2.32 in \cite{Adams2003}, one can prove the following proposition:
\begin{proposition}\label{xuyaoyongdade0110}
        Let $1<q<\infty$. If $\mathcal{L}\in \left[\widetilde{L}^{q}(\Omega_{R})\right]'$
        and $\big{\|}\mathcal{L}\big{\|}=1$, then there exists a unique $w\in \widetilde{L}^{q}(\Omega_{R})$
        such that $\|w\|_{q}=\mathcal{L}(w)=1$. Conversely, if $w\in \widetilde{L}^{q}(\Omega_{R})$
        is given and $\big{\|}w\big{\|}_{q}=1$, then there exists a unique $\mathcal{L}\in \left[\widetilde{L}^{q}(\Omega_{R})\right]'$
        such that $\big{\|}\mathcal{L}\big{\|}=\mathcal{L}(w)=1$.
        \end{proposition}
    \end{enumerate}
\end{remark}
\begin{proof}(\textbf{Proof of Proposition \ref{Rieszbiaoshiyinli0110}.})

    If $\mathcal{L}=0$, then we have $v=0$. If $\mathcal{L}\neq 0$, one can assume $\big{\|}\mathcal{L}\big{\|}=1$.
    Then from Proposition \ref{xuyaoyongdade0110},
    there exists a unique $w\in\widetilde{L}^{q}(\Omega_{R})$ such that
    $\big{\|}w\big{\|}_{q}=1$ and $\mathcal{L}(w)=1$. Let
    \begin{align*}
    v(x)=
    \begin{cases}
    |w(x)|^{q-2}w(x)-\frac{\int_{\Omega_{R}}|w(x)|^{q-2}w(x)dx}{|\Omega_{R}|},~~w(x)\neq 0,
    \\
    0,~~~~~~~~~~~~~~~~~~~~~~~~~~~~~~~~~~~~~~~~~~~~~~~~w(x)=0.
    \end{cases}
    \end{align*}
    Then $v(x)\in \widetilde{L}^{q'}(\Omega_{R})\subset L^{q'}(\Omega_{R})$. Furthermore, we define
    $\mathcal{L}_{v}:\widetilde{L}^{q}(\Omega_{R})\rightarrow \mathbf{R}$ as follows:
    \begin{align*}
    \mathcal{L}_{v}(u)=\int_{\Omega_{R}}u(x)v(x)dx.
    \end{align*}
    Obviously, a direct calculations gives $\big{\|}\mathcal{L}_{v}\big{\|}=|\mathcal{L}_{v}(w)|=1$, and by Proposition \ref{xuyaoyongdade0110} again, we have
$\mathcal{L}=\mathcal{L}_{v}$.

    Since $\big{\|}v\big{\|}_{q'}\leq 2$, then one can conclude that $\big{\|}v\big{\|}_{q'}\leq 2\big{\|}\mathcal{L}\big{\|}$.
    \end{proof}
    \begin{proof}(\textbf{Proof of Lemma \ref{dengjia}})
        Obviously, when $\mathbf{u}$ is $q-$generalized solution, then $\mathbf{u}$ satisfies
        \eqref{dengshichenglidengjia}.

        Conversely, suppose that $\mathbf{u}$ satisfies \eqref{dengshichenglidengjia}, then one can define a functional
     $\mathcal{F}:\mathbf{E}^{1,q'}(\Omega_{R})\rightarrow \mathbf{R}$
        as follows,
    \begin{align}\label{woxiwangnibeiaizhe1226}
          \begin{aligned}
          \mathcal{F}(\mathbf{v})&=\mu\int_{\Omega_{R}}\nabla\mathbf{u}\cdot\nabla\mathbf{v}dx
          -\int_{-R}^{R}[k_{1}u_{1}(x_{1},1)v_{1}(x_{1},1)
          \\
          &~~+k_{0}u_{1}(x_{1},0)v_{1}(x_{1},0)]dx_{1}
          -\int_{\Omega_{R}}\mathbf{f}\cdot\mathbf{v}dx,~\forall~\mathbf{v}
        \in \mathbf{E}^{1,q'}(\Omega_{R}).
          \end{aligned}
        \end{align}
        Then, it is easy to see that $\mathcal{F}$ is a linear and bounded on $\mathbf{E}^{1,q'}(\Omega_{R})$,
         and $\mathcal{F}\equiv 0$ on $\mathbf{G}^{1,q'}(\Omega_{R})$.

        We also define an operator $\mathbf{A}:\mathbf{E}^{1,q'}(\Omega_{R})\rightarrow
        \widetilde{L}^{q'}(\Omega_{R})$ in form
    \begin{align*}
      \mathbf{A}\mathbf{v}=\nabla\cdot\mathbf{v}.
    \end{align*}
    Note that $\mathbf{A}$ is linear and bounded. Moreover, from the excise III.3.5 in \cite{Galdi2011},
    We find that $\mathbf{A}$ is surjective. Thus, by Theorem 2.19 in \cite{brezis_functional_2011}, one has that
    $R(\mathbf{A}^{*})$ is closed,
    where $\mathbf{A}^{*}$ is the conjugate operator of $\mathbf{A}$. Thus, it is known that
    \begin{align*}
      [Ker(\mathbf{A})]^{\bot}=R(\mathbf{A}^{*}).
    \end{align*}
    Since $Ker(\mathbf{A})=\mathbf{G}^{1,q'}(\Omega_{R})$, thus,
    $\mathcal{F}\in R(\mathbf{A}^{*})$. Then there exists a linear and bounded functional $\mathcal{L}$
    on $\widetilde{L}^{q'}(\Omega_{R})$, such that
    $\mathbf{A}^{*}\mathcal{L}=\mathcal{F}$. Then by the Proposition \ref{Rieszbiaoshiyinli0110},
    there exists a unique
    $p\in \widetilde{L}^{q}(\Omega_{R})$ such that
    $\mathcal{L}(A\mathbf{v})=\int_{\Omega_{R}} p\nabla\cdot\mathbf{v}dx$,
    $\forall \mathbf{v}\in \mathbf{E}^{1,q'}(\Omega_{R})$.

    Since $\mathbf{A}^{*}$ may not be injective, we can assume that there exist
    $\mathcal{L}_{1}$ and $\mathcal{L}_{2}$ such that,
    $\mathbf{A}^{*}\mathcal{L}_{1}=\mathbf{A}^{*}\mathcal{L}_{2}=\mathcal{F}$, then there exist $p_{1}$
    and $p_{2}$ such that
    \begin{align*}
    \int_{\Omega_{R}} p_{1}\nabla\cdot\mathbf{v}dx=\int_{\Omega_{R}} p_{2}\nabla\cdot\mathbf{v}dx,
    \end{align*}
    that is, $\int_{\Omega_{R}} (p_{1}-p_{2})\nabla\cdot\mathbf{v}dx=0$, then $p_{1}-p_{2}=constant$. By
    $\int_{\Omega_{R}}p_{1}dx=\int_{\Omega_{R}}p_{2}dx=0$, one can conclude that $p_{1}=p_{2}$ almost everywhere, i.e., the uniqueness of
    $p$ is verified.

    By the definition of conjugate operator, one has
    \begin{align*}
      \mathcal{F}(\mathbf{v})=\mathbf{A}^{*}\mathcal{L}(\mathbf{v})=
      \mathcal{L}(\mathbf{A}\mathbf{v})=\int_{\Omega_{R}} p\nabla\cdot\mathbf{v}dx,
      ~~\forall~\mathbf{v}\in \mathbf{E}^{1,q'}(\Omega_{R}).
    \end{align*}
    Therefore, \eqref{dengshichengli} is verified.

    To derive \eqref{tidup0117},
    we need to investigate the following equation.
    \begin{align}\label{ningjingdelian0117}
    \begin{cases}
    \nabla\cdot \mathbf{v}=|p(x)|^{q-2} p(x)-\frac{\int_{\Omega_{R}}|p(x)|^{q-2} p(x)dx}{|\Omega_{R}|}\equiv g,
    \\
    \mathbf{v}\in \left[W^{1,q'}_{0}(\Omega_{R})\right]^2,
    \\
    \big{\|}\mathbf{v}\big{\|}_{W^{1,q'}_{0}(\Omega_{R})}\leq C(\Omega_{R})\big{\|}g\big{\|}_{L^{q'}(\Omega_{R})}.
    \end{cases}
    \end{align}
    The equation \eqref{ningjingdelian0117} has been studied by many researchers\cite{Bourgain2002,Russ2013,Panasenko2014} in which we find that
    there exists at least one
    solution $\mathbf{v}$ for equation \eqref{ningjingdelian0117}. And a direct computation gives that $\big{\|}g\big{\|}_{L^{q'}(\Omega_{R})}\leq 2\big{\|}p\big{\|}_{L^{q}(\Omega_{R})}^{q-1}$.
    Then, $\big{\|}\mathbf{v}\big{\|}_{W^{1,q'}_{0}(\Omega_{R})}\leq C(\Omega_{R})\big{\|}p\big{\|}_{L^{q}(\Omega_{R})}^{q-1}$.

    In \eqref{woxiwangnibeiaizhe1226}, let $\mathbf{v}$ be the solution of \eqref{ningjingdelian0117}. Then, one can
    obtain
    \begin{align*}
    \begin{aligned}
          \mu\int_{\Omega_{R}}\nabla\mathbf{u}\cdot\nabla\mathbf{v}dx
          -\int_{\Omega_{R}}\mathbf{f}\cdot\mathbf{v}dx=\int_{\Omega_{R}}p\nabla\cdot\mathbf{v}dx=\|p\|_{L^{q}(\Omega_{R})}^{q}.
          \end{aligned}
    \end{align*}
    Furthermore, by H\"{o}lder inequality, one can verify \eqref{tidup0117}.
    \end{proof}

Next, we will consider the Stokes estimate which is very useful for improving the regularity of weak solution.
\begin{proof}(\textbf{Proof of Lemma \ref{tishengzhengzexing}.})
    Since $\mathbf{u}$ is $q-$generalized solution, then for all $\mathbf{v}\in \widetilde{\mathbf{G}}^{1,q'}(\Omega_{R})
    =\left\{\mathbf{v}\in \left[W_{0}^{1,q'}(\Omega_{R})\right]^2\big{|}\nabla\cdot\mathbf{v}=0\right\}$, one has
        \begin{align*}
            \mu\int_{\Omega_{R}}\nabla\mathbf{u}\cdot\nabla\mathbf{v}dx=\int_{\Omega_{R}}\mathbf{f}\cdot
            \mathbf{v}dx.
        \end{align*}
We extend $\mathbf{u}$ in form
        \begin{align*}
    \widetilde{\mathbf{u}}(x)=
    \begin{cases}
        -3 \mathbf{u}(x_{1},-x_{2})+4 \mathbf{u}\left(x_{1},-\frac{x_{2}}{2}\right),~~~~~~x_{1}\in [-R,R],~
        x_{2}\in[-1,0],
        \\
        \mathbf{u}(x_{1},x_{2}),~~~~~~~~~~~~~~~~~~~~~~~~~~~~~~~~~~~x_{1}\in [-R,R],~
        x_{2}\in[0,1],
        \\
        -3\mathbf{u}(x_{1},2-x_{2})+4\mathbf{u}\left(x_{1},\frac{3-x_{2}}{2}\right),
        ~~~~x_{1}\in[-R,R],~
        x_{2}\in [1,2],
    \end{cases}
\end{align*}
and further extend $\widetilde{\mathbf{u}}(x)$ in following
\begin{align*}
\widehat{\mathbf{u}}(x)=
\begin{cases}
-3\widetilde{\mathbf{u}}(-2R-x_{1},x_{2})+4\widetilde{\mathbf{u}}\left(-\frac{x_{1}}{2}-\frac{3R}{2},x_{2}\right),~~~~x_{1}\in[-3R,-R],~x_{2}\in[-1,2],
\\
\widetilde{\mathbf{u}}(x_{1},x_{2}),~~~~~~~~~~~~~~~~~~~~~~~~~~~~~~~~~~~~~~~~~~~~~x_{1}\in[-R,R],~x_{2}\in[-1,2],
\\
-3\widetilde{\mathbf{u}}(2R-x_{1},x_{2})+4\widetilde{\mathbf{u}}\left(\frac{3R}{2}-\frac{x_{1}}{2},x_{2}\right),
~~x_{1}\in[R,3R],~x_{2}\in[-1,2].
\end{cases}
\end{align*}
Then, $\widehat{\mathbf{u}}\in W^{1,q}(T_{R})$, where $T_{R}=[-3R,3R]\times[-1,2]$, and
$\widehat{\mathbf{u}}|_{\Omega_{R}}\equiv\mathbf{u}$.
Thus, one still has
   \begin{align}\label{chentianyi0116}
            \mu\int_{\Omega_{R}}\nabla\widehat{\mathbf{u}}\cdot\nabla\mathbf{v}dx=\int_{\Omega_{R}}\mathbf{f}\cdot
            \mathbf{v}dx,~\forall~\mathbf{v}\in \widetilde{\mathbf{G}}^{1,q'}(\Omega_{R}).
   \end{align}
   For $h>0$ small enough and fixed, we define $\overline{\Omega}_{R,h}:=\left\{(x_{1},x_{2})\in\mathbf{R}^{2}\big{|}
   x_{1}\in[-R+h,R],x_{2}\in[0,1]\right\}$,
   $\widetilde{\Omega}_{R,h}$ $:=\left\{(x_{1},x_{2})\in\mathbf{R}^{2}\big{|}x_{1}\in[-R,R-h],x_{2}\in[0,1]\right\}$. Moreover,
   for $0<\varepsilon<h$,
   we introduce these
   operators $\Delta_{\varepsilon}$ and $\Delta^{\varepsilon}$ as follows.
   \begin{align}\label{chafensuanzi0127}
        \Delta_{\varepsilon}\mathbf{u}=\frac{\mathbf{u}(x_{1}+\varepsilon,x_{2})
        -\mathbf{u}(x_{1},x_{2})}{\varepsilon},~~
        \Delta^{\varepsilon}\mathbf{u}=\frac{\mathbf{u}(x_{1},x_{2}+\varepsilon)
        -\mathbf{u}(x_{1},x_{2})}{\varepsilon}.
   \end{align}
   When $\mathbf{v}\in \widetilde{\mathbf{G}}^{1,q'}\left(\overline{\Omega}_{R,h}\right)$, one can see that
   $\Delta_{-\varepsilon}\mathbf{v}\in \widetilde{\mathbf{G}}^{1,q'}(\Omega_{R})$.
   Thus, from \eqref{chentianyi0116},
   one obtains
   \begin{align*}
          \mu\int_{\Omega_{R}}\nabla\widehat{\mathbf{u}}\cdot\nabla\Delta_{-\varepsilon}\mathbf{v}dx=
          \int_{\Omega_{R}}\mathbf{f}\cdot
          \Delta_{-\varepsilon}\mathbf{v}dx.
 \end{align*}
 Therefore, conducting a computation gives
     \begin{align}\label{chentianyi01161}
          \int_{\Omega_{R}}\nabla\Delta_{\varepsilon}\widehat{\mathbf{u}}\cdot\nabla\mathbf{v}dx
          =-\frac{1}{\mu}\int_{\Omega_{R}}\mathbf{f}\cdot
          \Delta_{-\varepsilon}\mathbf{v}dx.
 \end{align}
 And there exists a constant $C>0$, such that
 \begin{align}\label{chentianyi01162}
 -\frac{1}{\mu}\int_{\Omega_{R}}\mathbf{f}\cdot
          \Delta_{-\varepsilon}\mathbf{v}dx\leq
          C\big{\|}\mathbf{f}\big{\|}_{L^{q}(\Omega_{R})}\big{\|}\nabla\mathbf{v}\big{\|}_{L^{q'}\left(\overline{\Omega}_{R,h}\right)}.
 \end{align}
 Define a functional $\mathcal{B}_{\varepsilon}$ on $\widetilde{G}^{1,q'}\left(\overline{\Omega}_{R,h}\right)$ as follows,
   \begin{align}\label{5170217}
   \mathcal{B}_{\varepsilon}(\mathbf{v})= \int_{\Omega_{R}}\nabla\Delta_{\varepsilon}\widehat{\mathbf{u}}\cdot\nabla\mathbf{v}+
\Delta_{\varepsilon}\widehat{\mathbf{u}}\cdot\mathbf{v}dx,~\forall~\mathbf{v}\in \widetilde{G}^{1,q'}\left(\overline{\Omega}_{R,h}\right).
   \end{align}
   Since $\mathbf{v}\in \left[W_{0}^{1,q'}(\overline{\Omega}_{R,h})\right]^2$, then
   $\big{\|}\mathbf{v}\big{\|}_{W_{0}^{1,q'}\left(\overline{\Omega}_{R,h}\right)}$ and
   $\big{\|}\nabla\mathbf{v}\big{\|}_{L^{q'}\left(\overline{\Omega}_{R,h}\right)}$ are equivalent.
   From \eqref{chentianyi01162} and
$\big{\|}\Delta_{\varepsilon}
\widehat{\mathbf{u}}\big{\|}_{L^{q}(\Omega_{R})}$
$\leq C\big{\|}\mathbf{u}\big{\|}_{W^{1,q}(\Omega_{R})}$, one can conclude that $\mathcal{B}_{\varepsilon}$ is bounded,
   and $\big{\|}\mathcal{B}_{\varepsilon}\big{\|}\leq C\left(\big{\|}\mathbf{f}\big{\|}_{L^{q}(\Omega_{R})}+\big{\|}\mathbf{u}\big{\|}_{W^{1,q}(\Omega_{R})}\right)$.

   Now we consider $\left[\widetilde{\mathbf{G}}^{1,q'}\left(\overline{\Omega}_{R,h}\right)\right]'$.
   Define a map $\mathbf{P}$ as follows,
   \begin{align*}
   \mathbf{P}(\mathbf{v})=(v_{1},v_{2},\partial_{1}v_{1},\partial_{2}v_{1},\partial_{1}v_{2}),
   ~\forall~\mathbf{v}\in \widetilde{\mathbf{G}}^{1,q'}\left(\overline{\Omega}_{R,h}\right).
   \end{align*}
   It is easy to discover that $\mathbf{P}:\widetilde{\mathbf{G}}^{1,q'}\left(\overline{\Omega}_{R,h}\right)
   \rightarrow \mathbf{W}\subset \left[\mathbf{L}^{q'}\left(\overline{\Omega}_{R,h}\right)\right]^{5}$ is isometric isomorphism, where
   $\mathbf{W}$ is the range of $\widetilde{\mathbf{G}}^{1,q'}\left(\overline{\Omega}_{R,h}\right)$
   under the map $\mathbf{P}$.
   And, one can show that $\mathbf{W}$ is a closed subspace of $\left[\mathbf{L}^{q'}\left(\overline{\Omega}_{R,h}\right)\right]^{5}$.
   By a similar argument as the proof of Proposition \ref{Rieszbiaoshiyinli0110}, one can verify that
   $\forall $ $\mathcal{L}\in \mathbf{W}'$, there exists a unique $\mathbf{r}=(r_{1},r_{2},r_{3},r_{4},r_{5})
   \in \left[\mathbf{L}^{q'}\left(\overline{\Omega}_{R,h}\right)\right]^{5}$,
   such that
   $\mathcal{L}(\mathbf{v})=\int_{\overline{\Omega}_{R,h}}\mathbf{r}\cdot \mathbf{P}(\mathbf{v}) dx$.
   From \eqref{5170217}, one can define a linear and bounded functional on $\mathbf{W}$ in form
   $\mathcal{B}_{\varepsilon}(\mathbf{v})=\mathcal{L}(\mathbf{P}(\mathbf{v}))$, $\forall$ $\mathbf{v}\in \widetilde{\mathbf{G}}^{1,q'}\left(\overline{\Omega}_{R,h}\right)$.
   Since $\mathbf{P}$ is isometric isomorphism, then
   $\big{\|}\mathcal{B}_{\varepsilon}\big{\|}=\big{\|}\mathcal{L}\big{\|}=\big{\|}\mathbf{r}\big{\|}_{L_{5}^{q}\left(\overline{\Omega}_{R,h}\right)}
   \leq C\left(\big{\|}\mathbf{f}\big{\|}_{L^{q}(\Omega_{R})}+\big{\|}\mathbf{u}\big{\|}_{W^{1,q}(\Omega_{R})}\right)$. Then $\Delta_{\varepsilon}\widehat{u}_{1}=r_{1}$, $\Delta_{\varepsilon}\widehat{u}_{2}=r_{2}$
   , $\partial_{1}\Delta_{\varepsilon}\widehat{u}_{1}=-\partial_{2}\Delta_{\varepsilon}\widehat{u}_{2}=r_{3}$,
   $\partial_{2}\Delta_{\varepsilon}\widehat{u}_{1}=r_{4}$, and $\partial_{1}\Delta_{\varepsilon}\widehat{u}_{2}=r_{5}$.
Thus, $\big{\|}\nabla\Delta_{\varepsilon}\widehat{\mathbf{u}}\big{\|}_{L^{q}\left(\overline{\Omega}_{R,h}\right)}
   \leq C \big{\|}\mathbf{r}\big{\|}_{L^{q}\left(\overline{\Omega}_{R,h}\right)}\leq C\left(\big{\|}\mathbf{f}\big{\|}_{L^{q}(\Omega_{R})}+\big{\|}\mathbf{u}\big{\|}_{W^{1,q}(\Omega_{R})}\right)$.
   Let $\varepsilon\rightarrow 0$, then one can obtain
   $\partial_{11}u_{1}$, $\partial_{12}u_{1}$, $\partial_{11}u_{2}$
   and $\partial_{12}u_{2}$$\in$$L^{q}\left(\overline{\Omega}_{R,h}\right).$
   Moreover, $\big{\|}\partial_{11}u_{1}\big{\|}_{L^q\left(\overline{\Omega}_{R,h}\right)}$$+$$\big{\|}\partial_{12}u_{1}\big{\|}_{L^q\left(\overline{\Omega}_{R,h}\right)}$ $+$$\big{\|}\partial_{11}u_{2}\big{\|}_{L^q\left(\overline{\Omega}_{R,h}\right)}$$+$$\big{\|}\partial_{12}u_{2}\big{\|}_{L^q\left(\overline{\Omega}_{R,h}\right)}$ $\leq C\left(\big{\|}\mathbf{f}\big{\|}_{L^{q}(\Omega_{R})}+\big{\|}\mathbf{u}\big{\|}_{W^{1,q}(\Omega_{R})}\right)$.

   When $\mathbf{v}\in \widetilde{\mathbf{G}}^{1,q'}\left(\widetilde{\Omega}_{R,h}\right)$, one can verify that
   $\Delta_{\varepsilon}\mathbf{v}\in \widetilde{\mathbf{G}}^{1,q'}(\Omega_{R})$. Thus, a similar argument gives
   that $\partial_{11}u_{1}$, $\partial_{12}u_{1}$, $\partial_{11}u_{2}$
   and $\partial_{12}u_{2}$$\in$$L^{q}\left(\widetilde{\Omega}_{R,h}\right)$.
   Moreover, $\big{\|}\partial_{11}u_{1}\big{\|}_{L^q\left(\widetilde{\Omega}_{R,h}\right)}$$+$$\big{\|}\partial_{12}u_{1}\big{\|}_{L^q\left(\widetilde{\Omega}_{R,h}\right)}$ $+$$\big{\|}\partial_{11}u_{2}\big{\|}_{L^q\left(\widetilde{\Omega}_{R,h}\right)}$$+$
$\big{\|}\partial_{12}u_{2}\big{\|}_{L^q\left(\widetilde{\Omega}_{R,h}\right)}$ $\leq C\left(\big{\|}\mathbf{f}\big{\|}_{L^{q}(\Omega_{R})}+\big{\|}\mathbf{u}\big{\|}_{W^{1,q}(\Omega_{R})}\right)$.

   Thus, one can conclude that $\partial_{11}u_{1}$, $\partial_{12}u_{1}$, $\partial_{11}u_{2}$
   and $\partial_{12}u_{2}$$\in$$L^{q}(\Omega_{R})$.
   Moreover, $\big{\|}\partial_{11}u_{1}\big{\|}_{L^q(\Omega_{R})}$$
   +$
   $\big{\|}\partial_{12}u_{1}\big{\|}_{L^q(\Omega_{R})}$ $+$$\big{\|}\partial_{11}u_{2}\big{\|}_{L^q(\Omega_{R})}$$+$
$\big{\|}\partial_{12}u_{2}\big{\|}_{L^q(\Omega_{R})}$ $\leq C\left(\big{\|}\mathbf{f}\big{\|}_{L^{q}(\Omega_{R})}+\big{\|}\mathbf{u}\big{\|}_{W^{1,q}(\Omega_{R})}\right).$

   Similarly, in the $x_{2}$ direction, the operator $\Delta^{\varepsilon}$ is used to improve the
   regularity. One can conclude that  $\partial_{12}u_{1}$, $\partial_{22}u_{1}$, $\partial_{12}u_{2}$
   and $\partial_{22}u_{2}$$\in$$L^{q}(\Omega_{R})$.
   Moreover, $\big{\|}\partial_{12}u_{1}\big{\|}_{L^q(\Omega_{R})}$$
   +$   $\big{\|}\partial_{22}u_{1}\big{\|}_{L^q(\Omega_{R})}$ $+$$\big{\|}\partial_{12}u_{2}\big{\|}_{L^q(\Omega_{R})}$
   $+$
$\big{\|}\partial_{22}u_{2}\big{\|}_{L^q(\Omega_{R})}$ $\leq C\left(\big{\|}\mathbf{f}\big{\|}_{L^{q}(\Omega_{R})}+\big{\|}\mathbf{u}\big{\|}_{W^{1,q}(\Omega_{R})}\right)$.

   Thus, one can get that
   \begin{align}\label{yongdejiushini0127}
    \big{\|}\nabla^2\mathbf{u}\big{\|}_{L^{q}(\Omega_{R})}\leq C\left(\big{\|}\mathbf{f}\big{\|}_{L^{q}(\Omega_{R})}+\big{\|}\mathbf{u}\big{\|}_{W^{1,q}(\Omega_{R})}\right),
   \end{align}
   where $C$ is independent of $R$.

   Since $\mathbf{u}\in \left[W^{2,q}(\Omega_{R})\right]^2$, then $\nabla p=\mathbf{f}+\mu\Delta\mathbf{u}$. Thus,
   one can easily obtain that
   \begin{align}\label{moranhuishou0127}
   \big{\|}\nabla p\big{\|}_{L^{q}(\Omega_{R})}\leq C\left(\big{\|}\mathbf{f}\big{\|}_{L^{q}(\Omega_{R})}+\big{\|}\mathbf{u}\big{\|}_{W^{1,q}(\Omega_{R})}\right).
   \end{align}

   Thus, in view of \eqref{yongdejiushini0127} and \eqref{moranhuishou0127}, one can conclude that
   \eqref{xiangyaodejieguo} is valid.

   Next, we show that $\mathbf{u}\in\mathbf{V}^{1,q}(\Omega_{R})$. Multiply $\eqref{tuiguanjie}_{1}$ by
   $\mathbf{v}\in \mathbf{G}^{1,q'}(\Omega_{R})$ and integrate over $\Omega_{R}$, then we get
   \begin{align*}
    \mu\int_{\Omega_{R}}\nabla\mathbf{u}\cdot\nabla\mathbf{v}dx-
    \mu\int_{-R}^{R}\left[\partial_{2}u_{1}(x_{1},1)v_{1}(x_{1},1)-\partial_{2}u_{1}(x_{1},0)v_{1}(x_{1},0)\right]dx_{1}=
    \int_{\Omega_{R}}\mathbf{f}\cdot\mathbf{v}dx.
   \end{align*}
From the definition of $q-$generalized solution, one can conclude that
   $\forall \mathbf{v}\in \mathbf{G}^{1,q'}(\Omega_{R})$,
   \begin{align*}
    \int_{-R}^{R}[\mu\partial_{2}u_{1}(x_{1},1)v_{1}(x_{1},1)-\mu\partial_{2}u_{1}(x_{1},0)v_{1}(x_{1},0)
    -k_{1}u_{1}(x_{1},1)v_{1}(x_{1},1)-k_{0}u_{1}(x_{1},0)v_{1}(x_{1},0)] dx_{1}=0.
   \end{align*}
   Thus, $\partial_{2}u_{1}(x_{1},1)=\frac{k_{1}}{\mu}u_{1}(x_{1},1)$
   and $\partial_{2}u_{1}(x_{1},0)=-\frac{k_{0}}{\mu}u_{1}(x_{1},0)$ indicating that
   $\mathbf{u}\in\mathbf{V}^{1,q}(\Omega_{R})$.
\end{proof}

\subsection{Appendix C}\label{appendixb1226}
\begin{proof}\textbf{Proof of Lemma \ref{tezhenzhidewenti}}.

    (1) In order to obtain the weak solution of problem \eqref{eigen1}, we first consider the following variational
    problem:
    \begin{align*}
        \inf\limits_{\mathbf{u}\in \mathbf{E}_{1}}F(\mathbf{u}),
    \end{align*}
    where $F(\mathbf{u})=\mu \int_{\Omega_{R}}|\nabla\mathbf{u}|^2dx-\int_{-R}^{R}
    \left[k_{1}|u_{1}(x_{1},1)|^2+k_{0}|u_{1}(x_{1},0)|^2\right]dx_{1}$ and $
    \mathbf{E}_{1}=\left\{\mathbf{u}\in \mathbf{G}^{1,2}(\Omega_{R})\big{|}\big{\|}\mathbf{u}\big{\|}_{L^2(\Omega_{R})}=1\right\}.$

    From conditions $k_{1}\leq 0$ and $k_{0}\leq 0$, then $F(\mathbf{u})\geq 0$.
    Thus, there exists a sequence $\{\mathbf{u}_{n}\}\subset E_{1}$, such that $F(\mathbf{u}_{n})\rightarrow
    \inf\limits_{\mathbf{u}\in E_{1}}F(\mathbf{u}).$ Here, we assume $F(\mathbf{u}_{n})\leq
    \inf\limits_{\mathbf{u}\in E_{1}}F(\mathbf{u})+1$ implying that $\{\mathbf{u}_{n}\}
    \subset\left[H^{1}(\Omega_{R})\right]^2$ is bounded.
    Due to the reflexivity of $\left[H^{1}(\Omega_{R})\right]^2$,
    there exists $\mathbf{e}^{1}\in \left[H^1(\Omega_{R})\right]^2$ such that
    \begin{align*}
        \begin{aligned}
        &\mathbf{u}_{n} \rightarrow \mathbf{e}^1~\text{weakly~in}~\left[H^{1}(\Omega_{R})\right]^2,
        \\
        &\mathbf{u}_{n} \rightarrow \mathbf{e}^1~\text{strongly~in}~\left[L^{2}(\Omega_{R})\right]^2.
        \end{aligned}
    \end{align*}

    Since $G(\mathbf{u}):=\mu\int_{\Omega_{R}}|\nabla \mathbf{u}|^2dx$ is convex, thus
    $F(\mathbf{u})$ is weakly lower continuous. Then,
    \begin{align*}
        F\left(\mathbf{e}^{1}\right)\leq\lim\limits_{n\rightarrow\infty}F(\mathbf{u_{n}})=
        \inf\limits_{\mathbf{u}\in \mathbf{E}_{1}}F(\mathbf{u}).
    \end{align*}
    Next, we show that $\mathbf{e}^{1}\in \mathbf{E}_{1}.$ First, $ \mathbf{e}^{1}\cdot
    \mathbf{n}|_{\Gamma_{1R}\cup\Gamma_{2R}}=0,
    ~ \mathbf{e}^{1}|_{\Gamma_{3 R}\cup\Gamma_{4R}}=\mathbf{0}$ and
    $\big{\|}\mathbf{e}^{1}\big{\|}_{L^2(\Omega_{R})}=1$.
    From the definition of weak derivative, one has $\nabla\cdot \mathbf{e}^{1}=0.$
    Then $\mathbf{e}^{1}\in \mathbf{E}_{1}$.
    Thus, $F\left(\mathbf{e}^{1}\right)=\inf\limits_{u\in \mathbf{E}^{1}}F(\mathbf{u}).$

    Let $\lambda_{1}=F\left(\mathbf{e}^{1}\right)$, i.e., $\lambda_{1}=\mu\int_{\Omega_{R}}
    |\nabla \mathbf{e}^{1}|^2 dx-\int_{-R}^{R}\left[k_{1}|\mathbf{e}_{1}^{1}(x_{1},1)|^2+
    k_{0}\left|\mathbf{e}^{1}_{1}(x_{1},0)\right|^2\right]dx_{1}>0.$
    Assume that $I(t,r)=\int_{\Omega_{R}}|\mathbf{e}^{1}+t\mathbf{u}+r\mathbf{e}^{1}|^2dx$,
    where $\mathbf{u}\in \left[H^{1}(\Omega_{R})\right]^2,~ \nabla\cdot \mathbf{u}=0$ and $t,r\in\mathbf{R}$. Then,
    \begin{align*}
        \begin{aligned}
            &I(0,0)=\big{\|}\mathbf{e}^{1}\big{\|}_{L^2(\Omega_{R})}^{2}=1,~I\in C^{\infty},
            \\
            &\partial_{t}I(0,0)=2\int_{\Omega_{R}}\mathbf{e}^{1}\cdot \mathbf{u} dx,
            ~\partial_{r}I(0,0)=2\big{\|}\mathbf{e}^{1}\big{\|}_{L^2(\Omega_{R})}^{2}=2\neq 0.
        \end{aligned}
    \end{align*}
By the theorem of implicit function,
in the neighbourhood of $(0,0)$,
there exists a smooth curve $r=r(t)$ satisfying $r(0)=0$ and
$r'(0)=-\int_{\Omega_{R}}\mathbf{e}^{1}\cdot \mathbf{u}dx $. Let $J(t)
=F\left(\mathbf{e}^{1}+t \mathbf{u}+r(t)\mathbf{e}^1\right)$.
Since $\mathbf{e}^{1}$ is extreme point, then
\begin{align}\label{dingding}
    \begin{aligned}
        0=J'(0)=2\mu \int_{\Omega_{R}}&\nabla\mathbf{e}^1\cdot\nabla\mathbf{u}dx-2
        \int_{-R}^{R}\bigg{[}k_{1}\mathbf{e}^{1}_{1}(x_{1},1)u_{1}(x_{1},1)
        \\
        &+k_{0}\mathbf{e}^{1}_{1}(x_{1},0)u_{1}(x_{1},0)\bigg{]}dx_{1}-2\lambda_{1}\int_{\Omega_{R}}
        \mathbf{u}\cdot\mathbf{e}^{1}dx.
    \end{aligned}
\end{align}
If take $\mathbf{u}\in \mathbf{G}^{1,2}(\Omega_{R})$, then we have
\begin{align*}
    \mu\int_{\Omega_{R}}\nabla \mathbf{e}^{1}\cdot\nabla \mathbf{u}dx-
    \int_{-R}^{R}\left[k_{1}\mathbf{e}^{1}_{1}(x_{1},1)u_{1}(x_{1},1)
        +k_{0}\mathbf{e}^{1}_{1}(x_{1},0)u_{1}(x_{1},0)\right]dx_{1}=\lambda_{1}\int_{\Omega_{R}}
    \mathbf{e}^{1}\cdot\mathbf{u}dx,
\end{align*}
By the Lemma \ref{dengjia}, there exists a $p^{1}\in L^{2}(\Omega_{R})$ such that
when $\lambda=\lambda_{1}$, $\left(\mathbf{e}^{1},
p^1\right)$ is the weak solution of problem \eqref{eigen1}.

One can apply the method in Lemma \ref{tishengzhengzexing} to improve the regularity of weak solution, and further
obtain
$\mathbf{e}^1\in [H^{\infty}(\Omega_{R})]^{2}$ and $p^1\in H^{\infty}(\Omega_{R})$.
In \eqref{dingding}, let $\mathbf{u}\in \mathbf{E}_{1}$. It follows that
\begin{align*}
    \begin{aligned}
    &-\mu\int_{\Omega_{R}}\Delta\mathbf{e}^1\cdot \mathbf{u}dx+
    \mu\int_{-R}^{R}\left[\partial_2\mathbf{e}^{1}_{1}(x_{1},1)u_{1}(x_{1},1)-
    \partial_{2}\mathbf{e}^{1}_{1}(x_{1},0)u_{1}(x_{1},0)\right]dx_{1}
    \\
    &-\int_{-R}^{R}\left[k_{1}\mathbf{e}^{1}_{1}(x_{1},1)u_{1}(x_{1},1)
    +k_{0}\mathbf{e}^{1}_{1}(x_{1},0)u_{1}(x_{1},0)\right]dx_{1}-\lambda_{1}\int_{\Omega_{R}}
    \mathbf{u}\cdot\mathbf{e}^{1}dx=0.
    \end{aligned}
\end{align*}
Then, one has
\begin{align*}
    \begin{aligned}
    &\mu\int_{-R}^{R}\left[\partial_{x_{2}}\mathbf{e}^{1}(x_{1},1)u_{1}(x_{1},1)-
    \partial_{2} \mathbf{e}^{1}_{1}(x_{1},0)u_{1}(x_{1},0)\right]dx_{1}
    \\
    &-\int_{-R}^{R}\left[k_{1}\mathbf{e}^{1}_{1}(x_{1},1)u_{1}(x_{1},1)
    +k_{0}\mathbf{e}^{1}_{1}(x_{1},0)u_{1}(x_{1},0)\right]dx_{1}=0.
    \end{aligned}
\end{align*}
Therefore, $\partial_{2}\mathbf{e}^{1}_{1}(x_{1},1)=\frac{k_{1}}{\mu}\mathbf{e}^{1}_{1}(x_{1},1)$
and $\partial_{2}\mathbf{e}^{1}_{1}(x_{1},0)=-\frac{k_{0}}{\mu}\mathbf{e}^{1}_{1}(x_{1},0).$

In order to obtain a sequence of solution, we continue to consider the following variational problem,
\begin{align*}
    \inf_{\mathbf{u}\in \mathbf{E}_{2}}F(\mathbf{u}),
\end{align*}
where $\mathbf{E}_{2}=\left\{\mathbf{u}\in \mathbf{E}_{1}\big{|}\int_{\Omega_{R}}\mathbf{u}\cdot \mathbf{e}^1dx=0\right\}$.

Similarly, one can obtain
$\left(\mathbf{e}^2,\nabla p^2, \lambda_{2}\right)$, where
$\lambda_{2}=F\left(\mathbf{e}^2\right)=\inf\limits_{\mathbf{u}\in \mathbf{E}_{2}}F(\mathbf{u})$,
$-\mu \Delta\mathbf{e}^2+\nabla p^2=\lambda_{2}\mathbf{e}^2$, $\mathbf{e}^2\in \mathbf{E}_{2}\cap \left[H^{\infty}\right]^2$,
 $\partial_{2}\mathbf{e}^{2}_{1}(x_{1},1)=\frac{k_{1}}{\mu}\mathbf{e}^{2}_{1}(x_{1},1)$
 and $\partial_{2}\mathbf{e}^{2}_{1}(x_{1},0)=-\frac{k_{0}}{\mu}\mathbf{e}^{2}_{1}(x_{1},0).$

Thus, from the following variational problem,
\begin{align*}
    \inf_{\mathbf{u}\in \mathbf{E}_{n}}F(\mathbf{u}),
\end{align*}
one can get the solution $\left(\mathbf{e}^n,\nabla p^n,\lambda_{n}\right)$ of problem \eqref{eigen1}, where
$\lambda_{n}=F\left(\mathbf{e}^n\right)=\inf\limits_{\mathbf{u}\in \mathbf{E}_{n}}F(\mathbf{u})$,
 $\mathbf{E}_{n}=\left\{\mathbf{u}\in \mathbf{E}_{n-1}\big{|}\int_{\Omega_{R}}\mathbf{u}\cdot \mathbf{e}^{n-1}dx=0\right\}$.
Hence, the first conclusion is valid.

(2) Since $\mathbf{E}_{n}\subset \mathbf{E}_{n-1}$, it is easy to find that $\lambda_{n-1}\leq \lambda_{n}$.

Assume that there exists a constant $C>0$ such that $\lambda_{n}\leq C$,
then $\lambda_{n}=F(\mathbf{e}^{n})=\mu\int_{\Omega_{R}}|\nabla\mathbf{e}^{n}|^2dx
 -\int_{-R}^{R}\left[k_{1}|\mathbf{e}^n_{1}(x_{1},1)|^2+k_{0}|\mathbf{e}^{n}_{1}(x_{1},0)|^2\right]dx_{1}\leq C,$
 there exists $\mathbf{e}^{0}$ such that
 \begin{align}\label{maodundeqiyuan}
     \begin{aligned}
     &\mathbf{e}^{n} \rightarrow \mathbf{e}^0~\text{weakly~in}~\left[H^{1}(\Omega_{R})\right]^2,
     \\
     &\mathbf{e}^{n} \rightarrow \mathbf{e}^0~\text{strongly~in}~\left[L^{2}(\Omega_{R})\right]^2.
     \end{aligned}
 \end{align}
 Let $m>n$. Then from $\int_{\Omega_{R}}\mathbf{e}^{m}\mathbf{e}^n dx=\delta_{mn}$ and
 $\big{\|}\mathbf{e}^{n}\big{\|}_{L^2(\Omega_{R})}=1$, one has that
\begin{align*}
 \int_{\Omega_{R}}\big{|}\mathbf{e}^m-\mathbf{e}^n\big{|}^2dx=2.
\end{align*}
Passing $m\rightarrow\infty$, then we have
\begin{align*}
 \int_{\Omega_{R}}\left|\mathbf{e}^0-\mathbf{e}^n\right|^2dx=2\neq 0,
\end{align*}
which contradicts \eqref{maodundeqiyuan}. Thus, when
$n\rightarrow\infty$, one gets $\lambda_{n}\rightarrow+\infty$.

(3) Verify that $\left\{\mathbf{e}^n\right\}$ is an orthogonal basis of $\mathbf{G}^{1,2}(\Omega_{R})$.

Firstly, we prove that $\{\mathbf{e}^n\}$ is orthogonal in $\mathbf{G}^{1,2}(\Omega_{R})$. In fact,
\begin{align*}
    \begin{aligned}
        \left((\mathbf{e}^n,\mathbf{e}^m)\right)&=
        \mu\int_{\Omega_{R}}\nabla\mathbf{e}^n \nabla\mathbf{e}^m dx
        -\int_{-R}^{R}\bigg[k_{1}\mathbf{e}^{n}_{1}(x_{1},1)\mathbf{e}^{m}_{1}(x_{1},1)+
        k_{0}\mathbf{e}^{n}_{1}(x_{1},0)\mathbf{e}^{m}_{1}(x_{1},0)\bigg{]}dx_{1}
        \\
        &=-\mu\int_{\Omega_{R}}\Delta\mathbf{e}^{n}\cdot \mathbf{e}^m dx
        \\
        &=\lambda_{n}\int_{\Omega_{R}}\mathbf{e}^n\cdot \mathbf{e}^mdx=\lambda_{n}\delta_{mn}.
    \end{aligned}
\end{align*}
Next, one can verify that $\forall \mathbf{u}\in \mathbf{G}^{1,2}(\Omega_{R})$
, $\big{\|}\mathbf{u}-\mathbf{S}^m\big{\|}_{L^2(\Omega_{R})}
\rightarrow 0$ as $m\rightarrow\infty$,
where $\mathbf{S}^{m}=\sum\limits_{k=1}^{m}c_{k} \mathbf{e}^{k}$,
$c_{k}=\frac{\left((\mathbf{u},\mathbf{e}^{k})\right)}{\lambda_{k}}=\left(\mathbf{u},\mathbf{e}^{k}\right)$. Indeed, a direct computation gives that
\begin{align*}
    \begin{aligned}
        \big{\|}\mathbf{u}-\mathbf{S}^m\big{\|}_{\mathbf{G}^{1,2}(\Omega_{R})}^{2}=
        \big{\|}\mathbf{u}\big{\|}_{\mathbf{G}^{1,2}(\Omega_{R})}^2-\sum\limits_{k=1}^{m}\lambda_{k}|c_{k}|^2\geq 0.
    \end{aligned}
\end{align*}
Thus, one has
\begin{align}\label{shoulian}
\sum\limits_{k=1}^{m}\lambda_{k}|c_{k}|^2\leq \big{\|}\mathbf{u}\big{\|}_{\mathbf{G}^{1,2}(\Omega_{R})}^2,
\end{align}
which implies that $\left\{\sum\limits_{k=1}^{m}\lambda_{k}|c_{k}|^2\right\}$ is convergent, since
$\left\{\sum\limits_{k=1}^{m}\lambda_{k}|c_{k}|^2\right\}$ is increasing with respect to $m$.
Moreover, for $n=1,2,\cdots,m$, we discover that
\begin{align*}
    \lambda_{n}\left(\mathbf{u}-\mathbf{S}^m,\mathbf{e}^n\right)_{L^2}=\left((\mathbf{u}-\mathbf{S}^m,\mathbf{e}^n)\right)=0,
\end{align*}
indicating that $\frac{\mathbf{u}-\mathbf{S}^m}{\|\mathbf{u}-\mathbf{S}^m\|_{L^2}} \in \mathbf{E}_{m+1}$. Then
\begin{align*}
    \big{\|}\mathbf{u}\big{\|}_{\mathbf{G}^{1,2}(\Omega_{R})}^2
    -\sum\limits_{k=1}^{m}\lambda_{k}|c_{k}|^2 =
   \big{\|}\mathbf{u}-\mathbf{S}^m\big{\|}_{\mathbf{G}^{1,2}(\Omega_{R})}^2=F(\mathbf{u}-\mathbf{S}^m)\geq
    \lambda_{m+1}\big{\|}\mathbf{u}-\mathbf{S}^m\big{\|}_{L^2(\Omega_{R})}^2.
\end{align*}
Thus, one has
\begin{align*}
    \big{\|}\mathbf{u}-\mathbf{S}^m\big{\|}_{L^2(\Omega_{R})}^2\leq \frac{1}{\lambda_{m+1}}
    \big{\|}\mathbf{u}\big{\|}_{\mathbf{G}^{1,2}(\Omega_{R})}^2,
\end{align*}
which implies that
\begin{align}\label{l2shoulian}
    \big{\|}\mathbf{u}-\mathbf{S}^m\big{\|}_{L^2(\Omega_{R})}^2\rightarrow 0,~~\text{as}~m\rightarrow\infty.
\end{align}
Finally, we show that $\big{\|}\mathbf{u}-\mathbf{S}^m\big{\|}_{G^{1,2}(\Omega_{R})}^2\rightarrow 0$ as $m\rightarrow\infty$. In fact,
take $m>n$, then from \eqref{shoulian}, one can obtain that
\begin{align*}
\big{\|}\mathbf{S}^m-\mathbf{S}^n\big{\|}_{G^{1,2}(\Omega_{R})}^2=
\sum\limits_{k=n+1}^{m}\lambda_{k}|c_{k}|^{2}\rightarrow 0,~~\text{as}~n\rightarrow\infty.
\end{align*}
Thus, $\{\mathbf{S}^n\}$ is a Cauchy sequence in $\mathbf{G}^{1,2}(\Omega_{R})$, which indicates that there exists
$\mathbf{S}\in \mathbf{G}^{1,2}(\Omega_{R})$ such that
$\mathbf{S}^n\rightarrow \mathbf{S}$ in $\mathbf{G}^{1,2}(\Omega_{R})$. Furthermore, from \eqref{l2shoulian}, one can obtain
$\mathbf{S}=\mathbf{u}$ in $\left[L^2(\Omega_{R})\right]^2$. Then we can conclude that
$\mathbf{S}=\mathbf{u}$ in $\mathbf{G}^{1,2}(\Omega_{R})$.
\end{proof}
\section*{Acknowledgments}
This research is supported by the Natural Science Foundation
of China (No. 11801108), the Guangdong Basic and Applied Basic Research Foundation (No.
2023A1515030107), the Science and Technology Planning Project of Guangzhou
City (No. 202201010111), the Guangzhou City-College-Eenterprise Joint Funding Project (No. SL2023A03J00375), the Tertiary Education Scientific Research Project of Guangzhou Municipal
Education Bureau (No. 202235103), and the Guangzhou Education Scientific Research Project (No. 202214066), the Natural Science Foundation of Sichuan Province Grant (No. 2022NSFSC1818).

\end{document}